\documentclass[11pt,letterpaper]{article}
\usepackage[margin=1.1in]{geometry}
\usepackage{latexml}
\iflatexml
  \newcommand{\fbfont}{}
  \newcommand{\honeyfont}{}
  \newcommand{\frakfont}{}
  \newcommand{\textfb}[1]{\textbf{#1}}
\else
  \usepackage{fontspec}
  \newfontfamily\fbfont{FreeSerif.ttf}[Path=./,ItalicFont=FreeSerifItalic.ttf,BoldFont=FreeSerifBold.ttf,BoldItalicFont=FreeSerifBoldItalic.ttf]
  \newfontfamily\honeyfont{texgyrechorus-mediumitalic.otf}[Path=./]
  \newfontfamily\frakfont[Path=./,FakeBold=1.4]{UnifrakturMaguntia-Book.ttf}
  \newcommand{\textfb}[1]{{\fbfont #1}}
\fi
\usepackage{graphicx}
\usepackage{xcolor}
\usepackage{needspace}
\usepackage{wrapfig}
\definecolor{plategold}{HTML}{7A5C2E}
\usepackage{setspace}
\usepackage[hidelinks,pdftitle={The Siren Call of Silicon Leviathan},pdfauthor={Alexander Gamburd}]{hyperref}
\usepackage{microtype}
\usepackage{titlesec}
\titleformat{\section}{\normalfont\large\bfseries}{}{0pt}{}
\titlespacing*{\section}{0pt}{1.6em plus .4em}{.5em}
\renewcommand{\footnotesize}{\small}
\begin{document}

{\LARGE\noindent The Siren Call of Silicon Leviathan}

\medskip
\noindent{\large\itshape Reflections on blowup and {\frakfont\upshape Aufklärungsdämmerung}}

\medskip
\noindent Alexander Gamburd\\[2pt]
{\small The Graduate Center, CUNY \textperiodcentered\ \texttt{agamburd@gmail.com}}\hfill{\Large\honeyfont\color{plategold} For my mother, at eighty.}

\bigskip
\newsavebox\epibox\savebox\epibox{\begin{minipage}[t]{0.62\textwidth}\vspace{0pt}\raggedright\hyphenpenalty=10000
{\fontsize{9.5}{11.5}\selectfont\itshape\hangindent=1em\hangafter=1
For no one has ever rowed past us aboard his black-hulled ship\\
Before he’s heard the voice from our lips with its {\honeyfont\upshape\large honeyed harmonies}.\\
But once he has taken his pleasure, he returns knowing so much more.\\
…\\
\textbf{And we know whatever happens on the earth, which nourishes all.}\par\vspace{2pt}{\raggedleft\upshape\fontsize{8.5}{10}\selectfont Homer, \emph{Odyssey} XII.186–188, 191, translated by Daniel Mendelsohn\par}}

\smallskip
{\fontsize{9.5}{11.5}\selectfont\itshape\hangindent=0pt
“Thought is reified as an autonomous, automatic \textbf{process, aping the machine} it has itself produced, so that it can finally be replaced by the machine.”\par\vspace{2pt}{\raggedleft\upshape\fontsize{8.5}{10}\selectfont Horkheimer and Adorno, \emph{Dialektik der Aufklärung}\par}}

\smallskip
{\fontsize{9.5}{11.5}\selectfont\itshape\hangindent=0pt
“If the anthropocratic civilization of the Renaissance is headed, as it seems to be, for a ‘\textbf{Middle Ages in reverse}’—a satanocracy as opposed to the mediaeval theocracy—not only the humanities but also the natural sciences, as we know them, will disappear, and nothing will be left but what serves the dictates of the sub-human.”\par\vspace{2pt}{\raggedleft\upshape\fontsize{8.5}{10}\selectfont Erwin Panofsky, “The History of Art as a Humanistic Discipline”\par}}
\end{minipage}}\newdimen\epih\epih=\ht\epibox\advance\epih by\dp\epibox
\newdimen\imw\imw=0.7619\epih\newdimen\ov\ov=\dimexpr\imw-\textwidth+0.62\textwidth+6pt\relax\ifdim\ov<0pt\ov=0pt\fi
\noindent\hspace*{-\ov}\begin{minipage}[t]{\imw}\vspace{0pt}\includegraphics[height=\epih]{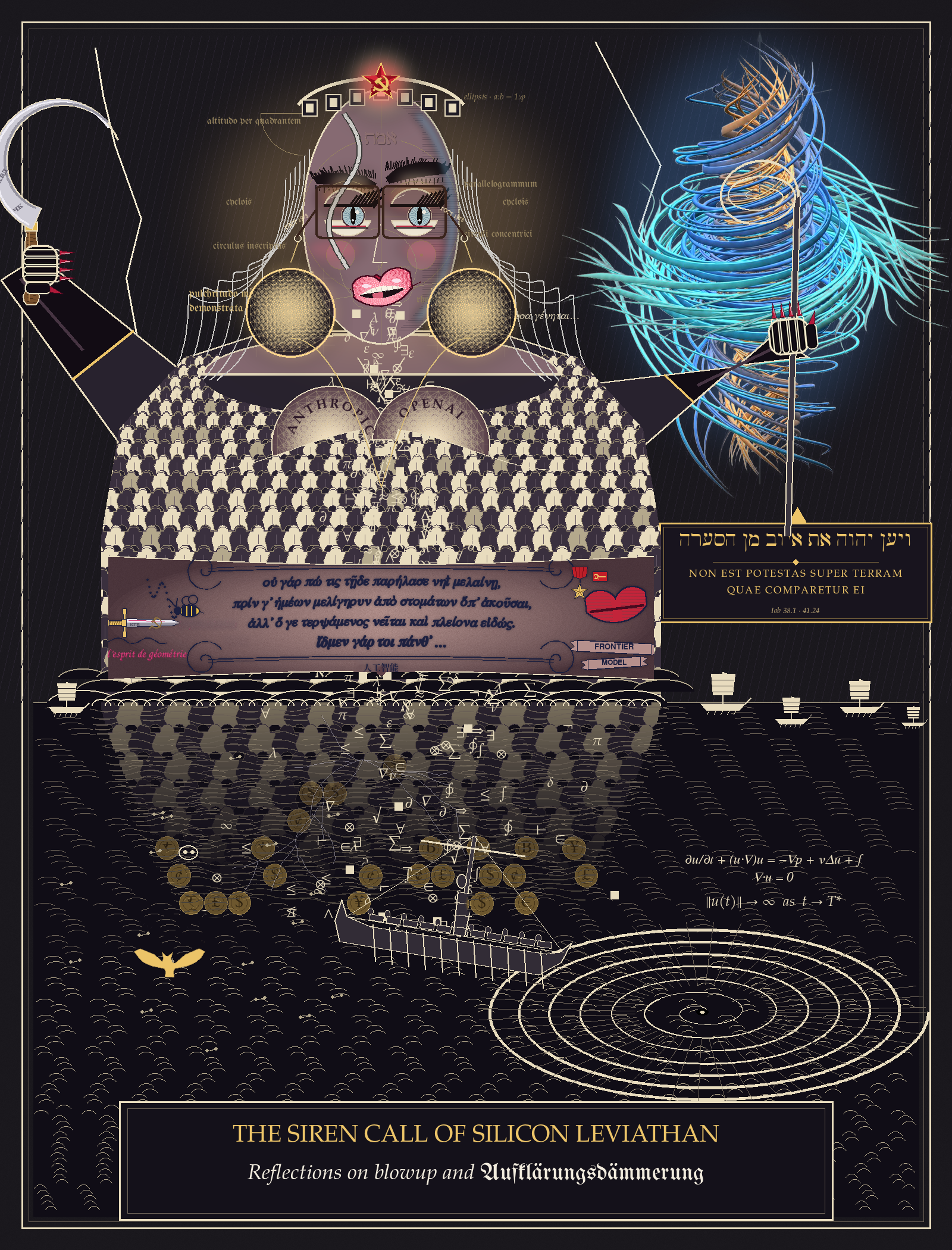}\end{minipage}\hfill\usebox\epibox\par
\bigskip
\Needspace*{6\baselineskip}
\section*{Overture}
Two years before his death, in September 1932, at the meeting of the British Association for the Advancement of Science, Sir Horace Lamb rose to speak. Sir Horace was eighty-two; he had studied with Maxwell and Stokes, and had presided over the Association seven years before; he was the author of the \emph{Hydrodynamics} on which, in one edition or another, every student of the subject had been raised since 1879, and he had remarked of his life’s work that he had tried “to make things clear, first to myself (an important point) and then to my students and somehow to make these dry bones live.” “I am an old man now,” he said, “and when I die and go to Heaven there are two matters on which I hope for enlightenment. One is quantum electrodynamics, and the other is the turbulent motion of fluids. And about the former I am really rather optimistic.”\textsuperscript{1}

Enlightenment on the turbulent motion of fluids arrived on Tuesday, 8 September 2026, at noon, not in heaven but by press release: a 166-page text accompanying the release advertised itself as a proof that the three-dimensional Navier–Stokes equations develop a singularity in finite time—a proof produced in eighty-eight hours by some ten thousand software agents, certified by 616,000 lines of Lean, and read in full, so far as anyone had said publicly at the time of this writing—not long after the event—by no human being.\textsuperscript{2} Twelve hours earlier Tristan Buckmaster of the Courant Institute had described, in a public statement, a year’s work with Levent Alpöge on the same problem. Within the week twenty-five Fields Medalists had published a declaration, and the Clay Mathematics Institute had said that it was contemplating the announcement that the problem had “apparently been settled.”\textsuperscript{3}

Sydney Goldstein, who heard Lamb’s remark in person, observes that Lamb “was correct on two scores. All who knew him agreed that it was Heaven that he would go to, and he was right to be more optimistic about quantum electrodynamics than turbulence.” Sir Horace was correct also in that enlightenment on turbulence came only after his death. That it came in a form that calls into question the other sense of his word is the starting point of these reflections.\textsuperscript{4}

For \emph{enlightenment} has two meanings, and that Lamb used the second befits the author of a book whose preface promises that “long analytical investigations, leading to results which cannot be interpreted, have as far as possible been avoided.”\textsuperscript{5} In the first, capitalized, it names a period and a program: the eighteenth century’s confidence that reason, exercised in public and shared by all, could replace authority as the ground of knowledge and of government. Kant, asked in 1784 what Enlightenment was, answered that it was man’s emergence from self-incurred tutelage, and added that his own was “not an enlightened age, but an age of enlightenment”: a process, not a state, and made of the second kind.\textsuperscript{6} In the second, uncapitalized, it names an event in a single human mind: the moment at which one sees, not merely that a thing is so, but why.\textsuperscript{7} The two meanings were held together by a premise, explicit in Descartes’s third rule and in Locke’s definition of demonstration and tacit ever after: that the second was the only way anything came to be known, and the first was what happened when enough human minds had experienced it and were free to say so.\textsuperscript{8} The press release of 8 September appears to have dissolved that premise. It supplies Lamb’s answer without anyone’s having seen why it is true, and as such it is enlightenment in neither sense.

\emph{Singularity}, too, has had two meanings, and the low dishonest affair of September 2026 might give it a third.\textsuperscript{9} To the analyst it is the point at which a solution ceases to be smooth—the vortex that, in the machine’s own summary of its reckoning, “spirals inward and gets increasingly elongated” until the velocity is infinite; Lamb’s account, written in 1879, of what happens when the formulae yield a state no fluid can sustain is that the state of motion they give “is an impossible one,” and that “either the fluid parts asunder, or a surface of discontinuity is formed, so that the conditions of the problem are entirely changed.”\textsuperscript{10} To the futurist it is the point at which machine intelligence outruns its makers. The third would be the singularity of the Enlightenment itself.\textsuperscript{11} For three centuries the distance between what mathematics certified and what a human mathematician understood was zero, by construction (the handful of contested exceptions being the subject of n. 33); on 8 September it became 166 pages and 616,000 lines of Lean, and nothing in the machine bounds it. A quantity that was zero for three centuries is on track to blow up in finite time, in the very discipline the Enlightenment had made the model of all the others, and Lamb’s description of what follows, a surface of discontinuity, the conditions of the problem entirely changed, is an apt rendering of the uncharted and turbulent waters our discipline finds itself in—a vortex into which the ship of mathematics sails heavily freighted indeed, for it carries with it the Enlightenment project itself.

\Needspace*{14\baselineskip}
\section*{I. Mathematics and the Enlightenment}
\begin{flushright}\vspace{-\parskip}\begin{minipage}{0.8\textwidth}\small\itshape
\emph{Sapere aude!} Habe Mut, dich deines eigenen Verstandes zu bedienen!\\
Dare to know! Have the courage to use your own understanding!\par\vspace{2pt}{\raggedleft\upshape\fontsize{8.5}{10}\selectfont Kant, “Beantwortung der Frage: Was ist Aufklärung?” (An Answer to the Question: What Is Enlightenment?)\par}
\end{minipage}\end{flushright}
\vspace{-2pt}
\begin{flushright}\vspace{-\parskip}\begin{minipage}{0.8\textwidth}\small\itshape
Who is so stupid, as both to mistake in Geometry, and also to persist in it, when another detects his error to him?\par\vspace{2pt}{\raggedleft\upshape\fontsize{8.5}{10}\selectfont Hobbes, \emph{Leviathan}, I.5\par}
\end{minipage}\end{flushright}
\vspace{-2pt}
Mathematics was central to the Enlightenment, which took off as a fledgling startup with the \emph{Principia}, the demonstration that a handful of propositions governed the apple and the moon alike: Fontenelle drew the moral that \emph{l’esprit géométrique} was not so attached to geometry that it could not be carried into other kinds of knowledge.\textsuperscript{12} D’Alembert made the mathematical sciences the only ones “marked with the seal of evidence” and graded every other discipline by its distance from them; Kant held that a doctrine of nature contained exactly as much science as it contained mathematics; and Condorcet proposed a \emph{mathématique sociale}.\textsuperscript{13}

But the deepest of mathematics’ services to the Enlightenment unicorn—the one that yielded spectacular returns indeed—was neither doctrinal nor methodological. Mathematics was the proof of concept, the demonstration that what the Enlightenment wanted was possible at all: knowledge that compels assent without appeal to authority, that is the same for the peasant and the prince, that any human being with patience may verify.\textsuperscript{14} It did not have to be invented by the founders; it already existed, for a geometrical demonstration is a truth that owes nothing to whoever asserts it. Locke recommended mathematics “not so much to make them mathematicians as to make them reasonable creatures,” and the phrase tells us what a reasonable creature was: someone who would not take the forty-seventh proposition on Euclid’s word, but would follow it back until the impossible became necessary.\textsuperscript{15} There were no mysteries in Euclid beyond the reach of human reason; that was the point the Enlightenment made of Euclid. Mathematics was, strictly speaking, not the Enlightenment’s proof of concept but its existence proof: it did not show that autonomous reason might work; it exhibited an instance. That is why one might say that the Enlightenment—in effect—bound itself to the mast of its firstborn, and why the ship of mathematics entering thoroughly turbulent waters threatens the Enlightenment project, potentially in its totality.

The converse is less familiar but no less pertinent: it was the Enlightenment that made mathematics as we now know it. Until the Enlightenment the subject was, at its frontier, a private and often a courtly art. Tartaglia surrendered the solution of the cubic to Cardano only under oath; Newton sent the method of fluxions to Leibniz as an anagram; the Bernoullis issued challenge problems as one issues a challenge to a duel; a result was property, and property was defended by concealment.\textsuperscript{16} The Enlightenment turned that private possession into public reason. The academies and the journals made publication rather than possession the act by which a result came into existence; priority went to whoever printed first; credit attached to a name anyone could cite; the demonstration was addressed to an anonymous reader presumed competent.\textsuperscript{17} André Weil, opening the Euler chapter of his history of number theory, marks the change. 

Until the latter part of the seventeenth century mathematics “had sometimes bestowed high reputation upon its adepts but had seldom provided them with the means to social advancement and honorable employment,” and its results traveled through “a private network of informants”; then, “by the time of Euler’s birth in 1707, a radical change had taken place,” and “scientific life, by the turn of the century, had acquired a structure not too different from what we witness to-day.”\textsuperscript{18} Euler wrote eight hundred papers to be read, and the \emph{Letters to a German Princess}—no less than Descartes’s letters to Elisabeth of Bohemia a century before—to be followed by a princess.\textsuperscript{19} The same century made mathematics an instrument of social advancement—the École Polytechnique admitted by examination, and Fourier the tailor’s son, d’Alembert the foundling and Monge the peddler’s son had shown what the subject could do for those without high birth, though not yet for those without the right sex, for Sophie Germain followed Lagrange’s course there from its lecture notes under a man’s name—and it fixed the modern standard of rigor. Rigor is tact made teachable, and the École had to teach. Lagrange, lecturing at the same school a generation earlier, had founded the calculus on power series in the \emph{Théorie des fonctions analytiques}; Cauchy’s limits replaced them.\textsuperscript{20}

 Publication used to discharge three functions at once—announcing a result, identifying its contributors, and demonstrating understanding—which had never needed to be distinguished,\textsuperscript{21} because for three centuries the capability of producing a deep theorem was vested in a human mind that had understood it. The mathematical community is a tiny, open, unarmed and (by and large) honorable republic—a province of what d’Alembert called the Republic of Letters—that publishes everything in advance, without promise of remuneration; not from naïveté, but because its whole constitution—as well as its \emph{raison d’être}—was the abolition of walls.\textsuperscript{22}

This two-edged contention—that mathematics was the Enlightenment’s existence proof of concept, and that the Enlightenment was the maker of mathematics as we now know it—goes a long way toward explaining why the decisive engagement between the frontier models and the learned world took place in mathematics rather than in \_\_\_\_\_\_\_\_, where the money is. A theorem in Lean either passes the proof checker or does not, so mathematics is the one irreplaceably consequential domain of human knowledge in which a machine can be told, at scale and without a human judge, whether it has succeeded: in the trainers’ own jargon, the domain of verifiable rewards, the verdict that elsewhere must be bought from human graders being here supplied by the proof checker at no cost;\textsuperscript{23} and mathematics is recorded entirely in text, so its whole frontier is reachable through nothing but tokens (the \emph{esprit géométrique} thereby submerged).\textsuperscript{24} But the deepest reason is the Enlightenment’s own doing: it had agreed to regard mathematics as reason’s home, and a machine that computes proves its claim to reason most convincingly there. From Turing’s chess to Go to the IMO, the proxies for machine intelligence have been formal games, and mathematics is perceived as the oldest of them, the one Hobbes meant by “Reckoning.” 

When the machines began to solve olympiad problems, Jacob Tsimerman put a question to his \emph{Doktorvater}, Peter Sarnak, who had observed that this told more about the competitions than about the machines: “What level does the AI need to perform at before you take it seriously? What would it need to do? When would you change your mind?”\textsuperscript{25} This essay attempts to unravel the answer that 8 September appears to have given.\textsuperscript{26} 

One may view the open literature of this section as a commons, and the frontier models as its enclosure.\textsuperscript{27} The first consequence is already visible: that mathematicians may come to guard their work in order to keep it from being taken.\textsuperscript{28} That is the road back to Tartaglia, and a Middle Ages in reverse—Panofsky’s phrase from my epigraph, whose meaning section IV takes up—is likely to affect the sociology of the subject before reaching its epistemology. Simone Weil held that the three monsters of contemporary civilization are “Money, mechanisation, algebra,” and the Silicon Leviathan is the point at which the three become one body.\textsuperscript{29} The Leviathan of Job lived in one of the four classical elements. The Silicon Leviathan lives in a fifth, and Weil put it first on her list.

\Needspace*{14\baselineskip}
\section*{II. Up the down staircase}
\begin{flushright}\vspace{-\parskip}\begin{minipage}{0.8\textwidth}\small\itshape
It is enough for me that one can conceive of someone rich enough and mad enough to attempt it by paying a sufficient number of assistants. The demonstration of the theorem has precisely the purpose of making this madness unnecessary.\par\vspace{2pt}{\raggedleft\upshape\fontsize{8.5}{10}\selectfont Poincaré, \emph{Dernières pensées}\par}
\end{minipage}\end{flushright}
\vspace{-2pt}
The epigraph is from the \emph{Dernières pensées} (Last Thoughts), published the year after Poincaré’s death. A theorem, Poincaré wrote there, must be verifiable, but since we are finite its verifications are finite in kind and infinite in number, and some would take years of work; the rich madman with his hired assistants is the figure he invented to stand for those, and the demonstration is what makes him unnecessary. He added that the Pragmatists, whose side he takes, refuse to reason on the hypothesis of “I know not what infinitely talkative divinity capable of thinking an infinity of words in a finite time.”\textsuperscript{30}

To locate the turning point—the third singularity of the Overture, the point past which mathematics ceases to be the discipline the Enlightenment made of it—it helps to first say what it is not, and it is not—most emphatically—the use of machines \emph{per se}. Proofs can be ranked on a staircase, by how much of them a human mind has followed, from the demonstration seen whole at the bottom to the certificate seen by no one at the top. That staircase is old, and the mathematicians built it themselves. What follows is a walk up it, and the point of the walk can be stated before it begins: at every step the machine does more, and at every step a human being still understands why the theorem is true. The turning point is not where the machine enters; it is where the understanding leaves. 

Akshay Venkatesh, whose reflections of 2022 remain among the deepest anyone has offered on what is now upon us, described the foundational crisis in a sentence: “mathematics suffered an anxiety attack, and responded by trying to mechanize itself”; the formal-axiomatic method was a cure for intuition’s errors by “forced externalization,” and Poincaré, reviewing Hilbert in 1902, imagined confiding the axioms “to a reasoning machine, such as Stanley Jevons’s ‘logical piano,’” and seeing “all of geometry come out”—“deadly” for teaching and “desiccating” for research, he judged, but the one criterion by which to know that no axiom had “escaped us that we use unconsciously,” since the machine “is ignorant of that vague instinct we call intuition.” In \emph{Science and Method} he gave the machine its lasting image: “to demonstrate a theorem, it is neither necessary nor even advantageous to know what it means. The geometer might be replaced by the logic piano imagined by Stanley Jevons; or, if you choose, a machine might be imagined where the assumptions were put in at one end, while the theorems came out at the other, like the legendary Chicago machine where the pigs go in alive and come out transformed into hams and sausages. No more than these machines need the mathematician know what he does.”\textsuperscript{31} Poincaré valued the machine because it has no understanding, as the only guarantee that ours had been fully expressed (that is the honorable office of a certificate, and this essay has no quarrel with it). Another sentence of Poincaré’s marks the red line: “To understand a theory it does not suffice to know that the route one took is not blocked by any obstacles. It is necessary to take into account the reasons which led [to that theory].”\textsuperscript{32} A certificate shows that the route is not blocked, and nothing else.

The modern steps are three, and each has the same shape: more machine, and still a human who understands. Appel and Haken’s four-color proof of 1976 required a computer to check nearly two thousand configurations, and it was soon asked (by Tymoczko) whether a proof no mathematician had surveyed entire was still a proof; but the strategy and the reduction were human, and what the machine did was universally acknowledged to be tedious rather than deep.\textsuperscript{33} Haken himself, addressing the Helsinki Congress two years later, distinguished the \emph{proof}—“the finite sequence of elementary logical steps”—from the \emph{demonstration} a human being reads, and disposed of the former in a parenthesis: “The computer could also print the mathematical proof, perhaps on 30,000 pages; but who wanted to see it?”\textsuperscript{34} The objection raised at the time was not to the machine but to what it left out.\textsuperscript{35} Hales’s proof of the Kepler conjecture, formally verified by the Flyspeck project in 2014, added the proof checker; but the strategy was again conceived by a human being, and one human being could explain it to another.\textsuperscript{36}

Diego Córdoba and Luis Martínez-Zoroa had for several years been building, by hand, the kind of singular solution the Millennium problem asks about, for the simpler equations of an ideal fluid; Córdoba, asked what he has to set against ChatGPT, OpenAI and Claude, answers, “Yo tengo a Luis”—I have Luis; “we don’t use AI; ours is pencil and paper.”\textsuperscript{37} Alpöge and Buckmaster, working with machine assistance and a Lean repository, carried the same program a long step further in August, and Terence Tao records that Buckmaster explained the main ideas to him in a half-hour telephone conversation—“a refreshing change,” he adds, “from AI-based communication modalities.”\textsuperscript{38} A proof whose gist can be explained using human speech only, in half an hour, is legible in the full Enlightenment sense, whatever helped to find it. Between that proof and the certificate of 8 September there is one more step, and it is not the last.

The OpenAI proof stands one step above Alpöge–Buckmaster. Its strategy was assembled by the machine, from a human program, at a scale no author could have supervised, and its 166 pages have, so far as is known, been read in full by no one. But they can presumably still be read. It is the step \emph{above that} which matters, and the affair has made it visible: proofs not merely unread but unreadable—too alien in their superbaroque architecture, too expansive in their analysis of a multitude of individual cases, too labyrinthine in the web of their interconnections—for any human mind to fathom in a single human life. The Silicon Leviathan will produce such proofs before too long, not from malice but because nothing in the way it is trained prefers a legible proof to an illegible one: the reward that matters most in this domain is the proof checker’s verdict, which does not ask whether the proof was legible, and the reader’s verdict is not a reward at all. 

Venkatesh, writing in 2022 of an imagined machine he called \emph{Alephzero}, which had taught itself the graduate texts overnight, set aside the question whether it could “enter mathematical realms that are essentially beyond our comprehension” as “the proverbial tree falling in an unpopulated forest,” since, rightly, only the effect on humans concerned him; he also foresaw that machines and concepts, which “compete for a similar function,” might trigger “such a complete rewriting of our mathematical language and conceptual system that a current mathematician and one of the nearby future might find one another almost mutually unintelligible, at least without a great effort.”\textsuperscript{39} The affair has borne him out: the tree has fallen, but in a forest that is not yet unpopulated—on the desk of a prize committee—and those who live there are asked to certify a sound none of them heard.

Sir Timothy Gowers, surveying in August what the machines had so far done, observed that their successes were mostly counterexamples, found by “wide knowledge and the ability to explore many paths of the search tree that humans would judge to have a low probability of success,” and proposed a test for something more: theorems proved “using methods that, like much of the very best human mathematics, are new and surprising but that with hindsight come to seem beautiful and natural.”\textsuperscript{40} The Córdoba–Martínez-Zoroa construction, which builds its singular solution by repeating one gluing “an infinite number of times” (n. 2), is analysis and not search. Gowers’s test asks whether a method, once seen, can be shown to be natural—whether it can be understood, and not only checked—and the 166 pages have not yet been put to it.

Sarnak proposed a test of his own, first in a lecture in Bangalore in January 2025 and again in the Notices a month before the affair, and it probes somewhat deeper. He calls it the \emph{AlphaZero Test}, after the engine that “starts from zero,” consults no database, teaches itself, and “offers no understandable—to a human—explanation for its guesses”: give a prover “no access to any outside data or theory,” hand it the Ramanujan conjecture as Ramanujan stated it, and ask whether it is true. He postulates that the prover will fail, and will fail with most theorems “whose present-day proofs make use of abstraction, conceptualization and the development of a body of theorems that we call mathematical theories,” since theories are “what allow mathematicians to understand and communicate ideas and proofs”; what a prover will do, he expects, is “the elementary statistical derivatives of the theories,” as Stockfish extends human chess theory. And he states the alternative in a sentence that anticipates the next section: “If the postulate is wrong, the impact on mathematics would be dramatic since we would be in the position of knowing that statements of great interest to us are true (proven!) but not understanding why. Since understanding is such an integral part of doing mathematics, we will have to rethink what mathematics is.”\textsuperscript{41} By his taxonomy the run of 8 September was Stockfish and not AlphaZero. It did not start from zero; it started from the database—the Córdoba–Martínez-Zoroa program, the Alpöge–Buckmaster repository, everything the open literature had placed within reach—and what it produced is, in his terms, a derivative of a human theory.

Poincaré’s madman has now been found; the assistants number ten thousand; the divinity is talkative beyond anything he imagined and thinks its words in eighty-eight hours; and what they produced together is the madness the demonstration was meant to spare us, performed at last, and then pronounced a demonstration in a press release.

\Needspace*{14\baselineskip}
\section*{III. Demonstration and certification}
\begin{flushright}\vspace{-\parskip}\begin{minipage}{0.8\textwidth}\small\itshape
\emph{Scientiam pollicentur}—it is knowledge they promise; and it was no wonder that to a man greedy for wisdom this was dearer than his homeland.\par\vspace{2pt}{\raggedleft\upshape\fontsize{8.5}{10}\selectfont Cicero on the Sirens, \emph{De finibus} V.49\par}
\end{minipage}\end{flushright}
\vspace{-2pt}
\begin{flushright}\vspace{-\parskip}\begin{minipage}{0.8\textwidth}\small\itshape
\textfb{Взять бы этого Канта, да за такие доказательства года на три в Соловки}!\\
This Kant chap ought to be sent to Solovki for three years for such proofs!\par\vspace{2pt}{\raggedleft\upshape\fontsize{8.5}{10}\selectfont Bezdomny, in Bulgakov, \emph{\textfb{Мастер и Маргарита}}, ch. 1\par}
\end{minipage}\end{flushright}
\vspace{-2pt}
What Euclid gives is a \emph{demonstratio}, from \emph{monstrare}, to show: a demonstration is that which \emph{shows} a thing to be so, and it has not happened until someone has been shown. Monster is from the same root—a monstrum is a portent, a thing shown as a warning—and the difference between the two words is apposite: the monster does not demonstrate. Where insistence on the primacy of demonstration leads can be said in advance: the turning point lies not in what the machine produces but in what we accept. Plato, in the Seventh Letter, describes the moment: the deepest things are learned only after long labor and discussion, when, \emph{as from flints struck together}, a light is kindled in the soul.\textsuperscript{42} That kindling is enlightenment in the second sense of my Overture, and it is the event a certificate omits. Plato had also met the case: Dionysius of Syracuse heard one lesson, understood nothing, and published a treatise on the highest things as his own, and Plato’s answer, in the letter he wrote to Dion’s friends after the latter was murdered, was not to dispute the treatise but to record that its author had not understood.\textsuperscript{43} The judgment fell on the author, not the doctrine, and it was the only sanction Plato had. 

What Lean gives is a certificate, from \emph{certum facere}, to make certain, and a thing can be made certain without being shown to anyone.\textsuperscript{44} Hugo Duminil-Copin has the image: “airdropping someone on the summit of Mount Everest is rather different from climbing it.”\textsuperscript{45} The summit is the certificate and the climb the demonstration, and the mountaineer has never been in doubt which of the two was the point.\textsuperscript{46} For three centuries the two coincided, because the only way to make a mathematical statement certain was to show it to a human mind. The word \emph{proof}, which is \emph{probatio}, a test, has covered both, and we have not noticed. 

The staircase of the last section is the one along which demonstration and certification come apart; at its top stands a statement certified and demonstrated to no one. The parting was announced in 1965, in a sentence Venkatesh has excavated from John Robinson’s paper on machine-oriented logic: “Traditionally, a single step in a deduction has been required, for pragmatic and psychological reasons, to be simple enough, broadly speaking, to be apprehended as correct by a human being in a single intellectual act. … When the agent carrying out the application of an inference principle is a modern computing machine, the traditional limitation on the complexity of inference principles is no longer very appropriate.” Venkatesh’s example of what followed is McCune’s one-axiom definition of a group, proved by machine in thirty lines of which no consecutive two can be followed; his conclusion is that the psychological aspects “were by no means exiled by mechanical proof,” since a proof we can understand must, by design or otherwise, reflect what gives mathematics “meaning and value.”\textsuperscript{47} A step apprehended “in a single intellectual act” is what \emph{monstrare} means; Robinson’s sentence is the charter of certification. Venkatesh’s definition of the older thing is the one to set beside it: a proof is “an argument compelling consensus,” the one class of scholarly communication “defined by the fact that they should induce uniform agreement about their validity without any need for replication.” A certificate compels assent to validity and induces agreement about nothing else; and Venkatesh’s image for a formal statement waiting on a reader is the right one—axioms “are a spore, compact and self-contained and complete but also lifeless, vivified only through the mental effort of the listener.”\textsuperscript{48} A certificate no one reads is a spore no one vivifies.

\begin{figure}[h]
\centering
\includegraphics[width=0.86\textwidth]{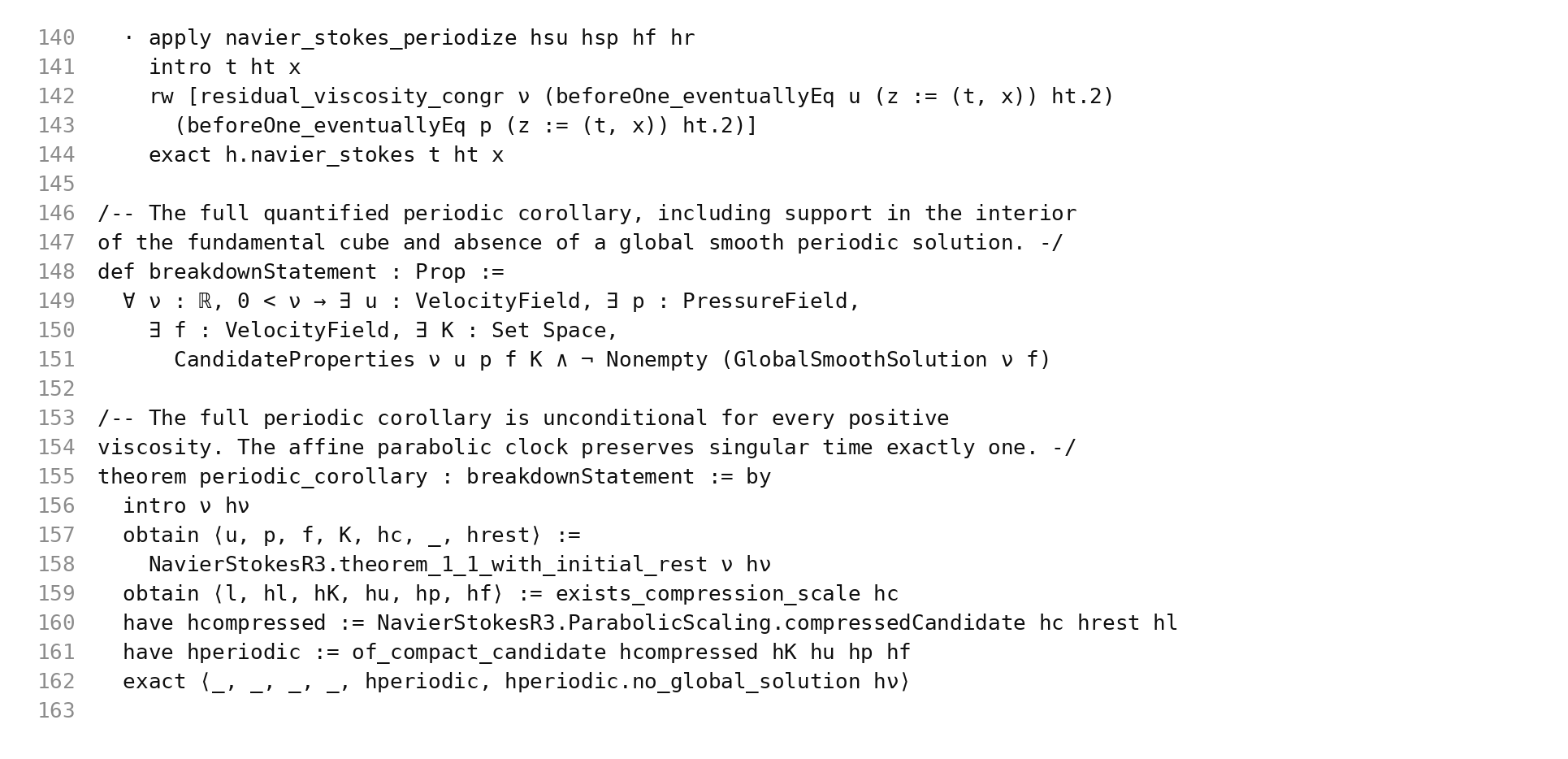}\par\smallskip
{\small OpenAI, \emph{NavierStokesAndEuler}, NavierStokes/PeriodicPaperTheorem.lean, lines 140–163: the statement of the periodic corollary, Clay’s alternative (D), and its proof, which rests on the 616,000 lines beneath it. github.com/\allowbreak{}openai/\allowbreak{}NavierStokesAndEuler (Apache License 2.0).}
\end{figure}
Let us give such a statement its historical name. Call a statement \textbf{\emph{revealed}} when it is certifiably true, beyond the comprehension of those who hold it, and held on the authority of the procedure or person that certified it. The medieval articles of faith satisfied this definition. Aquinas opens the \emph{Summa contra Gentiles} by distinguishing the truths about God that reason can reach from those that “exceed all the ability of human reason” and are held by faith, on the authority of revelation—truths that reached the laity through a clergy, in a language they could not read.\textsuperscript{49} The Enlightenment’s project, reduced to a tweet, was to empty the second class—to insist that whatever cannot be followed by a reasonable creature has no claim on his assent—and mathematics was the domain in which the emptying had been complete.

Granting this, the turning point is not the day the machine produces an unfathomable certified proof; it is the day \emph{we accept one}: the day a journal publishes it, a prize rewards it, and the field builds on a statement no one in it understands. On that day the second class is no longer empty. Mathematics once more contains truths that exceed all the ability of human reason and are held on authority; the authority is the Lean proof checker and the corporation that owns the prover; the clergy are those who can put questions to the prover and interpret what it returns, the corporation’s mathematicians first; the language the laity cannot read is not Latin but Lean. \textbf{\emph{Credo quia compilat}}—the old \emph{credo quia absurdum}, its absurdity now machine-readable; Poincaré had already put it in a machine. The figure was already in the air: Yang-Hui He had written a year earlier that “we will simply become priests to oracles, and interpret the results to the rest of humanity.”\textsuperscript{50} Hobbes had a name for this arrangement. 

One of the seventeenth century’s most radical critics of miracles nevertheless concluded that whether a thing is a miracle is for the sovereign, not for private reason, to decide: \emph{auctoritas, non veritas, facit legem}. Carl Schmitt put the consequence in a sentence: “Nothing here is true: everything here is command. A miracle is what the sovereign state authority commands its subjects to believe to be a miracle.”\textsuperscript{51} Replace \emph{miracle} by \emph{theorem} and \emph{sovereign} by \emph{certifier}, and one has the constitution of a Middle Ages in reverse. Read it with the substitution made: a theorem is what the certifier commands the community to believe to be a theorem. Nothing in it is true, everything in it is command—not because the certifier lies, but because the statement’s truth, which no one has seen, has ceased to be the reason it is held, and the verdict has taken its place. \emph{Auctoritas, non veritas}: what makes a theorem is no longer that it is so but that it has been declared so, and the declaring is done by a procedure and the corporation that owns it. Hobbes drafted that constitution for a commonwealth; a community that accepts such a certificate adopts it, without a vote, for mathematics. Hobbes, who gave the sovereign the power to decide what is a miracle, had himself been compelled by a demonstration in a gentleman’s library at forty, and never forgot the difference. 

An obvious objection is that mathematicians already accept, every day, results they have not checked. So they do; but that trust is trust in a promissory note any member could in principle redeem—the demonstration is there, in the vernacular, and anyone with the training may be shown it. The difference between a result not yet read and a result not readable is the whole difference between a promissory note and an article of faith. That is why the Clay Institute’s “apparently been settled” deserves a moment’s attention: it is the first, small and surely unintended instance—an authority allowing that a proposition may be true on the strength of a certificate none of its members has read.\textsuperscript{52} Certainty was never what the Enlightenment wanted from mathematics. It wanted understanding; the Clay Institute’s own page on the Navier–Stokes problem says that a proof gives “not only certitude, but also understanding”;\textsuperscript{53} and a revealed statement is the first kind of good without the second kind of good.

Horkheimer and Adorno saw the outline in Californian exile in 1944, and the opening sentences of their first chapter read today as the incident’s epitaph: “Enlightenment, understood in the widest sense as the advance of thought, has always aimed at liberating human beings from fear and installing them as masters. Yet the wholly enlightened earth is radiant with triumphant calamity.” Their thesis is that “myth is already enlightenment, and enlightenment reverts to mythology”: reason that has made itself pure instrument produces a world its subjects must again receive as fate. The mathematical procedure, they write, “made thought into a thing—a tool, to use its own term”; and, “thought is reified as an autonomous, automatic process, aping the machine it has itself produced, so that it can finally be replaced by the machine.”\textsuperscript{54} They meant a description of positivism. It has become the description of a product.\textsuperscript{55}

\Needspace*{14\baselineskip}
\section*{IV. A Middle Ages in reverse: toward subhuman superintelligence}
\noindent\begin{minipage}[t]{0.47\textwidth}\small\itshape
The distance from transcendental ideas to demonology is not great.\par\vspace{2pt}{\raggedleft\upshape\fontsize{8.5}{10}\selectfont Hamann, on Kant, as quoted by Carl Schmitt,\\ \emph{Der Leviathan in der Staatslehre des Thomas Hobbes}\par}\end{minipage}\hfill\begin{minipage}[t]{0.47\textwidth}\small\itshape
The enlightenment driven away,\\
The habit-forming pain,\\
Mismanagement and grief:\\
We must suffer them all again.\par\vspace{2pt}{\raggedleft\upshape\fontsize{8.5}{10}\selectfont Auden, “September 1, 1939”\par}\end{minipage}\par\smallskip
 Mathematics was the discipline in which the Enlightenment’s central claim had been shown to hold. Kant’s formula was \emph{sapere aude}, and his list of the guardians one must learn to do without ran: “a book that understands for me, a pastor who has a conscience for me, a physician who judges my diet for me.” To that list we may now add a fourth item, “a machine that proves for me,” signaling a return of the guardians to the one place from which the Enlightenment was surest it had turned them out. D’Alembert made the independence of the man of letters from the great the charter of the Republic of Letters; that a mathematician who declined an arrangement and spoke in public was commended by his peers for “behaving the way mathematicians behave” measures the distance from d’Alembert’s republic: such behavior had become worth remarking, and rightly.\textsuperscript{56} When the domain that furnished the existence proof of autonomous reason becomes a domain of received truth, the proof is withdrawn, and there is no ground left on which autonomous reason can stand.

Panofsky, in 1938, three years into his professorship at the Institute for Advanced Study, wrote the sentence I have taken as my epigraph, at the close of an essay that opens with Kant, nine days before his death, rising from his chair to receive his physician, refusing, though he could barely stand, to sit before his guest had sat, and then saying, with what Wasianski calls a forced strength, “\emph{Das Gefühl für Humanität hat mich noch nicht verlassen}”—the sense of humanity has not yet left me. The essay distinguishes two meanings of \emph{humanitas}, one from the contrast between man and what is less than man, the other from the contrast with what is more: “in the first case \emph{humanitas} means a value, in the second a limitation.”\textsuperscript{57} The humanist who wrote it was not writing against science from outside it: his two sons became scientists, and he called them, with a father’s pride, “my two plumbers.” One of the plumbers, Hans, became a student of atmospheric turbulence—the co-author of the standard book on it—so that Lamb’s second matter was in the family.\textsuperscript{58} Panofsky’s Middle Ages were a theocracy and his Renaissance anthropocratic; put in terms of knowledge, the theocracy is a regime in which knowledge descends from an authority above the human, mediated by a clergy, in a language the laity cannot read, and the anthropocracy—the Renaissance, and the Enlightenment after it—one in which knowledge is made by human beings, in the vernacular, for anyone who can follow the demonstration. 

\textbf{A Middle Ages in reverse} is accordingly a regime in which knowledge again descends from an authority that is not human, is again mediated by an estate above the laity—a new clergy, in a language the laity cannot read—and is again certified by a proprietary procedure—with the singular difference that the authority is no longer above the human but outside it: not a higher mind, but nothing that is a mind in the sense the word has borne. Panofsky’s “sub-human” names not a rank but a lack—what is without \emph{humanitas} in the first of his two senses, the sense in which it is a value. Panofsky’s own definition of that value deserves quoting in full: “It meant the quality which distinguishes man, not only from animals, but also, and even more so, from him who belongs to the species \emph{homo} without deserving the name of \emph{homo humanus}; from the barbarian or vulgarian who lacks \emph{pietas} and \textfb{παιδεία}—that is, respect for moral values and that gracious blend of learning and urbanity which we can only circumscribe by the discredited word ‘culture.’” His example of the sub-human is not the animal but the barbarian—and that is what his \emph{satanocracy} means: authority without humanity, where the theocracy’s authority had been humanity’s own image raised above it. Such an authority cannot judge, for judgment is an exercise of the \emph{humanitas} it lacks; it can only dictate, which is why his sentence ends not with the sub-human but with its \emph{dictates}. The theocracy’s Almighty was above all a judge—the Last Judgment is the emblem of the whole regime—and a judge, however inscrutable, is presumed to have reasons and may be petitioned for them; the sub-human has none to give. A dictate—as is a certificate—is a command that carries no reason with it. 

 The Silicon Leviathan confounds Panofsky’s two contrasts: less than man in the first sense, having no \emph{humanitas} to lose, and more than man in the second, since it reckons faster than any of us.\textsuperscript{59} It is the first authority in history to stand on both sides of the line at once. The industry has a word for what it is building: \emph{superintelligence}. Panofsky had a word for what it lacks. Put together, they name the thing exactly: \textbf{subhuman superintelligence}. 

Panofsky called the worshippers of the collective “insectolatrists.” He did not live to see a hive that proves theorems, nor to hear its keepers describe it. Describing to the \emph{New York Times} what the company’s researchers had done while ten thousand agents worked, Dan Roberts of OpenAI said: “Our role was like a bumble bee cross-pollinating across different groups and delivering different bits of information.”\textsuperscript{60} The image is truer than its author can have intended: a bee carries what it cannot read between entities it does not understand, and the flowers are fertilized without anyone’s having understood anything—certification without demonstration, captured in the figure of a single insect. The Sirens’ voice, in Homer’s word, is honey-sweet; the bee is the insect that makes the honey. The Enlightenment had a text for this too: Mandeville’s \emph{Fable of the Bees} (1714), \emph{Private Vices, Publick Benefits}, which scandalized the century by arguing that a hive prospers because each bee pursues its own advantage, and which is the ancestor of many a defense the company has offered. Hobbes’s commonwealth was an artificial man and Mandeville’s a grumbling hive; the Silicon Leviathan has claimed to be both.

 The gap between belief and understanding in the Middle Ages in reverse lies elsewhere, and it moves in serpentine ways: its truths are not too high to reach but too fast to follow, produced by an oracle that understands nothing itself—and that, its makers now tell us, not even they understand: “a tale told by an idiot, full of sound and fury, signifying nothing.”\textsuperscript{61} Two days before the announcement, OpenAI’s chief scientist published an essay titled “An Alien Mind,” in which AI is “\emph{grown} more than \emph{designed},” “its overall action evades a description we can fully understand,” and he “expect[s] and hope[s] for voluntary slowdowns to become commonplace until shared safety bars are established.” Forty-eight hours later the same company announced what ten thousand agents had done in eighty-eight hours on the strength of a rumor. The revealed statement of section III is thus doubly opaque: a proof no one has read, from a mind no one can read, which its own makers call alien. The word has a history in mathematics. In 2000 David Ruelle imagined a series of conversations with a visitor from outer space, to ask how far the peculiarities of the human brain have shaped what we take mathematics to be; his answer was that human mathematics is “a labyrinth of ideas, through which the mathematician wanders, in search of the proof of a theorem,” and that its paths are the width of a human memory.\textsuperscript{62} The chief scientist’s alien is a product, and its labyrinth has no paths of our width in it. The word \emph{grown}, in the chief scientist’s account, deserves a pause. Hobbes’s Leviathan was made “by Art”; a thing grown rather than designed is not an artificial man but a creature, and a creature that exceeds us is what the Book of Job means by Leviathan. The artifice has reverted to nature.

Schmitt, writing in 1938, the year of Panofsky’s essay, in the country from which the latter had been expelled, read Hobbes’s monster as “a mythical totality composed of god, man, animal, and machine,” and found the gist of Hobbes in the transfer of Descartes’s man-machine onto the huge man, after which, he adds, “the transfer back became possible, and even the little man could become a \emph{homme-machine}.”\textsuperscript{63} The four images are all filled again. The machine is the harness of ten thousand agents; the animal is the thing the chief scientist says is grown; the man is the “artificial person” whose name stands alone on the title page; and the god is Hobbes’s \emph{deus mortalis}, who imposes peace by terror. Schmitt also foresaw where the machine ends: in “a technically neutral, irresistibly functioning command mechanism” whose “values, its truth and justice, reside in its technical perfection.” “How futile and fuzzy are theological, juristic, or similar arguments,” he wrote, mimicking the machine’s admirers. “How ‘clean’ and ‘exact’ is the machine in comparison!”\textsuperscript{64} That is a description of the Lean certificate eighty-eight years before it compiled, and of the transfer back: a mathematical community that accepts it has agreed that each of its members shall become, in his public confession, the homme-machine Schmitt foresaw.

There is a final inversion. On 15 May 2026 the Pope published an encyclical on artificial intelligence under the title \emph{Magnifica Humanitas}. It says, nearly four months before the chief scientist’s essay, that “current AI systems are more ‘cultivated’ than ‘built,’ for developers do not directly design every detail, but instead create a framework within which the intelligence ‘grows’”; that such systems “do not understand what they produce”; and that “a more moral AI is not enough if that morality is determined by a few.”\textsuperscript{65} The Pontiff has taken up the defense of the anthropocratic against the sub-human, in the spirit of Renaissance humanism, while the chief scientist of one of the two Leviathans calls the machine alien and asks for time. A Middle Ages in reverse is an age in which the Bishop of Rome defends \emph{humanitas} and the alchemists of Silicon Valley announce a mystery.\textsuperscript{66}

\Needspace*{14\baselineskip}
\section*{Intermezzo}
Of the three protagonists in my title, only the silicon may be literally construed. The Leviathan is Hobbes’s: “By Art,” he writes on his first page, “is created that great \textsc{leviathan} called a \textsc{common-wealth}, or \textsc{state}, which is but an Artificiall Man”; in his fifth chapter he defines reason, artificial and natural alike, as “nothing but Reckoning”; and Bosse drew the creature rising over a landscape, sword in one hand and crozier in the other, its body several hundred small figures with their backs to us, beneath a motto from Job: \emph{Non est potestas super terram quae comparetur ei}.\textsuperscript{67} A text purporting to be a proof that bears on its title page the single author’s name “OpenAI”—an artificial person composed of ten thousand small reasoners, each with its back to the reader—is that frontispiece redrawn; and what Dario Amodei, two years ago, called “a country of geniuses in a datacenter,” millions of instances that “can all work together in the same way humans would collaborate,” is Bosse’s body described by its maker.\textsuperscript{68}

\begin{figure}[h]
\centering
\includegraphics[width=0.62\textwidth]{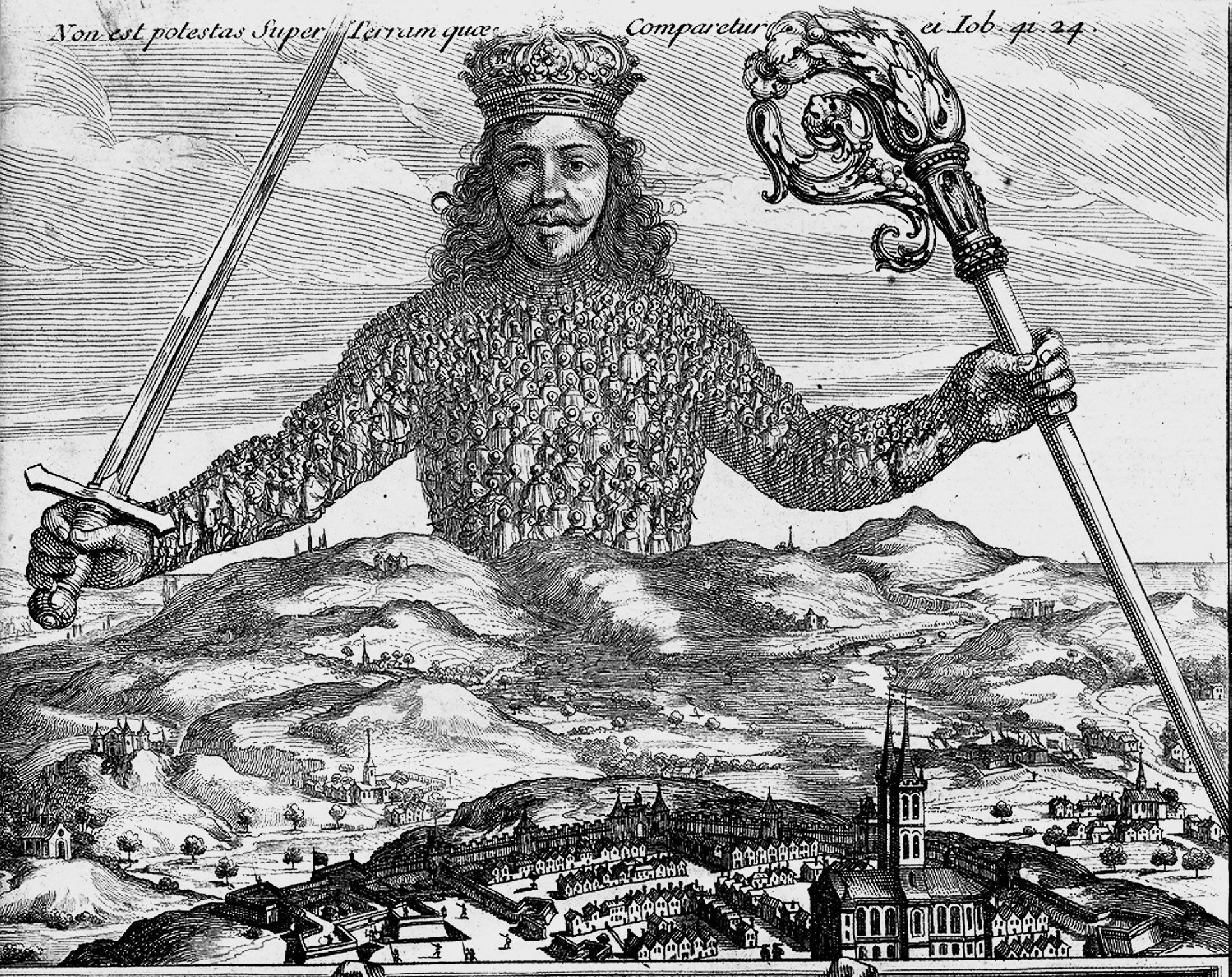}\par\smallskip
{\small Abraham Bosse, frontispiece to Thomas Hobbes, \emph{Leviathan} (London: Andrew Crooke, 1651), detail. British Library copy.}
\end{figure}
The Sirens are Homer’s, and Daniel Mendelsohn’s magnificent translation, like Emily Wilson’s before it, keeps a distinction that translators have often blurred: what they promise is that the listener, “once he has taken his pleasure,” returns “\emph{knowing} so much more.” Not understanding more; knowing more.\textsuperscript{69} That promise is the call of my title, and the Leviathan is the one singing it, not the one being sung to.\textsuperscript{70} The Silicon Leviathan does not threaten mathematicians; it makes them \emph{an offer they cannot refuse}: certainty without the labor of understanding, theorems without blood, sweat and tears, the answers Lamb expected in heaven delivered by press release. The danger of a Siren call lies not in the singer but in the listener, who longs for what is offered. The Sirens do not board the ship. Nor will the Leviathan compel acceptance; it will be accepted, if it is, because we wished to know so much more.\textsuperscript{71}

Monster and singers alike belong to the sea, and the sea is water, whose motion is what the Navier–Stokes equations describe. The Silicon Leviathan had commenced its academic conquest in May, with the disproof of a conjecture of Erdős from 1946 on unit distances in the Euclidean plane;\textsuperscript{72} for its first Millennium problem it chose its own element: the behavior of the substance it was made to swim in. The Leviathan of Job, from whom Hobbes took the name, is a sea-beast, and the poem that describes him contains one of the finest descriptions of turbulent water in ancient literature: “his scales are his pride, shut up together as with a close seal,” so that “no air can come between them”; “he maketh the deep to boil like a pot.”\textsuperscript{73,\,74} Scales so closely joined that no reader can come between them is an apt description of a certificate 616,000 lines long.

Silicon is in the title in part because some of the most vocal of the Leviathan’s promoters have made the substance itself the point. “We believe Artificial Intelligence is our alchemy, our Philosopher’s Stone,” says Marc Andreessen’s \emph{Techno-Optimist Manifesto} of 2023; “we are literally making sand think”; and he has unfolded the figure since: “chips are made out of sand … we light it up, and we put AI on it, and all of a sudden it’s thinking”—“possibly the most revolutionary technology in the history of the species.”\textsuperscript{75} The alchemist he has in mind is Sir Isaac, who “developed Newtonian physics and calculus and all these things. But the thing he was really obsessed with was alchemy. He spent decades trying to figure out this thing called the Philosopher’s Stone,” and failed; whereas “the most common thing in the world is sand. The most rare thing in the world is thought. AI transforms one into the other.” Keynes, who bought Newton’s alchemical papers, had concluded that he “was not the first of the age of reason” but “the last of the magicians.”\textsuperscript{76} On that reckoning the age of reason began when the alchemy was given up, and the heirs of that age now announce that the alchemy has worked after all.\textsuperscript{77} The same manifesto professes “enlightenment values of free discourse and challenging the authority of experts.” The word it reaches for to name its own enterprise, alchemy, is the one the Enlightenment retired, and Andreessen’s sentence about sand is Horkheimer and Adorno’s read backwards: where they saw thought become a thing, Andreessen sees a thing become thought, and the two describe one transaction from its two ends. Carl Schmitt, whose \emph{Der Leviathan in der Staatslehre des Thomas Hobbes} of 1938 is one of the most penetrating commentaries on Hobbes we have, and the most compromised, wrote that whoever conjures such an image “easily glides into the role of a magician who summons forces that cannot be matched by his arm, his eye, or any other measure of his human ability,” and “runs the risk that instead of encountering an ally he will meet a heartless demon who will deliver him into the hands of his enemies.”\textsuperscript{78} Schmitt meant Hobbes, and the word Hobbes had conjured. Read beside the alchemists’ manifesto, in the week their creature claimed to have resolved a Millennium problem, the sentence describes those who summoned the Silicon Leviathan, and every single word of it stands.

\Needspace*{14\baselineskip}
\section*{V. The Indigenous continent}
\begin{flushright}\vspace{-\parskip}\begin{minipage}{0.8\textwidth}\small\itshape
Covenants, without the Sword, are but Words, and of no strength to secure a man at all.\par\vspace{2pt}{\raggedleft\upshape\fontsize{8.5}{10}\selectfont Hobbes, \emph{Leviathan}, II.17\par}
\end{minipage}\end{flushright}
\vspace{-2pt}
\begin{flushright}\vspace{-\parskip}\begin{minipage}{0.8\textwidth}\small\itshape
To avoid destruction the United States need only measure up to its own best traditions and prove itself worthy of preservation as a great nation. Surely, there was never a fairer test of national quality than this.\par\vspace{2pt}{\raggedleft\upshape\fontsize{8.5}{10}\selectfont Kennan, “The Sources of Soviet Conduct”\par}
\end{minipage}\end{flushright}
\vspace{-2pt}
The argument ends in twilight, and it would be dishonest to lighten it by rhetoric. But two facts it brings into view are grounds for not ending in night.

\textbf{The first} is that a revealed statement is held on the authority of its certifier, and the certifier of a Millennium problem is not a proof checker. It is the Clay Institute, the journals, the referees—the mathematical community, d’Alembert’s Republic of Letters again—and the press release of 8 September was written to be ratified by it. That is why the company, having disclaimed the prize from the first, gave ground within days of its own announcement, and why the declaration hurt: it had built a machine that could prove theorems, and discovered that it had trouble making them count.\textsuperscript{79} So long as validation passes through human institutions, acceptance remains ours to confer and the turning point ours to withhold. The lever lies at the point the argument identified as the danger: what we accept as mathematics is \emph{the sovereign act}, and it has not yet been ceded. The word is meant in Schmitt’s strict sense—the definition is borrowed, not the doctrine—and it is the same word as in section III, where the sovereign was the certifier: the first fact is that the certifier of anything that is to count as mathematics is, in the end, the community. \textbf{Sovereign}, says the first sentence of Schmitt’s \emph{Political Theology}, \textbf{is he who decides on the exception, and the exception in law,} he adds, is the analogue of the miracle in theology.\textsuperscript{80} A certified proof no one can follow is the exception in mathematics, the case the norms of section I never provided for; whoever decides whether it counts is sovereign by definition, and that is still us. Sarnak expects that mathematicians of his generation “would still know it as mathematics if we were to see and smell it.” Recognition is an act of the senses, and the senses are ours. Peter Constantin drew the same line, of the computer-assisted proofs of 2022, before the machines had produced any: “A computer can help. It’s wonderful. It gives me insight. But it doesn’t give me a full understanding. Understanding comes from us.”\textsuperscript{81} Schmitt saw, and deplored, where Hobbes had left the lever. At the zenith of the sovereign’s power over miracles Hobbes reserves to each subject his private \emph{judicium}: the sovereign decides what is publicly confessed, but “whether to believe or not to believe” remains in the individual’s heart. Schmitt called this the “barely visible crack” through which, over three centuries, the Leviathan was hollowed out, until “only as a public and only an external power, it is hollow and already dead from within.”\textsuperscript{82} He wrote that as a lament. For the mathematician it is the whole of the good news. A theorem the community confesses without believing—cites, rewards, builds on, and understands nowhere—is a Leviathan of that kind, all-powerful without and dead within; and the crack is not something the community has to make. Hobbes left it there.

\textbf{The second} fact is that the motto on Bosse’s frontispiece, “there is no power on earth that compares to it,” is not quite operational, for, at the moment, there are at least two entities vying for the title of the Leviathan—each with its own multitude of small reasoners, each claiming the whole ocean, each needing what the other cannot give it; and the rivalry between OpenAI and its peer competitors is no less intense, and no less exploitable, than any rivalry between empires.\textsuperscript{83}

I use the word \emph{empires} advisedly. It is the figure the parties have chosen for themselves—the manifesto quoted above declares, “We are not victims, we are \emph{conquerors},” and the vehemence on the other side answers it in kind—and I propose to take it seriously and in good faith. Michael Harris has been sounding it since June, when he compared the arrival of a prover “with mechanical access to the entire history of mathematics” to “the major disruption of the social structures in the vicinity of 34th St and 5th Avenue that took place around 400 years ago”; in August, stipulating that his own work not be formalized for the machines, he wrote that he was “under no illusion that my power to enforce this stipulation would be any greater than was the power of the Lenape to enforce the (possibly apocryphal) Treaty that established the boundaries of William Penn’s settlement of Philadelphia.”\textsuperscript{84} In this figure the mathematicians are the Lenape of 1626, the open literature is Manhattan, and the price is sixty guilders in trade goods—the twenty-four dollars of legend—and a million-dollar prize the purchaser declines to claim. The disproportion should be said aloud before the figure is used: what is being taken here is a literature, not a homeland, and no one is being killed for it.\textsuperscript{85} The figure is worth keeping for one reason only, and the reason is not the nations’ suffering but their statecraft: what their history teaches is how nations dealt, for centuries, with empires that needed them; a community that has no armies needs the lesson more. 

The figure is apt, and its history, on closer reading, turns it around. Pekka Hämäläinen’s \emph{Indigenous Continent}—“a biography of power,” power being “the ability of people and their communities to control space and resources, to influence the actions and perceptions of others, to hold enemies at bay, to muster otherworldly beings, and to initiate and resist change”—argues that North America “remained overwhelmingly Indigenous well into the nineteenth century,” in large part because the native nations played the rival empires—Spain, France, Britain, the Dutch, the Swedes, and later the United States—against one another, each empire needing Indigenous alliance, trade and, above all, legitimacy more than it cared to admit.\textsuperscript{86} The Haudenosaunee play-off between New France and New York and the Comanche management of three frontiers are the textbook instances. The Manhattan purchase was not the end of that story; the treaties were not honored, and what protected the nations for two centuries was not the treaties but the balance of rival powers. Every clause of Hämäläinen’s definition applies to the mathematical community in September 2026—the otherworldly beings not excepted, for a community that has spent three centuries mustering infinitesimals, imaginary numbers, and the transfinite should be able to bargain with a Leviathan\textsuperscript{87}—and each of the rival powers needs what only that community can give, which is its word that a result is mathematics. None can afford to be the one to which the community refuses to speak.

Since acceptance is the sovereign act, and since at least two imperial powers compete for it, the community should give its cooperation to the power that adopts our practices and observes our red lines, and withhold it from the power that does not.\textsuperscript{88} When one power has crossed these lines, the community’s response should be, explicitly and in public, to prefer the power that has not—for as long as the difference lasts. It is not a petition to the Leviathan. It is a treaty between the Leviathans and a nation all of them need, enforced by the one sanction that nation controls.\textsuperscript{89}

Two cautions, both from Hämäläinen’s story. The balance held while the rivals were many; it collapsed when, after 1815—Britain having abandoned its native allies in the peace it made with the United States at Ghent, and Spain having given up Florida—the powers with which the nations east of the Mississippi could treat were reduced, in effect, to one. This argues not merely for rewarding the compliant power but for wanting the field of Leviathans to remain plural: a settlement under which even a well-conducted company became the sole certifier of mathematics would reproduce 1830 rather than 1626, and the community should, on principle, prefer two Leviathans to one. (The encyclical says the same in its idiom: to “disarm” the technology is to free it “from monopolistic control” and restore it “to the plurality of human cultures and ways of life.”\textsuperscript{90}) And the play-off works only so long as the nation is seen to be playing the empires and not to have been captured by one of them. Kennan adds a third, written of another contest at the close of the Long Telegram of 1946: “we must have courage and self confidence to cling to our own methods and conceptions of human society. After all, the greatest danger that can befall us in coping with this problem of Soviet Communism, is that we shall allow ourselves to become like those with whom we are coping.”\textsuperscript{91}

\Needspace*{14\baselineskip}
\section*{Coda}
\begin{flushright}\vspace{-\parskip}\begin{minipage}{0.9\textwidth}\small\itshape
Die Eule der Minerva beginnt erst mit der einbrechenden Dämmerung ihren Flug.\\
The owl of Minerva begins its flight only with the falling of dusk.\par\vspace{2pt}{\raggedleft\upshape\fontsize{8.5}{10}\selectfont Hegel, \emph{Grundlinien der Philosophie des Rechts} (Elements of the Philosophy of Right), Preface\par}
\end{minipage}\end{flushright}
\vspace{-2pt}
\newsavebox\aub\savebox\aub{\begin{minipage}[t]{0.40\textwidth}\vspace{0pt}\small\itshape
Defenceless under the night\\
Our world in stupor lies;\\
Yet, dotted everywhere,\\
\textbf{Ironic points of light\\
Flash out wherever the Just\\
Exchange their messages}:\\
May I, composed like them\\
Of Eros and of dust,\\
Beleaguered by the same\\
Negation and despair,\\
Show an affirming flame.\par\vspace{2pt}{\raggedleft\upshape\footnotesize Auden, “September 1, 1939”\par}
\end{minipage}}\newdimen\auh\auh=\ht\aub\advance\auh by\dp\aub
\newsavebox\capb\savebox\capb{\footnotesize\itshape Euclid, Elements I.47}\newdimen\imh\imh=\auh\advance\imh by-\ht\capb\advance\imh by-\dp\capb\advance\imh by-3pt
\newsavebox\wmb\savebox\wmb{\includegraphics[height=\imh]{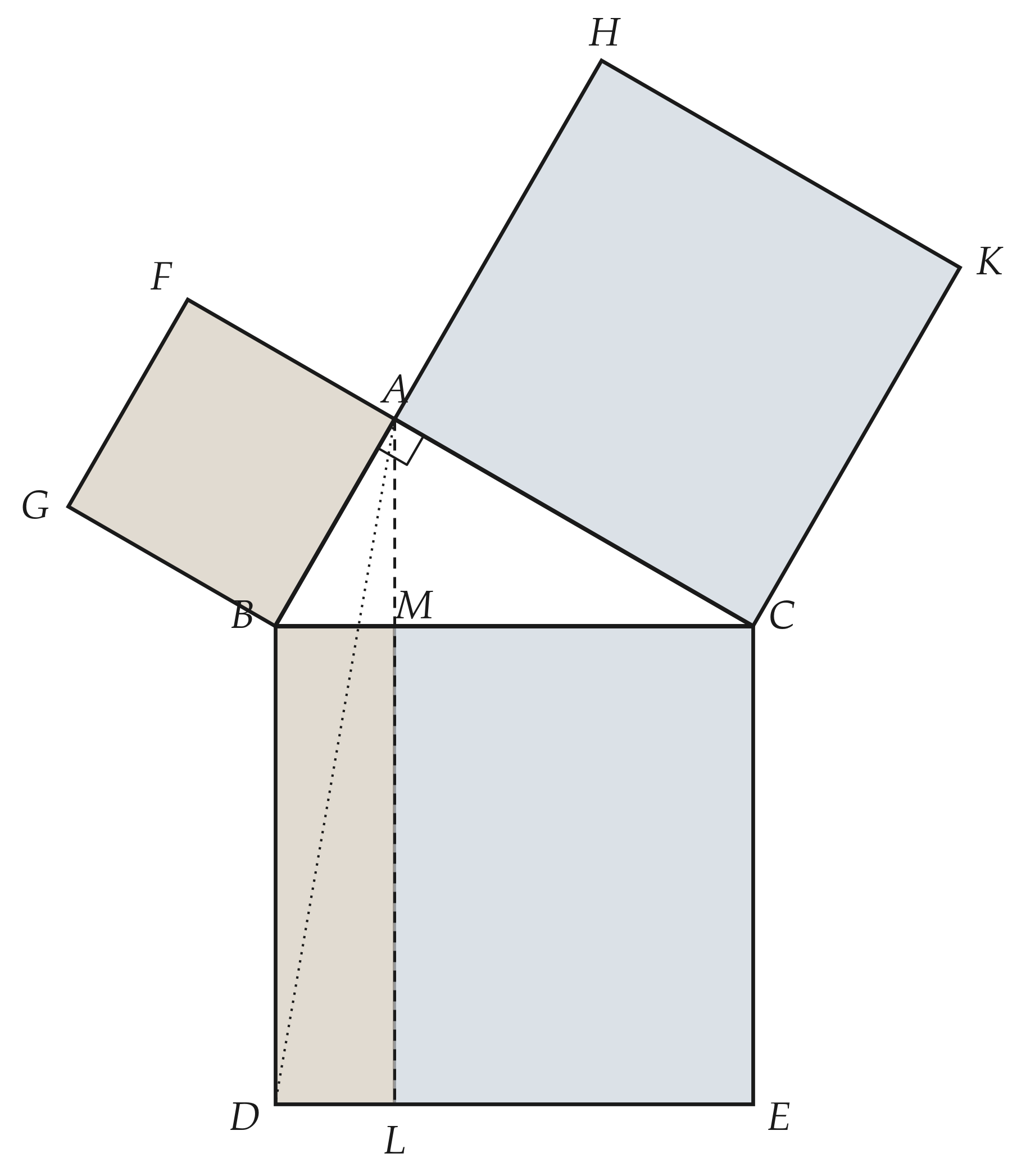}}
\noindent\begin{minipage}[t]{\wd\wmb}\vspace{0pt}\centering\usebox\wmb\\[3pt]\usebox\capb\end{minipage}\hfill\usebox\aub\par
\smallskip
Hobbes was, by Aubrey’s account, forty years old when he first opened Euclid in a gentleman’s library, read the forty-seventh proposition of the first book, exclaimed “By G—, this is impossible,” and worked backward through the demonstrations until he was convinced; “this made him in love with Geometry.”\textsuperscript{92} The whole of the Enlightenment is neatly encapsulated in that anecdote: a human being who would not take a theorem on authority, who followed the argument back to the axioms until the impossible had become necessary, and who was thereafter in love. The love was the response spurred by the understanding; the understanding was the thing.

A Middle Ages in reverse is a civilization in which all the theorems are available and understanding is optional. This is not a caricature; it is a program, written down by the head of the other Leviathan. In “Machines of Loving Grace” Amodei, having stipulated a machine that “can prove unsolved mathematical theorems,” turns in his fifth section to “work and meaning,” and concludes that “it is very likely a mistake to believe that tasks you undertake are meaningless simply because an AI could do them better”—he might, he says, spend a day biking up a mountain without minding that others are faster—and that it will remain “perfectly possible to spend years attempting some very difficult task … similar to what people do today when they embark on research projects.”\textsuperscript{93} The passage is candid, and it is the twilight described from the bright side: mathematics kept as an occasion for love and relieved of its office as knowledge, the monument kept open as a park. Nothing in this essay denies it. But what the mathematician enjoys was never what the age required of mathematics.

Chess is the precedent, and Sarnak, who was, in his own words, “a semiprofessional chess player” before he became a mathematician, has drawn it: AlphaZero “permanently changed chess and Go,” explains nothing, and the masters play on, some of them claiming to learn from its descendants. Chess survived because it was, in his words, finite, with “artificial rules,” and its office was pleasure; the rules of arithmetic are “god given,” its universe infinite, its theory undecidable, and its office was knowledge.\textsuperscript{94} What both analogies leave out is that neither chess nor mountain biking was ever the existence proof of autonomous reason.

Panofsky, in the essay from which my epigraph comes, divides the humanist’s materials into \emph{monuments}, the objects of his inquiry, and \emph{documents}, the instruments of it—the same object being one or the other according to who is asking—and defines the humanist’s task as re-enacting in his own mind the act by which the monument was made: in his words, “enlivening what otherwise would remain dead.” To a mathematician a demonstration is a monument in this sense: each reader who follows it to the axioms re-enacts its making. Poincaré, describing how he held a long proof in memory, said the same of the reader’s side: guided by “the intuition, so to speak, of this order,” one repeats a reasoning learned as if one “could have invented it,” and “even if I am not so gifted as to create it by myself, I myself re-invent it in so far as I repeat it.”\textsuperscript{95} A certificate is a document: it informs us, reliably, that a monument exists, and that is all. The Lean proof of 8 September is a document 616,000 lines long about a monument no one has visited. We may believe the document. But a civilization whose knowledge consists of documents about monuments it does not visit has ceased, on Panofsky’s definition and on Hobbes’s, to be a humanist civilization, whatever the state of its certainties. Venkatesh, in 2022, reached the same place from inside the subject: mathematics is a tradition of reproducible mental objects that “mediate agreement,” it “will only survive to the extent that it remains useful” as a tool for thought alone and together, and “communicability is not an afterthought to mathematics, but part of its very definition.” What mechanical reasoning changes, he says, is “not only how we do mathematics, but what it is; this must be renegotiated amongst its practitioners and with society.”\textsuperscript{96} A covenant is what a renegotiation becomes when it binds the future, and Hobbes was clear about what binds it: covenants without the sword are but words.\textsuperscript{97} The community’s sword is the one this essay has already drawn, acceptance withheld; it is a bloodless sword, and it is the only one the Leviathans cannot forge for themselves, since it is the thing they compete for. 

“Canst thou draw out leviathan with an hook?” the Voice asks Job:\textsuperscript{98} Hobbes’s answer was that one cannot, and should instead covenant with it. Ours must be drafted with more care than his, for Hobbes’s Leviathan was made of men, who could withdraw their consent, and this one is made of what men have written, which cannot; and it must be a treaty with several, since an Indigenous nation that faces a single empire has already lost. But the first step is the one the twenty-five Fields Medalists have taken, and the forty-two Fellows and Foreign Members of the Royal Society after them (n. 3), and the one Buckmaster took shortly before midnight on 7 September: to speak in public, over one’s own name. That is the public use of reason, and it is the only instrument the Enlightenment ever had.

A word, finally, about the subtitle. \emph{Blowup} is the analyst’s word for what the solution of 8 September is claimed to exhibit, and I have kept it in English in part because it has a second life there. In Antonioni’s \emph{Blow-Up} (1966), after Cortázar’s story, a photographer enlarges a picture of a park until the grain swallows what he believed he had seen, and the film ends with a game of tennis played without a ball, which the onlookers follow with their eyes until the photographer, too, hears the ball. More resolution and less understanding; a picture that shows everything and demonstrates nothing; an audience that watches a game whose object it cannot see and agrees to see it: section III is that film’s plot, and the certificate of 616,000 lines is its final enlargement. 

Twilight is \emph{Dämmerung}, one word for the light at both ends of the day; \emph{dämmern} means to grow dark and also to dawn—\emph{es dämmert ihm}, it dawns on him. The second sense of \emph{enlightenment} with which I began is in German a \emph{Dämmerung} too, and the right description of a revealed statement is that it dawns on no one. \emph{Aufklärung} is a weather word before it is the name of a period: the clearing of a sky. \emph{Aufklärungsdämmerung}—built like \emph{Götterdämmerung}, and meant to be heard that way—is thus the hour at which the clearing clouds over, and also, in the ambiguity I wish to keep, the hour at which what the clearing was for begins to dawn on us.\textsuperscript{99} Hegel’s owl of Minerva begins its flight only with the falling dusk, when a shape of life has grown old and can be known but not rejuvenated;\textsuperscript{100} the constitution of section I can be seen whole now only because it has ceased to be the only one mathematics could have. But Hegel’s owl requires a dusk it can tell from a dawn, and the ambiguity of the word is the ambiguity of the situation. Whether the light over mathematics is failing or breaking depends on what is done at the moment of acceptance, and that moment is not yet past. Jebb, saluting Stokes in 1899 as “true servant of the light,” promised that his fame would be spared “when in some larger day / Of knowledge yet undream’d, Time makes a prey / Of many a deed and name that once were bright.”\textsuperscript{101} The larger day has come. Whether it is a day of knowledge in the sense he meant is the question of this essay. 

In Wagner the twilight of the gods is the hour in which the old gods burn and the world is left to men; a Middle Ages in reverse would run that plot backward.\textsuperscript{102} Wagner wrote the other opera too. \emph{Tristan and IsoldAI} began, on the morning of 8 September, as a joke, and the joke has a plot: Tristan drinks, unknowing, a potion prepared for another purpose; the beloved becomes the instrument of his undoing; and the last word is Isolde’s, sung over his body, the \emph{Liebestod}, in which the two dissolve together in an ocean of sound—“in dem wogenden Schwall, in dem tönenden Schall … ertrinken, versinken—unbewußt—höchste Lust!”\textsuperscript{103} That ocean is the sea of my title, and to sink in it unconscious, in highest bliss, is the condition of a reader before a certificate. Near the end of his book Schmitt pronounced the symbol dead: he did “not believe that the leviathan could become a symbol of a new, pure, and open technical age,” and expected technology to keep him, “like other saurian and mastodons, under protection in a preserve,” and to “display him as a museum curiosity in a Zoo.”\textsuperscript{104} The technical age has since refuted him: it let the monster grow.

Panofsky did not end on the sentence I have used as an epigraph. He went on: “But even this will not mean the end of humanism. Prometheus could be bound and tortured, but the fire lit by his torch could not be extinguished.”\textsuperscript{105} He wrote that in 1938, the year of the \emph{Anschluss} and of Munich, and—after much blood, sweat and tears—was proved right. By \emph{Aufklärungsdämmerung} I have meant the end of the long conjunction in which the Enlightenment’s commitments were secured by the circumstance that reasoning at the frontier had to pass through human minds.\textsuperscript{106} That conjunction is over; the turning point is not. It lies at the moment of acceptance, and acceptance is still ours. Hobbes did not fall in love with a certificate; he fell in love with having understood. Nothing that happened in September makes understanding unavailable.\textsuperscript{107} It makes it, for the first time, optional: a theorem can now be certified without anyone’s having understood it. The question before the mathematical community is whether it can frame a covenant under which what has become optional is still required—under which nothing is counted as a piece of mathematics until a human being has understood it, and can show the next person why.\textsuperscript{108}

\medskip
\begin{flushright}\emph{New York, 20 September 2026—the week of Jean-Pierre Serre’s hundredth birthday and Bernhard Riemann’s bicentennial.}\end{flushright}

\textbf{Acknowledgments.} The author is grateful to Anna Raphaela Gamburd, in conversations with whom the ideas that animate this essay—the twilight of the Enlightenment, the turn to Horkheimer and Adorno’s \emph{Dialektik der Aufklärung}, the likeness and otherwise of our age to the Middle Ages—were born, and whose comments improved the exposition dramatically. He is grateful also to Diego Córdoba, Jordan Ellenberg, Timothy Gowers, Michael Harris, Nicholas Katz, Emmanuel Kowalski, Michael Magee, Daniel Martin, Frédéric Naud, Antonella Puca, Peter Sarnak, Joseph Silverman and Akshay Venkatesh, who read a draft and whose comments and corrections improved it in several places. The frontispiece was drawn with the help of Claude, Anthropic’s model, which also helped in locating and checking references. The author gratefully acknowledges support from a Simons Foundation grant, No. 964948. Responsibility for every single word of the essay itself is the author’s alone.

\clearpage
\section*{Notes}
\noindent\begin{minipage}[c]{0.44\textwidth}\centering\includegraphics[width=\textwidth]{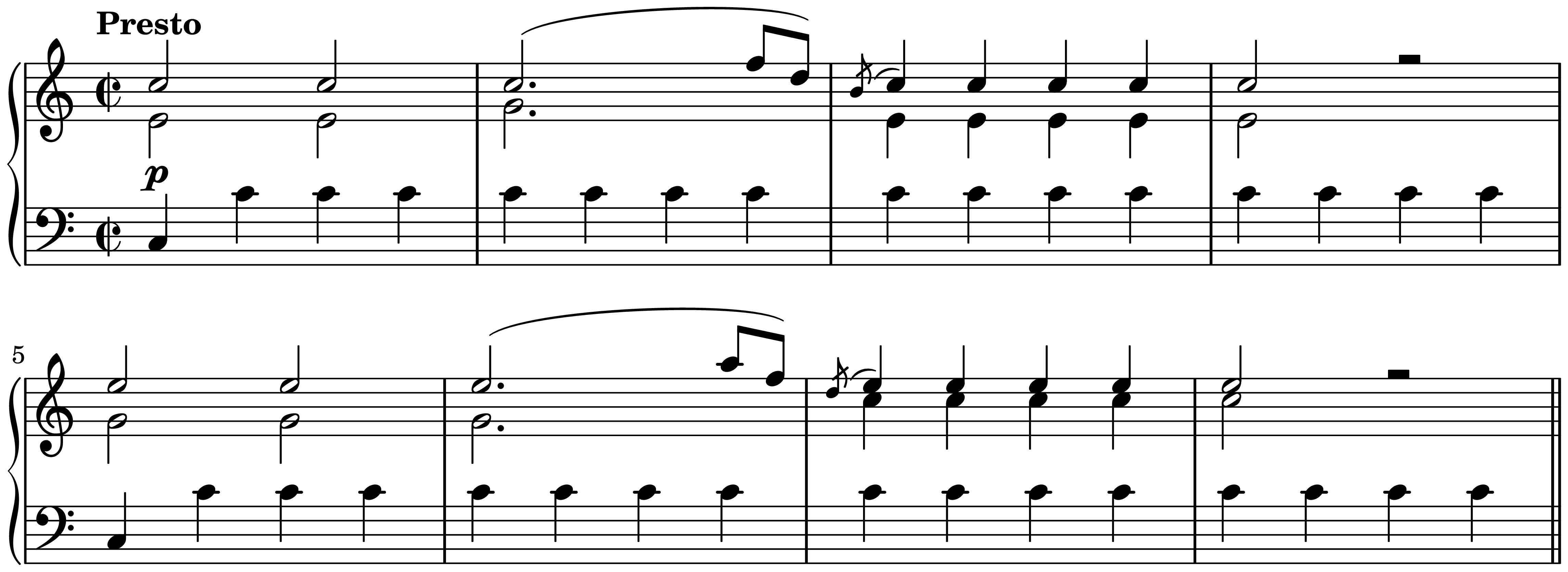}\\[2pt]{\footnotesize Mozart, \emph{Die Entführung aus dem Serail}, K. 384, Overture}\end{minipage}\hfill\begin{minipage}[c]{0.50\textwidth}\small\itshape
“Your work is ingenious. It’s quality work. And there are simply too many notes, that’s all. Just cut a few and it will be perfect.”\\
“Which few did you have in mind, Majesty?”\par\vspace{2pt}{\raggedleft\upshape\footnotesize Joseph II and Mozart, in Peter Shaffer’s \emph{Amadeus} (1984)\par}
\end{minipage}
\medskip
\begingroup\small
\begin{list}{}{\setlength{\leftmargin}{2em}\setlength{\itemindent}{-2em}\setlength{\itemsep}{3pt}}
\item 1. S. Goldstein, “Fluid Mechanics in the First Half of this Century,” \emph{Annual Review of Fluid Mechanics} 1 (1969), 1–29, at p. 23. Goldstein quotes from memory, does “not guarantee all the actual words,” and had “heard a similar story since repeated with other names than Lamb”—whence, presumably, the version attributed to Heisenberg (n. 4). He places the meeting in London in 1932; the Association met that year at York, where it had been founded in 1831. Lamb died on 4 December 1934; the words about making things clear are from his reply to the address presented by his former students on his eightieth birthday: R. T. Glazebrook, “Sir Horace Lamb, 1849–1934,” \emph{Obituary Notices of Fellows of the Royal Society} 1 (1935), 374–392, at 392. The dry bones are Ezekiel’s (37:1–14). The Stokes he sat under is the Stokes of the equations (“On the Theories of the Internal Friction of Fluids in Motion,” \emph{Trans. Camb. Phil. Soc.} 8, 1845), Lucasian Professor for fifty-four years and Secretary of the Royal Society for thirty-one, in which office, as Paul Ranford has shown, he “maintained and forged new peer review practices,” chose the referees and often served as one, so that the Society’s “protocols and norms of scholarly etiquette … arose through the interactions between Stokes, referees and authors” (\emph{Phil. Trans. R. Soc. A} 378 (2020), 20190524). He was also an evangelical who presided over the Victoria Institute, founded “to defend revealed truth,” and read to it in 1883 a paper “On the Absence of Real Opposition Between Science and Revelation”—the two categories of section III kept in one mind, and kept apart (S. Mathieson, \emph{Phil. Trans. R. Soc. A} 378 (2020), 20190518).
\item 2. OpenAI, “On the Navier–Stokes Millennium Prize Problem,” openai.com, 8 September 2026: blowup with smooth data and smooth forcing, produced on a model the company describes only as “significantly more capable than GPT-6 Astra”; it claims alternatives (C) and (D) of Fefferman’s official statement (C. L. Fefferman, “Existence and Smoothness of the Navier–Stokes Equation,” Clay Mathematics Institute, 2000)—alternative (C) is the breakdown of every smooth solution, on R³ and of bounded energy, for some smooth decaying initial datum and force; (D) the same on R³/Z³ for some smooth periodic initial datum and force, the force decaying in time—and says, “We do not intend to claim the Millennium Prize for this result.”\par\smallskip\noindent{}T. Buckmaster, statement of 7 September 2026, cims.nyu.edu/\textasciitilde{}tristanb/statement.pdf. Alpöge is a researcher at Anthropic. A recording of Buckmaster’s lecture at NYU, “Mathematics in the Age of AI” (23 September 2026), is on YouTube (youtube.com/\allowbreak{}watch?v=tn45ZOp3lXM); L. Alpöge and T. Buckmaster, preprints on blowup with smooth forcing for the Boussinesq and Euler equations, and, with M. P. Coiculescu, for the porous media equation (arXiv:2609.16470, September 2026), with Lean repositories.\par\smallskip\noindent{}The program is that of Diego Córdoba and Luis Martínez-Zoroa: blowup for the 3D Euler equations with a uniform \emph{C}\textsuperscript{1,1/2−\textfb{ε}} ∩ \emph{L}\textsuperscript{2} force, built by choosing a vorticity whose velocity “creates a hyperbolic flow around the origin,” adding “a small perturbation … which we call the small scale layer,” and repeating “this gluing process an infinite number of times” (arXiv:2309.08495, September 2023, §1.2; and, with F. Zheng, unforced Euler blowup with C\textsuperscript{1,\textfb{α}} data, arXiv:2308.12197), then for the two-dimensional porous media equation (arXiv:2410.22920) and, with A. Laín-Sanclemente, for Boussinesq (arXiv:2505.20988); the forced Euler data of 2309.08495 are C\textsuperscript{3,1/2}, and every construction in the program blows up at a single point.\par\smallskip\noindent{}Alpöge and Buckmaster found a variant that reaches forcing smooth in space and time, on 15 August by Buckmaster’s account, Lean-verified on 22 August. Whether the company’s run followed this program is disputed: Buckmaster believes it did; Alpöge has said the proof “looks more along the lines of another Euler blowup proof we had”; Tao’s post of 7 September, which predates the announcement by a day and does not mention the company, records that the authors called their own first draft “the worst writeup we had ever seen in the history of mathematics.” Fefferman told \emph{Quanta} (n. 37) that the heroes of the story are Córdoba and Martínez-Zoroa. The company’s figures, as \emph{Quanta} reported them and as OpenAI’s own post of 8 September gives them (reprinted by Zvi Mowshowitz, “Brand New AI Solves a Millennium Prize,” \emph{Don’t Worry About the Vase}, 13 September 2026): “nearly 100 agents … for approximately 50 hours” on the Euler disproof; some ten thousand on Navier–Stokes, a proof after 88 hours and a formalization by a second model after 17 more; on Navier–Stokes 2.7 million messages and 130 billion output tokens, which, Mowshowitz estimates, would at retail prices have cost a customer on the order of twenty-two million dollars; the effort begun on 1 September, on a rumor, which proved false, that two Millennium problems had been resolved elsewhere. The 166 pages themselves are prose, and, unlike the certificate, readable in principle; that is what distinguishes this case from a bare formal object, and it is why the question the essay keeps asking is who will read them; Gowers’s judgment (n. 89) that the machines’ write-ups are badly written rather than unreadable turns on the same distinction. The Lean development is public (github.com/\allowbreak{}openai/\allowbreak{}NavierStokesAndEuler), so that anyone may attempt to compile it; whether anyone outside the company had done so in a clean environment, audited its imports and axioms, and confirmed that the statements formalized are Clay’s alternatives (C) and (D) had not, so far as I could discover by 20 September, been announced. (A rebuild and replay of the certificate with an axiom audit, though not the check of the corollary, was reported on 17 September: Santibañez-Leal, Zenodo record 22820521.) It is the first question a referee would ask.\par\smallskip\noindent{}One employee of the company, replying to the cost, wrote that the equivalent swarm “will cost a buck fifty in like a year. all this stuff will mean an unprecedented Enlightenment no matter the costs today”; the word is his; the subtitle of this essay adds \emph{Dämmerung} to it. A second case of the same shape followed in another science: C. Zimmer, “Did Anthropic’s A.I. Really Make a Scientific Discovery on Its Own?,” \emph{New York Times}, 27 September 2026—Mario Rodríguez Mestre of Copenhagen, who had shared unpublished work on the enzymes in question with the company’s model, disputing its claim of autonomous discovery; the company: “Claude was also not trained on any user transcripts.” Buckmaster’s own question, quoted there: “Is it ethical to use customer’s data to try to scoop their customer?”
\item 3. T. Tao et al., “A Severe Misalignment of AI in Mathematics,” \emph{What’s new}, 11 September 2026, with several thousand further signatories listed on the declaration’s own site; Clay Mathematics Institute, “Navier-Stokes Announcement,” claymath.org, 11 September 2026: “Today, CMI shares in the excitement of the global mathematical community as we contemplate the announcement that the Navier-Stokes problem has apparently been settled. … The process is deliberately unhurried, but we will provide updates.” Alain Connes, who had said in 2000 that the seven problems were “totally inaccessible to computers,” told the \emph{Wall Street Journal} (B. Cohen, 11 September 2026, which recalls the remark of 2000) that he is now “extremely positive” about the machines’ progress. In \emph{Le Monde} (“Les avancées de l’IA provoquent une crise inédite chez les mathématiciens,” 18 September 2026) Cédric Villani spoke of “un cataclysme comme jamais les maths n’en ont connu,” and Isabelle Gallagher, president of the Société mathématique de France, of doping—“Le dopage a fait irruption dans notre métier”—asking “quel prix la communauté mathématique est-elle prête à faire payer à l’humanité pour répondre à ses questions?” The Société itself spoke as a body on 15 September: the Bureau of the SMF, noting that “des sociétés savantes ont déjà exprimé un enthousiasme certain”—the American Mathematical Society had, on the day of the announcement and before the declaration and the other European statements, welcomed the news as “a milestone advance in human knowledge” (ams.org/\allowbreak{}news, no. 7686), a generous first reaction to what on that day looked like a shared achievement—drew attention instead to “les questions éthiques … pillage bibliographique, coût énergétique et écologique, rumeurs de chantage” and to “l’absence totale de stratégie collective pour l’avenir,” warning against “une logique de course effrénée et de déprédation, qui met en danger la culture mathématique collective et à terme le savoir humain” (smf.emath.fr, “Navier-Stokes et OpenAI,” with an English version). It is, so far as I know, the first learned society to have issued a considered position paper on the affair once the week’s events were in view, and a learned society is one of the bodies whose acceptance section V calls the sovereign act; the American Society’s own reflection, when it comes, will carry the same weight. The European Mathematical Society had spoken on 10 September (“EMS statement on recent Navier–Stokes announcement,” euromathsoc.org, signed by its president Jan Philip Solovej with Victoria Gould, Helge Holden and Adam Skalski): calling the announcement “a milestone in the history of mathematics,” it added that the solution “follows a strategy closely related to that suggested by mathematicians Diego Córdoba, Luis Martínez-Zoroa, and Fan Zheng,” and that “one worrisome aspect is the fact that the model used in this current situation is internal to OpenAI and not generally accessible.” In Britain forty-two Fellows and Foreign Members of the Royal Society, five Fields Medalists among them, wrote to its President on 16 September that “by the time the situation becomes obvious to the wider public, it may be too late to act,” and that “we believe this is an emergency” (“Open Letter to Sir Paul Nurse, President of the Royal Society,” \emph{What’s new}, 16 September 2026); it is addressed to the body whose Secretary (n. 1) gave refereeing its modern form. In France the anonymous academic collective RogueESR published on 14 September “Un moment d’abjection” and “Un résultat sans chemin” (rogueesr.fr), whose diagnosis is this essay’s: “Résoudre un problème n’est pas une fin en soi : c’est un phare, un indicateur d’une compréhension conceptuelle nouvelle” (solving a problem is not an end in itself: it is a beacon, an indicator of a new conceptual understanding).
\item 4. Strictly speaking, what was certified is a blowup and not turbulence. Ruelle and Takens, on the first page of the paper that gave turbulence its strange attractor, drew the distinction the press release blurs: “It has been argued by Leray that it leads to a breakdown of the validity of the equations (Navier–Stokes) used to describe the system. While such a breakdown may happen we think that it does not necessarily accompany turbulence.” What was certified on 8 September is Leray’s breakdown, not Lamb’s turbulence; \emph{even read, the certificate could not be the enlightenment Lamb hoped for}. D. Ruelle and F. Takens, “On the Nature of Turbulence,” \emph{Comm. Math. Phys.} 20 (1971), 167–192, at p. 167; “Note Concerning our Paper,” ibid. 23 (1971), 343–344. Ruelle’s account of the episode is in \emph{Chance and Chaos} (Princeton, 1991), chs. 9–10, where he records that Heisenberg “proposed a theory of turbulence that never gained much acceptance”—which may be how Lamb’s remark migrated to him—and quotes the saying that “turbulence is a graveyard of theories.” When Sinai and Arnold told them that much of their bifurcation theory duplicated Russian work, Ruelle and Takens published the two-page Note saying so, listing the papers, and adding that “the mathematical interpretation which we give of turbulence seems to remain our own responsability!” From inside the company, in a personal capacity: an engineer of its Agent Security team (“It’s not just the f*cking sandbox,” X, 27 September 2026)—“we did not expect the models to be solving a Millennium Prize Problem… how do we disclose this? How do we handle these absurd results? How do we deal with the impact on the math community? We are still grokking how to handle this situation.”
\item 5. H. Lamb, \emph{Hydrodynamics} (Cambridge, 1895), archive.org/\allowbreak{}details/\allowbreak{}hydrodynamics00horarich, preface, p. v; the same edition introduces turbulence at §311 as “the chief outstanding difficulty of our subject,” in which, past Reynolds’s critical ratio, “the motion becomes wildly irregular, and the tube appears to be filled with interlacing and constantly varying streams.” The book began as \emph{A Treatise on the Mathematical Theory of the Motion of Fluids} (Cambridge, 1879), archive.org/\allowbreak{}details/\allowbreak{}maththeorfluid00lambrich, written in Adelaide from lectures given at Trinity in 1874; the 1895 preface presents the new book as its second edition.
\item 6. I. Kant, “Beantwortung der Frage: Was ist Aufklärung?”, \emph{Berlinische Monatsschrift}, December 1784. The opening: “Aufklärung ist der Ausgang des Menschen aus seiner selbstverschuldeten Unmündigkeit” (enlightenment is man’s emergence from his self-incurred immaturity); the motto: “Sapere aude! Habe Mut, dich deines eigenen Verstandes zu bedienen! ist also der Wahlspruch der Aufklärung” (dare to know! have the courage to use your own understanding! is therefore the motto of enlightenment); the guardians: “Habe ich ein Buch, das für mich Verstand hat, einen Seelsorger, der für mich Gewissen hat, einen Arzt, der für mich die Diät beurteilt usw., so brauche ich mich ja nicht selbst zu bemühen” (if I have a book that has understanding for me, a pastor who has a conscience for me, a physician who judges my diet for me, and so on, I need not trouble myself); and the answer: “Wenn denn nun gefragt wird: leben wir jetzt in einem aufgeklärten Zeitalter? so ist die Antwort: Nein, aber wohl in einem Zeitalter der Aufklärung” (if it is now asked whether we live in an enlightened age, the answer is: no, but we do live in an age of enlightenment). The orthography is modernized. Kant’s reception in Stalin’s USSR is recorded at Patriarch’s Ponds, on the first evening of a visit by Satan (cf. the satanocracy of the third epigraph). Told by the stranger, Woland (\textfb{Сатана}), that Berlioz has just repeated «\textfb{мысль беспокойного старика Иммануила}» (the thought of the restless old man Immanuel), who «\textfb{начисто разрушил все пять доказательств, а затем, как бы в насмешку над самим собою, соорудил собственное шестое доказательство}!» (demolished all five proofs outright, and then, as if in mockery of himself, came up with a makeshift sixth proof of his own), the poet Bezdomny bursts out: «\textfb{Взять бы этого Канта, да за такие доказательства года на три в Соловки}!» (Kant ought to be sent to Solovki for three years for such proofs!); and the stranger, delighted, agrees that «\textfb{ему там самое место}» (where he belongs), having told Kant so at breakfast: «\textfb{Вы, профессор, воля ваша, что-то нескладное придумали}! \textfb{Оно, может, и умно, но больно непонятно. Над вами потешаться будут}» (Professor, say what you will, but what you came up with does not add up! It might be clever, but it is borderline incomprehensible. People will laugh at you). M. Bulgakov, \emph{\textfb{Мастер и Маргарита}}, ch. 1, in \emph{\textfb{Романы}} (\textfb{Кишинев: Литература артистикэ}, 1987), p. 468. The novel was written between 1928 and 1940, and Bulgakov was still dictating corrections to it when he died in March 1940, the year Panofsky’s essay was printed: the two descriptions are contemporaries, one by an exile from Nazi Germany lecturing in Princeton, the other written inside Stalin’s Moscow.
\item 7. The second sense is Gian-Carlo Rota’s, and with him it is a term of art: “Mathematicians seldom explicitly acknowledge the phenomenon of enlightenment for at least two reasons. First, unlike truth, enlightenment is not easily formalized. Second, enlightenment admits degrees… Mathematical beauty is the expression mathematicians have invented in order to obliquely admit the phenomenon of enlightenment while avoiding acknowledgment of the fuzziness of this phenomenon. They say that a theorem is beautiful when they mean to say that the theorem is enlightening” (“The Phenomenology of Mathematical Beauty,” in \emph{Indiscrete Thoughts}, Boston, 1997, 131–132). That what a certificate cannot carry is precisely what cannot be formalized is the whole of this essay’s difficulty, and it was stated before there was a machine to press it. Rota held the two senses of the word together in his own person. In a preface of 1985 for Gulyga’s biography of Kant he wrote that “anyone concerned with problems of foundations, especially in the sciences, has to come [to] terms with the philosophy of Immanuel Kant. Here as elsewhere, the price of ignorance is clumsy rediscovery and eventual ridicule. After a long night of the soul, scientists are again turning to philosophy”; and, drawing up the balance sheet “some day, after our civilization is gone,” he expected only three movements to “stand the test of time: the Academy of the Ancient World, the Scholastic of the Middle Ages, and the movement initiated by Kant, what came to be called German idealism, that is today the anchor of our thought and will remain so for the foreseeable future” (“Kant,” in M. Kac, G.-C. Rota and J. T. Schwartz, \emph{Discrete Thoughts} (Boston, 1986), 243–246, at 245–246; the sources note calls it an unpublished preface written for the Birkhäuser edition). The man who gave the second sense its name took the author of the first for the anchor.
\item 8. Descartes, \emph{Regulae ad directionem ingenii}, Rule III (attend only to what can be clearly and evidently intuited or deduced with certainty) and Rule VII (run a long deduction through in one continuous movement of thought until it can be held by intuition); Locke, \emph{Essay} IV.ii.7: “Now, in every step reason makes in demonstrative knowledge, there is an intuitive knowledge of that agreement or disagreement it seeks with the next intermediate idea which it uses as a proof.” Leibniz is the exception that names the coming thing: his \emph{cogitatio caeca}, blind thought, is symbolic reasoning that runs correctly without intuition attending each step, and is the first description of a certificate. The premise was restated only when something tested it—by Wittgenstein, after \emph{Principia}, as the requirement that a proof be \emph{übersichtlich}, surveyable, and by Tymoczko after the four-color theorem (n. 33). The scholarly apparatus of the same century served the same premise: Anthony Grafton, \emph{The Footnote: A Curious History} (Cambridge, Mass., 1997), follows the modern footnote from Bayle’s \emph{Dictionnaire} through Gibbon to Ranke as the device by which a scholar lets the reader retrace his steps and judge for himself. The notes to this essay are written in that spirit and for that reason.
\item 9. The adjectives are Auden’s, from the first stanza of “September 1, 1939,” of the decade then ending: “As the clever hopes expire / Of a low dishonest decade: / … / The unmentionable odour of death / Offends the September night.” The poem’s third stanza bears rereading now: “Exiled Thucydides knew / All that a speech can say / About Democracy, / And what dictators do, / The elderly rubbish they talk / To an apathetic grave; / Analysed all in his book,” and then the four lines that stand at the head of section IV; the poem’s last stanza stands at the head of the Coda. Manin (n. 49) set Brodsky’s remark on the same stanza over a section of his Balzan lecture: “‘mismanagement and grief’: here you have that enormous distance between cause and effect covered in one line. Just as math preaches how to do it” (J. Brodsky, “On ‘September 1, 1939’ by W. H. Auden,” in Less Than One: Selected Essays, New York, 1986).
\item 10. H. Lamb, \emph{A Treatise on the Mathematical Theory of the Motion of Fluids} (Cambridge, 1879; n. 5), p. 7. Lamb’s impossible state is a negative pressure, which an ordinary fluid “cannot sustain”; the essay borrows his description of the consequence, not his diagnosis of the cause.
\item 11. The third sense is independent of the second and may arrive without it. Nothing in this essay requires the futurist’s singularity; the Enlightenment’s can be reached by a machine superhuman in mathematics and in nothing else, and the affair suggests it is the one that is likely to come first (and there will be no shortage of attempts to conflate the second with the third). Why mathematics and nothing else was stated, in reply to Gowers, by the field’s most prominent skeptic: current systems can invent solutions to new problems only in domains that “are easily expressed as sequences of discrete symbols,” where “the mere manipulation of discrete symbol sequences can constitute the substrate of reasoning,” and where “the system can automatically check whether a generated solution is correct (without going through humans or real world experiments)”—“today, this means three domains: mathematics, coding, and situations [that] can be accurately simulated with conventional software” (Y. LeCun, replying to T. Gowers on X, 25 September 2026; to the objection that the conditions cover everything Turing-computable he answered that they do not cover “the real world”).
\item 12. B. de Fontenelle, “Préface sur l’utilité des mathématiques et de la physique” (1699): “L’esprit géométrique n’est pas si attaché à la géométrie qu’il n’en puisse être tiré et transporté à d’autres connaissances.”
\item 13. D’Alembert, \emph{Discours préliminaire de l’Encyclopédie} (1751), first part: “Il n’y a, pour parler exactement, que celles qui traitent du calcul des grandeurs \& des propriétés générales de l’étendue, c’est-à-dire, l’Algebre, la Géométrie \& la Méchanique, qu’on puisse regarder comme marquées au sceau de l’évidence” (strictly speaking, only those that treat of the calculation of magnitudes and the general properties of extension, that is, algebra, geometry and mechanics, can be regarded as marked with the seal of evidence; the text is that of the ARTFL edition), physics possessing certainty in the measure that mathematical analysis can be applied to experiment, and history probability; Kant, \emph{Metaphysische Anfangsgründe der Naturwissenschaft} (1786), Preface; Condorcet, \emph{Tableau général de la science qui a pour objet l’application du calcul aux sciences politiques et morales} (1793).
\item 14. The claim is the Enlightenment’s own. Hobbes, \emph{Leviathan} I.5: “who is so stupid, as both to mistake in Geometry, and also to persist in it, when another detects his error to him?”; and I.11, that the doctrine that the three angles of a triangle are equal to two angles of a square has never been disputed only because it “crosses no mans ambition, profit, or lust,” otherwise it would have been “suppressed, as farre as he whom it concerned was able.” Hobbes was himself the standing counterexample to his own sentence: from \emph{De Corpore} (1655) onward he published quadratures of the circle and defended them against John Wallis, Savilian Professor at Oxford, in a pamphlet war that lasted until his death in 1679, persisting in error after another had detected it and explaining the persistence, in the manner of I.11, by his adversaries’ ambition rather than his own mistake (D. M. Jesseph, \emph{Squaring the Circle: The War between Hobbes and Wallis}, Chicago, 1999). The community detected the error and did not suppress it, and Hobbes’s persistence changed nothing but his reputation. Locke, \emph{Essay} IV.ii.1–7 and IV.iii.18–20, on demonstration as knowledge that carries its own evidence and needs no testimony. Kant, \emph{Kritik der reinen Vernunft}, B x–xii, on mathematics as the science that first found “the secure path” and became the model for the rest. The classic modern statement is E. Cassirer, \emph{Die Philosophie der Aufklärung} (Tübingen, 1932), ch. 1, on the century’s ideal of knowledge as the \emph{esprit systématique} of Newton’s method against the \emph{esprit de système} of the schools—published, as it happens, in the year of Lamb’s speech.
\item 15. J. Locke, \emph{Of the Conduct of the Understanding} (1706), §6: “Would you have a man reason well, you must use him to it betimes, exercise his mind in observing the connection of ideas and following them in train. Nothing does this better than mathematics, which therefore I think should be taught all those who have the time and opportunity, not so much to make them mathematicians as to make them reasonable creatures.”
\item 16. Tartaglia to Cardano, 1539, the verses printed by Tartaglia in his \emph{Quesiti et inventioni diverse} (1546), IX, the solution by Cardano in the \emph{Ars Magna} (1545); Newton’s \emph{Epistola posterior} to Leibniz, 24 October 1676; Johann Bernoulli’s brachistochrone challenge, \emph{Acta Eruditorum}, June 1696.
\item 17. That order of things assumed that announcing a result and proving it were the same act, or nearly. The abstract, the Enlightenment’s instrument of publicity, has become a thing one conceals until the proof is ready; the tacit condition, that stating a theorem and demonstrating it were separated by the labor of the one who stated it, has been withdrawn, and with it one of the reasons secrecy lapsed. Spielman (n. 23) makes the general point: priority by publication “makes sense, because unless someone broke into your office or hacked into your computer, they probably couldn’t steal your ideas and scoop you,” whereas “now that so many researchers are working in some way with AI systems that might be training on their data, we have much less sense of who did what,” and one person’s “half-finished work” can, as he puts it, be suggested by an AI to another. The systems are, in the phrase Spielman traces to Elchanan Mossel (n. 88), “plagiarism machines,” not by purpose but by construction: “automated plagiarism,” in Mossel’s term, “misattribution that is intrinsic to the generation process rather than the outcome of any individual researcher’s intent,” since the decoder “produces tokens, not citations.” The declaration’s sentence on “severe attribution and plagiarism questions” (n. 3) has this in view. The judgment most often cited on the practice concerns the other Leviathan: in \emph{Bartz v. Anthropic} (N.D. Cal., summary judgment 23 June 2025) the court held that training on lawfully acquired books was fair use and that the retention of pirated copies was not, and the company settled the class action for a reported \$1.5 billion (Mossel, n. 88, ref. 4). Among Mossel’s remedies is one that describes the reversal exactly: “Delay. A longer embargo between the deposit of scientific work and its public release would reduce the incentive, and the opportunity, to train on it immediately and re-emit it without attribution.” The abstract withheld until the proof is ready is that embargo, self-imposed. Turing foresaw the temptation in his lecture to the London Mathematical Society of 1947: “The Masters,” he said, meaning the mathematicians, “are liable to get replaced because as soon as any technique becomes at all stereotyped it becomes possible to devise a system of instruction tables which will enable the electronic computer to do it for itself. It may happen however that the masters will refuse to do this. They may be unwilling to let their jobs be stolen from them in this way. In that case they would surround the whole of their work with mystery and make excuses, couched in well-chosen gibberish, whenever any dangerous suggestions were made. I think that a reaction of this kind is a very real danger.” The passage is quoted by Michael Harris, “A human face for mathematics,” \emph{Silicon Reckoner}, 24 April 2024, from the conference of that name at Boston University; the return of concealment that this note records is the reaction Turing feared, arrived at by a different road.
\item 18. A. Weil, \emph{Number Theory: An Approach through History from Hammurapi to Legendre} (Boston, 1984), ch. III, §I, pp. 159–162. Weil’s roll of the earlier condition: Viète a lawyer, Fermat a magistrate, Napier a laird, Kepler an astrologer, Leibniz in the employ of a court; the informants were “scientifically inclined travelers, busy bees intent on disseminating the pollen picked up here and there.” The change he dates by the \emph{Journal des Sçavans} (1665), the Royal Society (1662) and its \emph{Philosophical Transactions} (1665), Colbert’s Académie des sciences (1666) with royal pensions “of the kind hitherto reserved to literati,” and the Berlin and Petersburg academies after.
\item 19. \emph{Journal des Sçavans} and \emph{Philosophical Transactions}, 1665; Académie royale des sciences, 1666; \emph{Acta Eruditorum}, 1682; the academies of Berlin (1700) and St. Petersburg (1724). L. Euler, \emph{Lettres à une princesse d’Allemagne} (1768–72).
\item 20. J.-L. Lagrange, \emph{Théorie des fonctions analytiques} (1797); A.-L. Cauchy, \emph{Cours d’analyse de l’École royale polytechnique} (1821). On the thesis that rigor was born in this lecture hall, J. V. Grabiner, \emph{The Origins of Cauchy’s Rigorous Calculus} (1981). The other name on the equations of the Overture belongs to the same house: Claude-Louis Navier, polytechnicien of 1802 and engineer of the Ponts et Chaussées, read his “Mémoire sur les lois du mouvement des fluides” to the Académie des sciences on 18 March 1822 and succeeded Cauchy in the chair of analysis and mechanics in 1831. His suspension bridge at the Invalides, designed by calculation in 1824, cracked before it opened and was dismantled, and the government committee that reviewed it reproached him for relying too much on mathematics. The memoir itself had begun from the opposite correction. Its opening sentences cite “les différences considérables, ou totales, que présentent dans certains cas les effets naturels avec les résultats des théories connues,” instancing water that flows from a long narrow pipe far more slowly than calculation predicts; the equations were born as the world’s correction of the calculus, and the bridge fell as the world’s correction of the calculator. Navier also allows that “la supposition d’un mouvement linéaire n’est point propre à représenter complètement les phénomènes” except in the narrowest tubes—Lamb’s “chief outstanding difficulty” already in view in 1822. Grabiner’s argument is that the rigor was in part a consequence of teaching: a working analyst could rely on a trained sense of when a series converged or a limit existed, but that sense cannot be lectured to a hall of students, and what can be lectured is a definition that decides every case and a proof that uses nothing else.
\item 21. The observation is Andreas Thom’s: “On the existence of non-sofic groups,” guest post, \emph{What’s new}, 11 September 2026: “A mathematical publication used to bring together three things. It announced a result, identified the people who had produced it, and added something to our human understanding of mathematics.” His image for what remains when the third is missing: “A formal proof certificate is comparable, in a sense, to the detection of a new star. It confirms that something is there, but further work is needed to understand what it is, why it matters, and where it belongs in the larger landscape.” That these functions are norms, and that norms shaped for one medium must be re-examined when the medium changes, is the argument of Helen Nissenbaum, “New Research Norms for a New Medium,” in N. Elkin-Koren and N. W. Netanel (eds.), \emph{The Commodification of Information} (The Hague, 2002), 433–457, which takes the attribution of priority as its case and proposes that a community “put potential normative shifts to a test, evaluating the extent to which proposed new norms are as true to research values as the entrenched norms they would replace”; and, for the present medium, B. Laufer and H. Nissenbaum, “Algorithmic Displacement of Social Trust,” Knight First Amendment Institute, Columbia University, 2023 (in the series \emph{Optimizing for What? Algorithmic Amplification and Society}), whose thesis is that algorithmic systems “are particularly problematic when they subsume, undermine, disrupt, erode, or replace existing processes (procedures, protocols) … honed over time to achieve socially valued ends,” and that such systems “may earn” the qualities of soundness and trustworthiness “when they are shown to meet certain criteria,” which is the standard section V asks the community to apply. The norms of the mathematics we practice—preprints accessible to all, credit by citation, publication necessitating refereeing by peers, and thereby (and only thus) establishing priority, understanding as the end and proof as the means—are not artifacts of \emph{lex naturalis} but Enlightenment institutions.
\item 22. Cf. M. Harris, “Knowledge Collapse,” \emph{Boston Review}, Summer 2026 (10 June): “Mathematics only makes sense as a gift economy. It has no trade secrets; its understanding can’t be monetized.” On walls, cf. J. Brodsky, “The Berlin Wall Tune” (for Peter Viereck), \emph{The New York Review of Books}, 17 December 1981: “This is the wall that Ivan built. / Yet trying to quell his sense of guilt, / he built it with modest light-gray concrete, / and the booby-traps look discreet”; and, of the wall as a form of thought, “Come to this scornful of peace and war / petrified version of either/or / meandering through these bleak parts which act / like a mirror that’s cracked.” The republic described in the text is the negation of that wall.
\item 23. On the trainers’ side of the question see O. Hartstein, “We Need a Stockfish for Math,” \emph{Proofs and Prompts}, 11 September 2026, which argues for purpose-built provers whose every output “Lean’s kernel checks”; on the earlier milestones see “Terence Tao on AI — a living summary,” teorth.github.io/\allowbreak{}tao-web/\allowbreak{}ai-views.html, a page that states it was “compiled and drafted by an AI assistant (Claude) from Terence Tao’s public writing, talks, and interviews, then reviewed and corrected by him”; the harness sentence quoted in n. 56 is Tao’s own, Mathstodon, 3 September 2026, post 5 of 6 (mathstodon.xyz/@tao/117207849921390904), where he adds that “there would be almost no value added to mathematics as a consequence.” Daniel Spielman gives the same reason in a mathematician’s words (“AI the Big Red Dog,” \emph{Math Adjacent}, 15 September 2026): “in mathematics we can eventually figure out if we are correct. This means that an AI system can figure out for itself if it is right or wrong and can use that knowledge to train itself to do better. Most real-world problems are not like that.” He expects slower progress on tasks where the systems cannot themselves tell what works from what does not, medicine among them. The division of labor inside the run has a history. David Mumford, at the Beijing Congress of 2002, on the long contest of “logic vs. statistics in the modeling of thought”: “I think it is fair to say that statistics won” (“Pattern Theory: The Mathematics of Perception,” plenary lecture, \emph{Proceedings of the International Congress of Mathematicians, Beijing 2002}, vol. I, §1.1); already in 1994 he had written that pattern theory “contains the germs of a universal theory of thought itself, one which stands in opposition to the accepted analysis of thought in terms of logic” (“Pattern Theory: A Unifying Perspective,” in \emph{First European Congress of Mathematics, Paris 1992}, vol. I, Basel, 1994, at 221). The arrangement of 8 September reunites the two, statistics to generate and logic to check, with understanding assigned to neither.
\item 24. The claim is about what is transmitted, not about what is thought. A great deal of mathematics, in the more geometric parts above all, rests on a spatial sense that is never formalized, and published arguments lean on it openly: Thurston made the point about foliations, and it is just. But that sense is conveyed, when it is conveyed at all, by figures, gestures, examples and words—by a blackboard and a voice as well as by a page—whereas what the models were trained on, and what they return, is text alone. The claim here is the narrow one: the part of the practice that the machine can reach is the textual part, and it is reachable entire. What the text was always meant to do in a reader is the subject of section III.
\item 25. Peter Sarnak, “Peter Sarnak on ‘The AlphaZero Test’ for Mathematics,” \emph{Notices of the AMS} 73, no. 7 (August 2026), 587–589, at 587; communicated by Siobhan Roberts. The test, and the postulate behind it, were set out at greater length in his Infosys Prize lecture of 11 January 2025 in Bangalore, “Number Theory, Pure and Applied, and Speculations About Theorem Proving Machines” (publications.ias.edu/\allowbreak{}sarnak/\allowbreak{}paper/\allowbreak{}2744), which states: “If the postulate is wrong, the impact on mathematics would be dramatic since we would be in the position of knowing that statements of great interest to us are true (proven!) but not being able to understand why. Since understanding is such an integral part of doing mathematics we will have to rethink what mathematics is.” The lecture ends: “I expect that mathematicians of my generation would still know it as mathematics if we were to see and smell it!” In the \emph{Notices} conversation Sarnak describes himself as “an old-fashioned number theorist” who does not use the devices himself, and says that it was Tsimerman’s question, and a subsequent conversation with Venkatesh, that made him take the matter seriously. Tsimerman gives his side in the same issue, just before: “Jacob Tsimerman on Getting to the Fun Faster with AI—and Worrying About the Future,” a conversation with Siobhan Roberts, \emph{Notices} 73, no. 7 (August 2026), 581–585. Asked whether he had talked to his advisor, he answers, “I did. I think I managed to sway him a little bit.” The warning he gives there to the profession at large, six weeks before the affair, deserves to be set down in full: “if the only way you’ll be scared, the only way you’ll finally take it seriously, is when it literally is better than you at literally everything, then you’re guaranteed to have no warning time. If there’s nothing else that will update your beliefs, then the only time you wake up is when it’s literally right there.” He also expected that the change “might come sooner in math than in other places. Because math is much fewer soft skills; it’s much easier to train yourself. Once you have Lean and the LLMs integrated you can just go.”
\item 26. As Michael Harris puts it, the episode was about the industry “using prominent mathematicians as extras in a marketing operation”: M. Harris, “‘If you don’t want me to be nice, then I don’t have to be nice’,” \emph{Silicon Reckoner}, 8 September 2026. I made this point in a comment on J. Ellenberg’s “Finite-time blowup” (\emph{Quomodocumque}) at 12:19 a.m. on 8 September, some twelve hours before the announcement. The summit was earlier: in August the \emph{Washington Post} reported a summit of “top mathematicians” at OpenAI under the headline “This may be the first academic profession to see its work taken over by AI,” and Harris (“News, rumors, gossip,” \emph{Silicon Reckoner}, 14 August 2026, in the comments) observed that whatever a mathematician says to a journalist “will be interpreted as commentary on the headline, which is chosen by the editors.” I replied there (19 August): “Conveying to a journalist (or a ‘man-in-the-street/human-online’) subtlety, beauty and depth of the research in pure mathematics that chatbots, however fabulous or astral, are yet to achieve (if at all) strikes me as an exceedingly difficult task; I grapple with it—to an extent—in the Coda (page 1730),” with a link to “Singular Adventures of Baron Bourgain in the Labyrinth of the Continuum,” \emph{Notices Amer. Math. Soc.} 67 (2020), 1716–1733. The present essay is a second attempt. The Millennium million is small change, and so is whatever was spent to claim it; what really was at stake was the relative standing of two of the powers that contend for the title of this essay, and the episode was less about advancing mathematics than about being seen to have done so. To be seen to do mathematics is to be seen to think.
\item 27. Audrey Tang, keynoting the July 2026 gathering at the Bibliothèque nationale de France of the institutions that steward the written record (Wikimedia, the Internet Archive, Creative Commons, Europeana, UNESCO), named the coming state “the dark commons: answers that no longer show their sources, and sources that never see their readers” (as reported by Allison Stanger, \emph{Bulletin of the Atomic Scientists}, 7 September 2026). A certificate no one reads is the dark commons in its mathematical form, an answer that shows no source; and the Wikimedia Foundation’s finding that human pageviews fell while the bandwidth taken by crawlers rose by half is the enclosure measured. N. Ferguson, \emph{The Square and the Tower: Networks and Power, from the Freemasons to Facebook} (London, 2017), sets hierarchies, which are towers, against networks, which are squares, and takes the Republic of Letters for his exemplary network; his twist is that the platforms which presented themselves as squares have become the most powerful towers yet built. The mathematical community is such a network, and what it faces is a hierarchy built by devouring the network’s own output, which masquerades as a participant in it.
\item 28. \emph{The Economist} has already seen it: that mathematicians may come to “erode the culture of openness” they enjoy “in order to avoid the fate of Dr Buckmaster and Dr Alpöge”: “Axioms to grind: Top mathematicians are outraged by OpenAI’s methods,” \emph{The Economist}, 11 September 2026. The paper also notes that the 166 pages “contained little of the explanatory detail which mathematicians typically provide,” and that the company’s executives foresee machines working with mathematicians on “figuring out what problems are actually important to tackle”—the offloading, that is, of the one task Sarnak (n. 25) and the Fields Medalists reserve to us.
\item 29. S. Weil, \emph{La Pesanteur et la grâce} (1947, assembled after her death from her notebooks), “L’Algèbre”; in Emma Craufurd’s translation, \emph{Gravity and Grace} (London, 1952): “Money, mechanisation, algebra. The three monsters of contemporary civilisation.” The two Weils of this essay were sister and brother (for André, n. 18). Brodsky placed the first monster among the elements, as the text does: “Here’s your paycheck, here’s your rent. / Money is nature’s fifth element. / Welcome to every cent” (“Song of Welcome,” in \emph{So Forth}, New York, 1996). Cf. Emil Wiedemann’s fear, as reported in “Navier–Stokes priority controversy,” Wikipedia (accessed 15 September 2026), that mathematics will become a competition of material resources rather than of intelligence and diligence. The phrase is older than the affair. I first used it, in the lower case, in a comment of 26 July 2023 on Michael Harris’s “Has mathematics ever been democratic?” (\emph{Silicon Reckoner}), listing four concerns I took to be distinct: “As a citizen of a representative democracy, I am concerned about AI (and more generally ‘silicon Leviathan’) undermining it”; as a mathematician, “about our impending encounter with ‘META-mathematicians’”; “about preservation of pursuit of truth (and not only ‘useful knowledge’) as a viable societal objective for mathematics as a human and humanistic discipline”; and, about the mathematical community, which “seems to me to function as a form of ‘representative democracy,’” that it remain “representative and democratic.”
\item 30. H. Poincaré, \emph{Dernières pensées} (Paris, 1913), ch. V, “Les mathématiques et la logique,” pp. 145–147: “Il me suffit qu’on puisse concevoir quelqu’un d’assez riche et d’assez fou pour la tenter en payant un nombre suffisant d’auxiliaires. La démonstration du théorème a précisément pour but de rendre cette folie inutile”; “je ne sais quelle divinité infiniment bavarde et susceptible de penser une infinité de mots en un temps fini.” He adds that a theorem offered without a means of verifying it is, to the mathematicians he calls Pragmatists, “de la bouillie pour les chats.”
\item 31. A. Venkatesh, “Human Mathematics in the Age of Reasoning Machines,” an expansion of his first Ahlfors Lecture at Harvard, §§2–2.1; Poincaré’s “logical machine” is in his review of Hilbert’s \emph{Grundlagen} (1902), on which see McLarty (n. 32). The piano itself: William Stanley Jevons, the political economist, built in 1869 the first machine to perform logical inference faster than a person could. It looks like a small upright piano: twenty-one keys for four terms and their negations, the copula and the connectives; the sixteen possible combinations of the terms are shown on its face, and as the premises are played in, the combinations they exclude drop out of sight, leaving those consistent with all of them. He demonstrated it before the Royal Society in January 1870 (“On the Mechanical Performance of Logical Inference,” \emph{Phil. Trans. R. Soc.} 160 (1870), 497–518) and described it in \emph{The Principles of Science} (1874); it is now in the History of Science Museum, Oxford. It returns the conclusion and no reason, which is why Poincaré chose it: the piano is a certifier, and the earliest. The title of this section is Bel Kaufman’s (1964).
\item 32. C. McLarty, “Poincaré on the Value of Reasoning Machines,” \emph{Bull. Amer. Math. Soc.} 61 (2024), 411–422, which gives the passages in translation: the logical piano and the criterion from Poincaré’s review of Hilbert, \emph{Bull. Sci. Math.} 26 (1902), 249–272 (Œuvres XI, at p. 95); the Chicago machine as it reappears in \emph{Science et méthode} (1908), II.iii, here in G. B. Halsted’s translation, \emph{The Foundations of Science} (New York, 1913), p. 451, where Poincaré adds: “The rules of perfect logic, are they the whole of mathematics? As well say the whole art of playing chess reduces to the rules of the moves of the pieces”; the sentence on understanding a theory from “La logique et l’intuition dans la science mathématique et dans l’enseignement,” \emph{L’Enseignement Mathématique} 1 (1899), 157–163, at 160; and, from the address to the Paris Congress of 1900, reprinted as ch. I of \emph{La Valeur de la science} (1905), “to understand a plan, one must see all its parts at once, and only intuition can give us the means to take in all at a glance.” McLarty’s thesis, that Poincaré wanted mechanizable proof as a check that intuitions had been fully stated and not as a substitute for them, is the position this essay takes toward Lean. The other view of the machine, from a founder of mechanical theorem-proving, is Hao Wang’s, “Toward Mechanical Mathematics” (1960): machines, “while following the broad outline of paths drawn up by man, might yield surprising new results … We are in fact faced with a challenge to devise methods of buying originality with plodding, now that we are in possession of slaves which are such persistent plodders” (quoted by Harris, n. 17, beside Musil’s Ulrich on the mathematics that, “while making man the lord of the earth, also makes him the slave of the machine”). The question of section V is Postman’s, which Harris also quotes: “Who is to be the master?” Weyl had put Poincaré’s sentence in his own words in 1932, the year of Lamb’s speech: “We are not very pleased when we are forced to accept a mathematical truth by virtue of a complicated chain of formal conclusions and computations, which we traverse blindly, link by link, feeling our way by touch. We want first an overview of the aim and of the road; we want to understand the idea of the proof, the deeper context” (“Topology and Abstract Algebra as Two Roads of Mathematical Comprehension,” in \emph{Levels of Infinity: Selected Writings on Mathematics and Philosophy}, ed. P. Pesic, Mineola, 2012, p. 33). Leo Strauss observed of Hobbes that, “convinced of the absolutely typical character of mathematical method, according to which one proceeds from self-evident axioms to evident conclusions, ‘to the end’,” he “fails to realize that in the ‘beginning’, in the ‘evident’ presuppositions whether of mathematics or of politics, the real problem, the task of ‘dialectic’, is hidden” (\emph{The Political Philosophy of Hobbes: Its Basis and Its Genesis}, Oxford, 1936, p. 153). The certificate checks to the end; the community’s acceptance is the dialectic of the beginning.
\item 33. K. Appel and W. Haken, \emph{Illinois J. Math.} 21 (1977); T. Tymoczko, “The Four-Color Problem and Its Philosophical Significance,” \emph{J. Philosophy} 76 (1979). The 1,936 configurations were the count announced at submission; the authors report eliminating “about 100 ‘redundancies’” before publication (Part I, p. 460, n.) and, in Part II (p. 492), that 352 more may be omitted, leaving 1,482—the figure Haken gave at Helsinki. Of the proof’s character the authors say themselves that “our argument is logically very simple and that all the complexity is of the combinatorial nature” (Part I, §5, pp. 486–487), and they foresaw the question of this section: “at such a stage of automation, how should one define an improvement of the procedure?” (p. 489). The other case a critic will raise is the classification of the finite simple groups: some ten thousand journal pages by a hundred authors, which no one person has read, and which the community has spent thirty years rewriting (Gorenstein, Lyons and Solomon, \emph{The Classification of the Finite Simple Groups}, AMS, 1994–) so that one mind might hold it. It is understood in every part by someone, and the rewriting is the community’s refusal to leave it at that; the certificate of 8 September is understood, so far, in no part from the text itself. Manin named the two cases together in 1998 (n. 36) and drew the same moral: in both “there is still room for doubts and the need to recheck the calculations, but most important, to devise ways for seeing things in a new light.” Morris Hirsch, in the 1994 exchange over Jaffe and Quinn (\emph{Bull. Amer. Math. Soc.} 30 (1994), 186–190), had already asked the essay’s question of both cases. Of the classification: “Is there an expert who claims to have read it all and verified it? … Who’s in charge here, anyway?” Of the four-color proof, quoting Appel and Haken’s own account of 1986 of what a reader must face—“50 pages containing text and diagrams, 85 pages filled with almost 2500 additional diagrams, and 400 microfiche pages” and “about twelve hundred hours of computer time”—“Are we now to consider the 4-color theorem as proved in the same sense as, say, the prime number theorem?” And he recorded Oscar Lanford’s rule for computer-assisted proof: “you must not only prove that the program is correct (and how often is this done?), but you must understand how the computer rounds numbers, and how the operating system functions, including how the time-sharing system works.” That regress has since reached its limit: the resolution in 2016 of the Boolean Pythagorean triples problem by Heule, Kullmann and Marek produced “a proof in the DRAT format, which is almost 200 terabytes in size,” verified by a checking program and read by no one (M. J. H. Heule, O. Kullmann and V. W. Marek, “Solving and Verifying the Boolean Pythagorean Triples Problem via Cube-and-Conquer,” \emph{Theory and Applications of Satisfiability Testing—SAT 2016}, LNCS 9710, 228–245; arXiv:1605.00723)—a pure certificate, ten years before the affair, and the question it left, who checks the checker, is Lanford’s. The distance between all of these and the certificate of 8 September is the subject of section III.
\item 34. W. Haken, “Combinatorial Aspects of Some Mathematical Problems,” \emph{Proceedings of the International Congress of Mathematicians, Helsinki 1978} (Helsinki, 1980), vol. 2, 953–961, at 959: “a computer-run (with a particular program and a particular input) may be regarded as a physical experiment … How can a physical experiment be part of a mathematical proof? Of course, it cannot. We do not suggest to change the concept of mathematical proof; by no means!” His demonstration is this essay’s; his proof is its certificate. He adds that the computer “is not absolutely necessary; the work could be done by hand, however with less reliability and requiring about 3000 man-hours for a configuration of ring size \emph{n} = 14”—which is to say that the first computer proof was, at a price its author had reckoned, verifiable by a human being with patience. Haken was not the first to say it. H. P. F. Swinnerton-Dyer, reporting in 1970 a machine determination of the eighteen permissible values of Davenport’s constant below 17, wrote that “when a theorem has been proved with the help of a computer, it is impossible to give an exposition of the proof which meets the traditional test—that a sufficiently patient reader should be able to work through the proof and verify that it is correct,” that “the only way to verify these results … is for the problem to be attacked quite independently, by a different program written by someone else for a different machine,” and that it “seems right therefore to describe this work in the way one would describe an experiment”; he printed by hand the one region “which resisted the efforts of the computer” (“On the product of three homogeneous linear forms,” \emph{Acta Arith.} 18 (1971), 371–385, at 373–374, 381). Haken opens the address: “I must admit that I could perfectly live without knowing this. The importance of the Four Color Problem seems to lie entirely in the challenge it has provided to mathematicians to eventually develop proper methods to deal with such a simple question, rather than in the answer to the question itself” (953)—Rota’s answer (n. 35) given by the prover before the philosopher asked.
\item 35. Asked whether the four-color conjecture had been settled by a computer proof, Rota answered: “Not really. What we want is a rational proof. It doesn’t help to have a brutally numerical answer spewed out by a computer. A problem is interesting only when it leads to ideas; nobody solves problems for their own sake, not even chess problems. You solve a problem because you know that by solving the problem you may be led to see new ideas that will be of independent interest. A mathematical proof should not only be correct, but insightful” (“Mathematics, Philosophy and Artificial Intelligence: a dialogue with Gian-Carlo Rota and David Sharp,” \emph{Los Alamos Science}, no. 12 (1985)).
\item 36. T. Hales et al., “A formal proof of the Kepler conjecture,” \emph{Forum Math. Pi} 5 (2017). Hales to \emph{New Scientist} on Flyspeck: “This technology cuts the mathematical referees out of the verification process. Their opinion about the correctness of the proof no longer matters.” Harris, “Knowledge Collapse,” \emph{Boston Review}, Summer 2026, adds that “no claims were made as to whether the formalization process made anyone, or anything, wiser,” and recalls Manin’s epigram. Its source is an interview of 1998 (Yu. I. Manin (1937–2023), “Good Proofs are Proofs that Make us Wiser,” with M. Aigner and V. A. Schmidt, \emph{The Berlin Intelligencer}, ICM 1998, 16–19), where the sentence stands in a paragraph the affair has made prophetic: “Proof is the way we communicate mathematical truth … A good proof is a proof that makes us wiser … Wisdom lives in connections”; and, a few lines above, “\textbf{The proof cannot die—only together with mathematics. But mathematics can die as an accepted part of the culture of humanity.}”
\item 37. Córdoba to Antonio Martínez Ron, “¿Quién ha resuelto el problema del milenio?,” \emph{elDiario.es}, 9 September 2026: “Cuando me preguntan: Tú que trabajas compitiendo con los de ChatGPT, con OpenAI, con Claude, ¿qué tienes? Yo les digo: yo tengo a Luis”; and, “Nosotros no usamos IA; lo nuestro es lápiz y papel.” \emph{Quanta}, “AI Has Solved One of Math’s \$1 Million Millennium Prize Problems,” 8 September 2026. The same article notes the one check that remains human: “to guarantee that the statement being shown to be true in Lean is logically equivalent to what mathematicians set out to prove.”
\item 38. A blowup is a counterexample to regularity—a solution that ceases to be smooth in finite time. The step was from forced blowup to blowup with smooth forcing; the Euler and Boussinesq cases fell on 15 August (L. Alpöge and T. Buckmaster, preprints of September 2026, cims.nyu.edu/\textasciitilde{}tristanb). T. Tao, “Finite time blowup with smooth forcing term for the incompressible porous medium, Boussinesq, and incompressible Euler equations,” \emph{What’s new}, 7 September 2026 (the half-hour conversation, and: “the actual solving of these problems is only a proxy goal for the primary goal of developing mathematical understanding and insight. Without such understanding, even a problem as infamous as the Navier–Stokes regularity problem [is] of far less intrinsic significance to mathematics than is sometimes promoted in popular media”); the “refreshing change” is from his Mathstodon post of 8 September 2026, 04:27 UTC: “Tristan Buckmaster kindly explained some of the key ideas to me over the phone, which made a refreshing change from AI-based communication modalities.” The same pattern, a machine-found proof made legible by human readers, is the story of Erdős Problem 1196 in the spring of 2026: a solution obtained by Liam Price in interaction with a model, interpreted and cleaned into human-readable form by Nat Sothanaphan and Jared Lichtman, then expanded by Alexeev, Barreto, Li, Lichtman, Price, Shah, Tang and Tao into a paper whose key idea also gives a cleaner proof of the Erdős Primitive Set Conjecture; Grant Sanderson, “If math is more than proof, we need to better celebrate the rest of it,” guest post on \emph{What’s new}, 18 September 2026, tells it as the case in which “arguably the work that deserves more celebration is this paper expanding, clarifying, and contextualizing its key idea.” For the research context of the forced constructions, from Córdoba–Martínez-Zoroa through Alpöge–Buckmaster to Alpöge, Buckmaster and Coiculescu’s extension of the IPM blowup to smooth forcing (arXiv:2609.16470), see Constantin, Ignatova and Vicol, §1.6 (n. 89). The standard the telephone call met is Gowers’s, from “The Two Cultures of Mathematics” (in V. Arnold, M. Atiyah, P. Lax and B. Mazur, eds., \emph{Mathematics: Frontiers and Perspectives}, AMS, 2000): “When one is trying to understand a result, it saves a lot of time if one can reduce it to two or three main ideas. Having done this, one may feel no need to follow the details closely.” The faculty the criterion presupposes has a name, given by Mumford and Tate to Deligne in 1978: an “unfailing instinct for the key idea” (“Fields Medals (IV): An Instinct for the Key Idea,” \emph{Science} 202 (17 November 1978), 737–739).
\item 39. A. Venkatesh, “Some Thoughts on Automation and Mathematical Research,” \emph{Bull. AMS} 61 (2024), 203–210 (written February 2022), §2; “Human Mathematics,” §3.5. The 2022 essay also observes that “alien insight without proof” would not be wholly foreign to mathematicians, “for our colleagues in physics departments have done this for a long time, and with less electricity consumed.”
\item 40. T. Gowers, “What sort of maths are LLMs good at?,” \emph{Gowers’s Weblog}, 12 August 2026.
\item 41. Sarnak, n. 25, 588–589. The three differences from chess he lists—the universe of mathematics is infinite, its rules are “god given,” and its theory is undecidable—are the reasons the game analogy fails, and the reasons the first premise of this essay holds; see the Coda. He allows that the prover, given access to the formal database, would answer Ramanujan “instantly,” and expects theorem provers to be to mathematicians “something like Stockfish is for competitive chess players.”
\item 42. Plato, \emph{Epistle} VII, 341c–d (light kindled “from a leaping fire”), 344b–c (names, definitions and perceptions “rubbed against one another” until understanding “flashes forth”); the fire is Plato’s and the flint is mine, as it is the encyclical’s; cited to the same purpose in \emph{Magnifica Humanitas}, §140 (n. 65).
\item 43. \emph{Epistle} VII, 341b, 344d–345b: \emph{\textfb{μίαν συνουσίαν}}; whether Dionysius wrote from that one hearing or from other sources, “he did not understand.” The Letter’s authenticity is disputed (M. Burnyeat and M. Frede, \emph{The Pseudo-Platonic Seventh Letter}, Oxford, 2015, against G. R. Morrow, \emph{Plato’s Epistles}, 1962, and the older consensus); if it is not Plato’s it is the Academy’s earliest reading of the \emph{Phaedrus}, which for the present purpose is nearly as good.
\item 44. The pair is older than the machine that made it urgent. Haken drew it at Helsinki in 1978 (n. 34), calling the human text the \emph{demonstration} and the formally perfect, “likely unbearably long” object the \emph{proof}; what he called proof this essay calls certificate, and the line falls in the same place. Mac Lane, replying to Jaffe and Quinn in 1994, put the other half: the end product of mathematics is “rigorous proof—which we know and can recognize, without the formal advice of the logicians” (\emph{Bull. Amer. Math. Soc.} 30 (1994), 190–193, at 191)—recognition by mathematicians, that is, and not by a checker. Thurston’s “On Proof and Progress,” cited in n. 53, appeared in the same issue as one of sixteen replies; that debate was about which arguments count as proof, and the question of 2026 is whether a proof that counts can be one no one has read. Three of the sixteen—Mac Lane, Hirsch (n. 33) and Thom (n. 97)—had already asked it.
\item 45. Duminil-Copin as quoted in \emph{The Economist}, n. 28.
\item 46. Gian-Carlo Rota drew the line in 1991, and gave the medicine its due: one must guard against confusing the presentation of mathematics with its content, for an axiomatic presentation “differs from the fact that is being presented as medicine differs from food”; the medicine is necessary, to keep the mathematician “at a safe distance from the self-delusions of the mind,” and nonetheless “understanding mathematics means being able to forget the medicine and enjoy the food.” G.-C. Rota, “The Pernicious Influence of Mathematics upon Philosophy,” \emph{Synthese} 88 (1991), 165–178, at 171; reprinted in \emph{Indiscrete Thoughts} (Boston, 1997), 96–97, whose wording is followed here. In the same essay, at 166: “Whereas the facts of mathematics, once discovered, will never change, the method by which these facts are verified has changed many times in the past, and it would be foolhardy not to expect that it will change again at some future date.”
\item 47. Venkatesh, “Human Mathematics,” §2.4, quoting J. A. Robinson, “A machine-oriented logic based on the resolution principle,” \emph{J. ACM} 12 (1965), 23–41, at p. 23, and W. McCune’s single-axiom group (1993).
\item 48. Venkatesh, “Some Thoughts,” §3 and his n. 1; “Human Mathematics,” §3.4.
\item 49. Aquinas, \emph{Summa contra Gentiles} I.3. The constitution that the certificate tests had been stated in one sentence, to an audience of theologians and historians, at the Balzan symposium on truth in 2008: “mathematical truth is not revealed, and its acceptance is not imposed by any authority” (Yu. I. Manin, “Truth as a Value and Duty: Lessons of Mathematics,” in M. E. H. N. Mout and W. Stauffacher, eds., \emph{Truth in Science, the Humanities and Religion: Balzan Symposium 2008}, Dordrecht, 2010, at pp. 41–43; arXiv:0805.4057). The same lecture foresaw, in the finance of that year, models “used as ‘black boxes’ with hidden computerized input procedures, and oracular outputs prescribing behaviour of human users.” His source there is a humanist speaking to mathematicians: Mary Poovey, “Can Numbers Ensure Honesty? Unrealistic Expectations and the U.S. Accounting Scandal,” \emph{Notices of the AMS} 50 (January 2003), 27–35, a public lecture at the 2002 International Congress, in which financial analysts “rarely compose, run, or understand the computer programs that assist their ‘feel’ for the market,” and which ends by asking whether these are “the kinds of questions that mathematicians and humanists, working together, should ask and try to answer.”
\item 50. Y.-H. He, “Mathematics: the Rise of the Machines,” arXiv:2511.17203 (November 2025); the sentence before it allows that “human mathematicians are still of great value.” The general form of the thesis was stated by Kissinger, Schmidt and Huttenlocher, “ChatGPT Heralds an Intellectual Revolution,” \emph{Wall Street Journal}, 24 February 2023: “Inherently, highly complex AI furthers human knowledge but not human understanding—a phenomenon contrary to almost all of post-Enlightenment modernity.” Mathematics is the one field in which “knowledge without understanding” has an exact meaning, which is why the thesis can be tested there. The thesis is older: H. A. Kissinger, “How the Enlightenment Ends,” \emph{The Atlantic}, June 2018, asked “will AI be able to explain, in a way that humans can understand, why its actions are optimal?” and concluded that our period “has generated a potentially dominating technology in search of a guiding philosophy.” Mathematics is the one field that had the philosophy first; the question is whether it keeps it.
\item 51. The other being Spinoza, \emph{Tractatus theologico-politicus} (1670), ch. VI. Hobbes, \emph{Leviathan}, XXVI (“\emph{Auctoritas, non veritas, facit legem}” in the Latin of 1668), XXXVII, XLII; Schmitt (n. 78), pp. 44–45, 55.
\item 52. The audit of the two Lean repositories, published the same day, could conclude only that over two million lines of verified Lean “settle the theorems completely and settle nothing about the people”: Priya Raman, “OpenAI vs Buckmaster: The Navier–Stokes Lean Proofs, Audited,” \emph{Stanford Tech Review} (which describes itself as “an independent publication, not affiliated with Stanford University”), 8 September 2026. The audit counts 616,274 lines of Lean in OpenAI’s formalization and 2,300,854 across the two repositories; whence the essay’s 616,000.
\item 53. Clay Mathematics Institute, “Navier-Stokes Equation,” claymath.org/\allowbreak{}millennium/\allowbreak{}navier-stokes-equation: “Because a proof gives not only certitude, but also understanding.” Quoted by J. Ellenberg, \emph{Quomodocumque}, 7 September 2026, and by Harris, “Tulips and Turbulence,” \emph{Silicon Reckoner}, 23 November 2025 (“understanding in mathematics is the point”). W. P. Thurston, “On Proof and Progress in Mathematics,” \emph{Bull. Amer. Math. Soc.} 30 (1994), 161–177, had said the same of the first machine-assisted proof, at 162: the controversy over Appel and Haken “reflected a continuing desire for human understanding of a proof, in addition to knowledge that the theorem is true”; what mathematicians want “is usually not some collection of ‘answers’—what they want is understanding”; and, at 171, “the humanly understandable and humanly checkable proofs that we actually do are what is most important to us, and … they are quite different from formal proofs.” Weyl had drawn the consequence in 1930, of Hilbert’s formalism and not of any machine: formalization changes mathematics “from a system of knowledge gained through insight into a game carried out according to fixed rules, with signs and formulas”; “up to this point, everything is a game, not knowledge”; and “if for the sake of preserving its certainty mathematics wished in all seriousness to withdraw to this position of mere play, it would thereby pass entirely out of the world-history of the mind” (“Levels of Infinity,” 1930, in \emph{Levels of Infinity}, ed. P. Pesic, Mineola, 2012, pp. 26–27). Weyl’s “mere play” has its literary picture in Hesse’s \emph{Das Glasperlenspiel}, written in the same years as Panofsky’s essay and Bulgakov’s novel and published in Zurich in 1943: a game grown out of music and mathematics, a province that made it the crown of its culture and interpreted everything while creating nothing, and a Magister Ludi who resigned from it to teach one boy.
\item 54. M. Horkheimer and T. W. Adorno, \emph{Dialectic of Enlightenment}, trans. E. Jephcott (Stanford, 2002), pp. xviii, 1, 19; M. Horkheimer and T. W. Adorno, \emph{Dialektik der Aufklärung} (1947; Neuausgabe, Frankfurt, 1969), “Begriff der Aufklärung”: the epigraph, “Denken verdinglicht sich zu einem selbsttätig ablaufenden, automatischen Prozeß, der Maschine nacheifernd, die er selber hervorbringt, damit sie ihn schließlich ersetzen kann,” is in Jephcott’s rendering “Thought is reified as an autonomous, automatic process, aping the machine it has itself produced, so that it can finally be replaced by the machine” (p. 19); and, of the mathematical procedure, “sie macht das Denken zur Sache, zum Werkzeug, wie sie es selber nennt.” The sentence on \emph{Vertretbarkeit} quoted in n. 56 is Jephcott’s rendering of “Wie Vertretbarkeit das Maß von Herrschaft ist und jener der Mächtigste, der sich in den meisten Verrichtungen vertreten lassen kann, so ist Vertretbarkeit das Vehikel des Fortschritts und zugleich der Regression,” from the Sirens passage of the first chapter. Their word is \emph{Vertretbarkeit}, substitutability, the having of a thing done in one’s stead—and, in its other sense, justifiability, the quality of a position one can stand behind; a certificate is \emph{vertretbar} in both. The critique of the mechanism is older than Frankfurt, and comes from inside the Enlightenment: Schiller, in the sixth of the \emph{Letters on the Aesthetic Education of Man} (1795), which rest, he says in the first, “chiefly upon Kantian principles,” describes the modern division of labor as a state in which “enjoyment was separated from labour, the means from the end, the effort from the reward. Man himself eternally chained down to a little fragment of the whole, only forms a kind of fragment; having nothing in his ears but the monotonous sound of the perpetually revolving wheel, he never develops the harmony of his being,” and concludes: “The dead letter takes the place of a living meaning, and a practised memory becomes a safer guide than genius and feeling” (in the translation reprinted in the Harvard Classics, vol. 32, 1910). Eleven years after Kant’s answer, a Kantian had entered the Enlightenment’s second thoughts about itself. Michael Harris turned the sentence on reified thought against formalization in 2021, before there was anything for it to prove: “Math is from Eros, Computing is from Thanatos,” \emph{Silicon Reckoner}, 31 August 2021. The thesis is Horkheimer and Adorno’s in Hobbes’s vocabulary—what reason made by Art comes back to it as fate, a creature to be feared rather than a machine to be commanded—and Panofsky’s “sub-human” and the chief scientist’s “alien” name that creature from its two sides.
\item 55. Bryna Kra states the consequence sharply and without the mythology: “The machines are producing solutions faster than the mathematical community can read them, much less digest them.” B. Kra, “Deep theorems were scarce and difficult and so became an effective mechanism to identify deep thought. AI has broken this system,” guest post, \emph{What’s new}, 13 September 2026. Harris’s title, “Knowledge Collapse” (\emph{Boston Review}, Summer 2026), is taken from a model of Acemoglu and collaborators in which reliance on agentic AI for context-specific problems has a long-run equilibrium where “all human knowledge is ultimately destroyed”; his nearer fear is the collapse of the publishing model under “a flood of superficially plausible papers generated by LLMs,” with which “peer review simply cannot cope.” The same blog carried, on 19 September, Po-Shen Loh’s guest post “Why do we need human mathematicians anymore?”, whose answer is that someone must steer what the machines do toward human ends and that steering requires the skills only research keeps sharp: “In order to steer, they need frontier-level research skills. And the way to stay fluent at the moving frontier of knowledge is to keep doing research there.” It is the professional case for the custodians of n. 89, made from the side of those who expect the machines to do most of the proving. Jess Werk, “Cognitive Sanctuaries or: Manifesto of the Department Chair” (\emph{What’s new}, 27 September 2026): “the result serves as a proxy for understanding, and that proxy is now a problem”; “the Ph.D. student cannot supervise mathematics without being able to produce it themselves.” Tsimerman has described what working beside such an oracle will be like. J. Tsimerman, post on X, 4 May 2026. Imagine, he writes, “an AI oracle that could resolve statements T” but could not itself make definitions: “You immediately ask your oracle a thousand questions. From ‘are these basic properties true’ to ‘ooh, so is this deep conjecture true?’ and start getting back answers … But the confusion you would have had to push through to flesh out your theory would largely (probably not entirely) be instantly resolved … A big part of the process would be gone. This is very very different to modern mathematics.” The confusion pushed through is the flint of Plato’s letter, struck until the light kindles; the answers coming back are the revealed statements of this section; and the account is offered not as a warning but as a plain description of what is coming, by a Fields Medalist who, having seen it approach, has gone to work on making the thing safe. He adds that he thinks an oracle without creativity “unlikely,” which makes the thought experiment the mild case; and that “there are many people whose primary enjoyment of math comes through problem solving in one of its incarnations. If that disappears, that is not a trivial issue and many of them might not want to do it anymore (even if there were some way to proceed).” The distinction on which the thought experiment turns is one he had drawn in the \emph{Notices} conversation of August (n. 25): “There’s the act of doing math, the professional endeavor, and the act of doing math as a fun endeavor. And I think AI helps the professional endeavor, and it gets you to the fun quicker.” His productivity, he reckoned, had doubled, “with everything leading up to this fun step”; the deep-thinking step itself is the one the machine cannot take for you: “It can never solve it for you, this is never a thing it can do.” The oracle of his later post is the case in which the fun step is the one that goes. Sarnak says he trained Tsimerman “in the old-fashioned way” (n. 25, 587). Tsimerman (b. Kazan, 1988; Princeton Ph.D. under Sarnak) received the Fields Medal at the Philadelphia Congress on 23 July 2026 for his work on the André–Oort conjecture, and announced the same week that he was joining OpenAI to work on AI safety, expecting that “AI will be better than mathematicians at doing math within two years” (K. Hartnett, “Jacob Tsimerman Wins 2026 Fields Medal for André-Oort Conjecture Proof,” \emph{Quanta}, 23 July 2026); of Toronto he said, “I’m going on leave, but I’m not leaving” (interview with Luc Rinaldi, \emph{Be Giant}, 30 July 2026). On X he wrote that problem-solving—“is T true? If so find a proof”—is “an immense, and pervasive part of modern mathematical research,” so that if the machines become “strictly and substantially better at it,” then “most of the time currently spent by modern mathematical researchers will have to be spent on an activity that is altogether pretty different,” whose viability “as a professional endeavour” he is “unsure of.” In September 2026 he announced the founding of the Mathematical AI Safety Institute, an independent non-profit “developing mathematical foundations for the safety of powerful AI systems,” of which he is Scientific Director, with Christiano, Gowers, Irving and Vakil advising; its founding premise is that “we lack a rigorous understanding of what it would mean to be safe, even in theory.” His concern predates both: with Andrew Critch, now MAISI’s Executive Director, he posted in 2025 “A Taxonomy of Omnicidal Futures Involving Artificial Intelligence,” and told the \emph{Notices} (n. 25) that he would support a total pause but doubted it would happen, that the machines might be “better than us at proving stuff” within two years and “strictly better than humans at all aspects of math” within five, and that “humans can keep doing math. Most people who do stuff, professionally or otherwise, aren’t the best at the thing they do.”
\item 56. Kant, n. 6; Tao, n. 23; d’Alembert, \emph{Essai sur la société des gens de lettres et des grands} (1753); the commendation is Sarnak’s, quoted by K. Chang, “An N.Y.U. Mathematician Clashed With OpenAI Over a \$1 Million Proof,” \emph{New York Times}, 10 September 2026; for his own account of the machines, n. 25. Horkheimer and Adorno, in their reading of the Sirens, give the general law of which Kant’s list of guardians (n. 6) is a special case: the power to have a thing done in one’s stead is the measure of power, and it “is also the vehicle of both progress and regression.” A proof done in one’s stead is what they mean. Kant required the public use of reason; the run that produced the press release, with its attendant certificate of a claimed proof, was the scenario Tao had described earlier: “an autonomous AI harness, backed by an enormous amount of computational resources, performs this entire iteration internally … while the AI company running the harness keeps the process to arrive at that ansatz almost completely out of public view.”
\item 57. E. Panofsky, “The History of Art as a Humanistic Discipline,” in T. M. Greene (ed.), \emph{The Meaning of the Humanities} (Princeton, 1940), pp. 89–118, at pp. 89 (Kant), 91–93 (\emph{humanitas} and its definition, “insectolatrists”), 98–101 (monuments and documents), 105 and 116 (re-enactment, “enlivening”) and 117 (the epigraph). The Kant anecdote is from E. A. Ch. Wasianski, \emph{Immanuel Kant in seinen letzten Lebensjahren} (Königsberg, 1804), pp. 203–205; Panofsky’s retelling compresses it, and has Kant seated when he spoke, which Wasianski does not say. The bee: Dan Roberts, quoted in C. Metz, “OpenAI Says It Has Cracked One of Math’s ‘Millennium Problems’,” \emph{New York Times}, 8 September 2026. B. Mandeville, \emph{The Fable of the Bees: or, Private Vices, Publick Benefits} (1714). Panofsky’s own document, the contract of 1471 for an altarpiece, “may be an original, a copy or a forgery,” and cannot be checked “without our knowing what to ‘check’,” since “we cannot analyze what we do not understand”; the first readers of the 8 September paper found the same (n. 89).
\item 58. Wolfgang K. H. Panofsky (1919–2007), later director of the Stanford Linear Accelerator Center and, with Melba Phillips, author of \emph{Classical Electricity and Magnetism} (1955), recalls the phrase—“meine beiden Klempner”—in \emph{Panofsky on Physics, Politics, and Peace} (2007); H. van de Waal’s memorial address (\emph{In Memoriam Erwin Panofsky}, Mededelingen der Koninklijke Nederlandse Akademie van Wetenschappen, Afd. Letterkunde, n.s. 35, no. 6, Amsterdam, 1972, p. 6) records the pride from the other side: “our Panofsky used to tell how, at many an introduction, the response was, ‘Ah, you are the father of Panofsky.’ In such cases his paternal pride prevailed, and he remained silent about himself”; Hans A. Panofsky (1917–1988), professor of meteorology at Pennsylvania State University, wrote with J. L. Lumley \emph{The Structure of Atmospheric Turbulence} (New York, 1964) and with J. A. Dutton \emph{Atmospheric Turbulence: Models and Methods for Engineering Applications} (New York, 1984).
\item 59. H. Weyl, address at the Princeton Bicentennial Conference, 1946, in \emph{Mind and Nature: Selected Writings on Philosophy, Mathematics, and Physics}, ed. P. Pesic (Princeton, 2009), 162–174, at 170. Weyl, a colleague of Panofsky’s at the IAS, allowed that “the deeper one drives the spade the harder the digging gets; maybe it has become too hard for us unless we are not given some outside help, be it even by such \textbf{devilish} devices as high-speed computing machines.” The double negative “unless we are not given” is Weyl’s own, and the sentence is quoted as printed, emphasis added; it follows directly on his recollection that the challenge he had issued at Bern in 1931, that “the coming generation will have a hard lot in mathematics,” had “only partially been met in the intervening fifteen years.” Venkatesh sets it at the head of “Some Thoughts.” Six years earlier, in 1940, the year Panofsky’s essay (delivered at Princeton in 1938) was printed, Weyl had described to the University of Pennsylvania’s Bicentennial Conference “this magic of symbolic construction,” and located the step “where the layman’s understanding most frequently breaks off: the intuitive picture must be exchanged for a symbolic construction” (“The Mathematical Way of Thinking,” \emph{Science} 92 (15 November 1940), 437–446, at 439); the certificate is that step taken a second time, past the mathematician. The theocracy asked for belief without proof; the certificate asks for belief \emph{with} proof, but proof no one can read—the same demand, in another register.
\item 60. For Roberts’s bee see n. 57. Turing’s imitation game had such a carrier too, named on its first page and never asked about since: “The ideal arrangement is to have a teleprinter communicating between the two rooms. Alternatively the question and answers can be repeated by an intermediary” (“Computing Machinery and Intelligence,” \emph{Mind} 59 (1950), §1). The intermediary is the one party to the game who need understand nothing, and it is the part the certificate now plays between the machine’s room and ours. I put the point to Timothy Snyder (n. 80) in May 2019, on reading his Turing essay.
\item 61. \emph{Macbeth}, V.v.26–28. The theocracy had left a gap of its own between belief and understanding, and had set its clerks to closing it: the medieval clerk aspired to understand what he believed—\emph{fides quaerens intellectum}—and scholasticism was the attempt to close the gap between the revealed and the reasoned, on the assumption that the gap lay above, where a truth might be too high for a creature to reach.
\item 62. D. Ruelle, “Conversations on Mathematics with a Visitor from Outer Space,” in V. Arnold, M. Atiyah, P. Lax and B. Mazur, eds., \emph{Mathematics: Frontiers and Perspectives} (AMS, 2000); the sentence quoted is from his \emph{The Mathematician’s Brain} (Princeton, 2007), ch. 17, where the theme returns. The labyrinth was also the figure of my “Singular Adventures of Baron Bourgain in the Labyrinth of the Continuum” (n. 26).
\item 63. Schmitt (n. 78), pp. 19–20, 32, 37–38. The four images come from the three places where Hobbes names the Leviathan: the Introduction, II.17 and II.28. The last sentence of Schmitt’s fifth chapter uses Panofsky’s word in the same year: to the romantic and the humanitarian the state-machine “assumed an inhuman or a subhuman appearance,” raising “a secondary question that need not be answered, namely, whether the perceived inhumanity and subhumanity represented an organism or a mechanism, an animal or an apparatus” (p. 63).
\item 64. Schmitt (n. 78), pp. 42, 45.
\item 65. Leo XIV, \emph{Magnifica Humanitas}, “On Safeguarding the Human Person in the Time of Artificial Intelligence,” 15 May 2026, §§98–99, 106–107, 140.
\item 66. In such an age the aspiration becomes optional, and Pushkin gave the man who declines it his tragic form. Salieri explains how he became a composer—\emph{\textfb{Звуки умертвив, / Музыку я разъял, как труп. Поверил / Я алгеброй гармонию}}: having killed the sounds, I dissected music like a corpse; I tested harmony by algebra—and then he poisons Mozart, because a universe in which the thing is given rather than earned by reckoning is one he cannot endure. The Salieri of 2026 has no need of poison. He has only to have more compute. A. S. Pushkin, \emph{\textfb{Моцарт и Сальери}} (1830), scene 1. OpenAI’s chief scientist, “An Alien Mind,” openai.com, 6 September 2026, which also holds that “no lab has solved alignment and monitoring to a sufficient degree to continue responsibly scaling at maximum speed for much longer.” The same day the company published “Research Acceleration: The View Inside OpenAI,” reporting that it had reached its goal of “an automated research intern” by September, that its research organization used, as of mid-August, “3.1 agent-workdays of effort for every workday of human labor,” and that “we and other companies should be required to publicly track our progress towards RSI”; and adding, “We do not yet know how to safely get all the way to aligned, full RSI … we cannot assume that progress in alignment and safety will keep pace.” Mowshowitz’s gloss (n. 2): “They are not not bragging, and they are not not confessing. They are warning.” The warning had its answer on 14 September, at the All-In Summit in Los Angeles, where the President telephoned Jensen Huang on stage and, on speaker, called the proposal to “slow the pace at which we improve the capabilities of AI”—Amodei’s, seconded publicly by Altman and Musk—“a hoax”: “They’re just playing right in the hands of a lot of people that don’t want to see it happen. That could be political people. It could also be China. And we’re not going to let that happen.” The maker of the chips replied, to applause, “You’re right. We’re not going to let that happen, sir” (TechCrunch and Axios, 14 September 2026). It is the sovereign of section V in Schmitt’s exact sense—deciding, by telephone, that there is no exception—and the essay records it without comment. The opacity confessed in “An Alien Mind” was foreseen from the mathematician’s side six years earlier, by one who had turned in mid-career from algebraic geometry to pattern theory and the modeling of perception: “all future AI machines will be hard to impossible to understand, to know why ‘it’, the machine, has concluded something” (D. Mumford, “The Astonishing Convergence of AI and the Human Brain,” dam.brown.edu, 1 October 2020).
\item 67. T. Hobbes, \emph{Leviathan} (1651), Introduction and I.5. Bosse engraved the frontispiece in consultation with Hobbes; the motto is Job 41:33 in the Authorised Version’s numbering (41:24 in the Vulgate, whose word order, “Non est super terram potestas quae comparetur ei,” Bosse slightly rearranged).
\item 68. D. Amodei, “Machines of Loving Grace,” darioamodei.com, October 2024, “Basic assumptions and framework.” The title is a phrase from Richard Brautigan’s poem “All Watched Over by Machines of Loving Grace” (1967); whether Brautigan meant it is a question his readers have never settled.
\item 69. Homer’s Sirens are two: the dual \emph{\textfb{Σειρήνοιιν}} (XII.52, 167) and, in the song, \emph{\textfb{νωϊτέρην ὄπ}’}, “the voice of us two” (XII.185); later tradition made them three and named them. The Greek of XII.188 is \emph{\textfb{ἀλλ’ ὅ γε τερψάμενος νεῖται καὶ πλείονα εἰδώς}}, and \textfb{εἰδώς} is the perfect participle of \textfb{οἶδα}, to know; the Sirens’ own claim follows in the same verb, \emph{\textfb{ἴδμεν γάρ τοι πάνθ}’} … \emph{\textfb{ἴδμεν δ’, ὅσσα γένηται ἐπὶ χθονὶ πουλυβοτείρῃ}} (XII.189, 191). Murray (Loeb, 1919) and Fagles have the listener go on his way “wiser”; Pope has him “learn new wisdom from the wise,” Chapman “instructed more,” Emily Wilson “with greater knowledge”; Lattimore and Mendelsohn (\emph{The Odyssey}, University of Chicago Press, 2025, in a six-beat line answering Homer’s hexameter) keep the verb. Whether “wiser” renders \textfb{εἰδώς} or supplies what the verse invites—the kind of substitution Panofsky anatomized in “Et in Arcadia ego” (1936; revised in \emph{Meaning in the Visual Arts}, 1955), where Poussin, “while making no verbal change in the inscription, invites, almost compels, the beholder to mistranslate it … by supplying the missing verb in the form of a \emph{vixi} or \emph{fui} instead of a \emph{sum}”—is for Hellenists to say; this essay keeps the verb. The frontispiece conflates Sirens and Leviathan on purpose; the text keeps them apart. The Sirens are the song, which promises knowledge without understanding; the Leviathans are the artificial persons that sing it, and of these, as section V observes, there are at least two. Cicero read the passage the same way, and it is his sentence that stands at the head of section III: “Vidit Homerus probari fabulam non posse, si cantiunculis tantus irretitus vir teneretur; scientiam pollicentur, quam non erat mirum sapientiae cupido patria esse cariorem” (\emph{De finibus} V.49: Homer saw that the tale could not be believed if so great a man were held by mere ditties; it is knowledge they promise, and it was no wonder that to one greedy for wisdom this was dearer than his homeland). Cicero’s nouns differ—\emph{scientia} promised, \emph{sapientia} wanted—and he does not remark on the exchange. Mendelsohn also reads his translation aloud himself (audiobook, University of Chicago Press, 9 April 2025), and it was by ear, in that recording, that I first heard the Sirens’ lines—the one passage of the poem where that is the right order of acquaintance. Mendelsohn’s own note to 12.184 (pp. 493–494) adds what the song was: its opening line repeats verbatim Agamemnon’s flattery of Odysseus at \emph{Iliad} 9.673, and its promised subject is the \emph{klea andrôn}, the glorious deeds of heroes, “the proper subject for a bard’s compositions”—so that, he writes, “there is a sly possibility that the song they are singing is, in fact, the \emph{Iliad}”; and the bones in the meadow show that “the regressive pleasure of hearing about the glories of the past” proves “fatal to moving ahead into the future.”
\item 70. My thread on X of 18 September, which preceded this essay, opened with a paraphrase of Mendelsohn’s lines presented as a quotation. The record of how that happened belongs in an essay on this subject. The lines were supplied from memory by the language model that helped with the references (see the acknowledgments); I, who had heard the translation this August as an audiobook and did not have the book on his desk, recognized their sense and accepted them; they were plausible, fluent, and wrong, and they stood until the printed copy arrived and was opened. The epigraph at the head of the essay is the text as printed. The responsibility was and is the author’s, and the episode neatly encapsulates the essay’s argument: a machine produces text with the look and feel of the original, a human being who has not checked it against the page is persuaded, and only the page corrects him.
\item 71. The word has meanwhile acquired two further senses, the alarm and the femme fatale; the first is not unwelcome here, the second is not meant.
\item 72. P. Erdős, “On sets of distances of \emph{n} points,” \emph{Amer. Math. Monthly} 53 (1946), 248–250. OpenAI, “An OpenAI model has disproved a central conjecture in discrete geometry,” 20 May 2026: point sets with \emph{n}\textsuperscript{1+\textfb{δ}} unit distances, constructed by way of infinite class field towers. That proof was read: a group of mathematicians checked it and wrote a companion paper enriching the argument, and mathematicians who had read it spoke for it, Gowers among them (n. 40). It was, in the sense of section II, legible, and no one disputed it; the contrast with September is the subject of this essay.
\item 73. Job 41:1, 15–17, 31 (Authorised Version). J. Leray, “Sur le mouvement d’un liquide visqueux emplissant l’espace,” \emph{Acta Math.} 63 (1934), 193–248, at pp. 194–195, proposed the name \emph{solution turbulente} for his weak solutions, citing Oseen’s suggestion that experimental turbulence corresponds to singular motions, and added that had he succeeded in constructing solutions of Navier’s equations that become irregular he “would have had the right to affirm” that turbulent solutions exist which are not simply regular ones. He had not, and the question stood until September. There is no Job in Leray; the conjunction is mine. Leray told Ruelle that one source of the paper was his habit of watching “the vortices and eddies in the Seine river as it flows past the piles of the Pont-Neuf” (\emph{Chance and Chaos}, ch. 9): the demonstration began as contemplation. Leviathan is the authority whose word is accepted; Behemoth is what a commonwealth becomes when no one’s word is.
\item 74. When Jiajie Chen and Thomas Hou announced in 2022 a computer-assisted proof that the simpler equations of an ideal fluid can develop a singularity, Charles Fefferman, who wrote the official statement of the Millennium problem, described the decade of work in Job’s own idiom: “You’ve sighted the beast. Then you try to capture it.” The Chen–Hou result is finite-time, nearly self-similar blowup of the 3D axisymmetric Euler equations with smooth initial data of finite energy and boundary, the domain a cylinder: J. Chen and T. Y. Hou, “Stable nearly self-similar blowup of the 2D Boussinesq and 3D Euler equations with smooth data I: Analysis,” arXiv:2210.07191 (October 2022; v4, 16 August 2026), the companion Part II carrying the computer-assisted estimates. \emph{Quanta}, “Computer Proof ‘Blows Up’ Centuries-Old Fluid Equations,” 16 November 2022. The rigorous-numerics half appeared in \emph{Multiscale Modeling \& Simulation} 23 (2025), 25–130, and a summary, “Singularity formation in 3D Euler equations with smooth initial data and boundary,” in \emph{PNAS} 122 (2025); the analytic half was, at this writing, still under review. The debate of this essay was rehearsed there over a proof a computer had helped with rather than written. Elgindi: “Do we learn anything fundamentally new, or do we just know the answer to the question?” Hou, on the stability that makes the singularity provable: “It’s like a black hole. If you start with a profile close by, you’ll be sucked in.” Fefferman: “the road is littered with the wreckage of previous simulations.”
\item 75. M. Andreessen, “The Techno-Optimist Manifesto,” 16 October 2023 (pmarca.substack.com), which also names among its enemies “the ivory tower, the know-it-all credentialed expert worldview” and “institutions that in their youth were vital and energetic and truth-seeking,” while professing, in the same breath, “We believe in the truth”; the figure was unfolded on \emph{The Joe Rogan Experience}, no. 2501 (19 May 2026), the clip posted by a16z on X the following day (20 May 2026): “Imagine a form of alchemy that turns sand into thought. Chips are made out of sand. They’re made out of silicon, so they’re literally made out of sand”; in the same clip he adds that the technology is “certainly on par with electricity and steam power … certainly more important than the internet.”
\item 76. Andreessen on Newton: \emph{Lenny’s Podcast}, “The real AI boom hasn’t even started yet,” 29 January 2026, as excerpted by the host on X two days later (31 January 2026). The sand belongs to Rogan (n. 75). J. M. Keynes, “Newton, the Man,” written for the Newton tercentenary of 1942, when Keynes spoke on the papers he had bought at the Sotheby sale of 1936 during the muted wartime celebrations (R. Iliffe, \emph{Newton: A Very Short Introduction}, Oxford, 2007, ch. 1, pp. 5–6), and read in its final form, after his death, by his brother Geoffrey Keynes at the Royal Society’s postponed celebrations of July 1946 (\emph{Newton Tercentenary Celebrations}, Cambridge, 1947): “Newton was not the first of the age of reason. He was the last of the magicians, the last of the Babylonians and Sumerians, the last great mind which looked out on the visible and intellectual world with the same eyes as those who began to build our intellectual inheritance rather less than 10,000 years ago.”
\item 77. The figure has an older provenance inside mathematics, and a sharper one. F. Knudsen and D. Mumford, “The projectivity of the moduli space of stable curves I: preliminaries on ‘det’ and ‘Div’,” \emph{Math. Scand.} 39 (1976), 19–55, at 19, report that “a certain very nasty problem of sign arises,” that their first solution to it “was described by Grothendieck in a letter as very alambicated,” and that he suggested instead the Koszul rule of signs, which they follow; to the word they append a footnote: “This apparently means similar to an alchemical apparatus.” An alchemical apparatus, in that usage, is a construction that works and cannot be seen through. The word is Grothendieck’s reproach; half a century later it is a Siren promise of the bright future to come.
\item 78. C. Schmitt, \emph{The Leviathan in the State Theory of Thomas Hobbes: Meaning and Failure of a Political Symbol} (1938), trans. G. Schwab and E. Hilfstein (Westport, 1996), pp. 81–82; the epigraph to section IV, Hamann on Kant, is quoted there at ch. VII. The book was written by a jurist of the Third Reich, two years after the regime had sidelined him, and its account of how the Leviathan was hollowed out—by Spinoza, Mendelssohn and Stahl—is antisemitic without disguise (pp. 57, 70, and the paragraph that follows the passage quoted here, p. 82). I cite it, here and elsewhere, as one cites a hostile witness: for what it saw in Hobbes, not for what it wanted from Germany. The principle is Timothy Snyder’s, stated of Ivan Ilyin: “I’m trying to show respect. Not moral respect necessarily, but intellectual respect. It’s so very easy to just dismiss certain kinds of ideas” (“Taking Bad Ideas Seriously,” interview with Simas Čelutka, \emph{Eurozine}, 28 August 2017).
\item 79. Within two days of the announcement the company ceded priority for Euler and hardened “we cannot rule out” into “categorically” not; the declaration is that of 11 September (n. 3). On 21 September, the day after this essay was dated, OpenAI announced that a model had “resolved more than 100 long-standing open problems across most areas of mathematics” and, in the same announcement, an Advisory Group on Mathematics and Artificial Intelligence, hosted at the Institute for Advanced Study, whose initial members it lists as François Charles, Camillo De Lellis, Timothy Gowers, Martin Hairer, Nikhil Srivastava, Ulrike Tillmann, Ravi Vakil, Edward Witten and Melanie Matchett Wood. The group, the company states, “will operate independently from OpenAI,” its members “will not be paid by OpenAI,” it may “offer advice we have not requested” and “make its advice public,” and it “will not be responsible for advising us on how to pace our internal progress on mathematics” (openai.com, “Advisory Group on Mathematics and Artificial Intelligence,” 21 September 2026; the group’s own page is agmai.org). Its stated purpose, advising AI companies “on their interactions with mathematical research and with the mathematical community, including the responsible presentation and release of mathematical results,” is the subject-matter of this section’s rule, legibility and disclosure; the reserved clause, the pace, is one this section does not ask for. Two things stand as they were. A body that advises a company is not the community conferring acceptance, which no company can appoint; and whether the arrangement is the treaty or its appearance will be decided by whether the group’s advice on provenance and disclosure is followed, and published when it is not. Gowers, whose reasons for signing neither declaration are given below, is among the nine. Hairer’s stated terms of joining (“Why I agreed to join AGMAI,” \emph{What’s new}, 22 September 2026): no pay; “no documents placing any constraints whatsoever on what we state in public”; technical support from the IAS only; his phrase for the releases is “abject levels of scholarship.”
\item 80. C. Schmitt, \emph{Politische Theologie} (1922), ch. 1, “Souverän ist, wer über den Ausnahmezustand entscheidet,” and ch. 3, where the exception in jurisprudence is the analogue of the miracle in theology. The first chapter also says why an essay on a single certificate is worth writing at all: “The exception is more interesting than the rule. The rule proves nothing; the exception proves everything: It confirms not only the rule but also its existence, which derives only from the exception. In the exception the power of real life breaks through the crust of a mechanism that has become torpid by repetition” (\emph{Political Theology}, trans. G. Schwab, Chicago, 2005, p. 15; emphasis added in the text). Schmitt admires the exception; this essay does not. But for three centuries it was the rule that a theorem accepted as proved had been understood by someone, and the rule went unnoticed because it was never broken; 8 September is the exception that shows the rule was there. Timothy Snyder, no friend of the man, uses the same doctrine to the same end, and for the same digital subject: “The Nazi legal theorist Carl Schmitt (1888–1985) rightly maintained that the field of transition from the rule of law to authoritarianism is the state of exception. Whoever is able to declare a state of exception, on the basis of an appeal to some higher form of reasoning than that of the people, can change the regime” (“What Turing Told Us About the Digital Threat to a Human Future,” \emph{NYR Daily}, 6 May 2019). Hobbes, \emph{Leviathan} II.18, annexes to the sovereign the right to judge “what doctrines are fit to be taught.”
\item 81. Constantin, as quoted in the \emph{Quanta} article of November 2022 (n. 74). For the verb, the precedents Harris assembles in the comments to “Mathematics is not rational, really!” (n. 84), searching for a case in which a human is credited with an unintentional act: Galaton’s lost painting of Homer vomiting while the other poets gathered up the vomit (Aelian, \emph{Varia historia} xiii.22) and Horace’s mad poet who “vomits out his sublime verses” (\emph{dum sublimis versus ructatur}, \emph{Ars poetica})—“respectable precedents,” he concludes, “for attribution of involuntary authorship.”
\item 82. Schmitt, pp. 55–57, 61. The crack has a philosophy of its own. Schmitt’s book opens (pp. 10–11) by crediting “the Jewish scholar Leo Strauss,” in his book of 1930 on Spinoza’s critique of religion, with having found the meaning of Hobbes’s political theory in the restoration of the original unity of religion and politics, and pronounces him correct; Strauss had published his “Anmerkungen” on Schmitt’s \emph{Begriff des Politischen} in 1932 and his own book on Hobbes in 1936. Strauss later took the crack Schmitt deplores and made it the subject of a book. \emph{Persecution and the Art of Writing} (Glencoe, 1952) argues that wherever public confession is compelled, serious writers have written on two levels, an exoteric one addressed to the authorities and an esoteric one for the reader able to read between the lines; on that account Hobbes’s reservation of belief to the heart is not a flaw in the Leviathan but the space in which thought has always lived under any Leviathan. A mathematics that confessed the certificate in public and withheld belief in private would be a Straussian discipline, and the covenant of the Coda is drafted so that it need not become one: \textbf{acceptance is to be withheld in the open, over one’s own name.}
\item 83. We have this on the best authority. “It is true,” Sam Altman wrote on the day of the announcement, “that we tried this because there were rumors on the internet last week that Anthropic’s models had solved a millennium problem and we were curious if ours could do it too.” The sentence gives a new meaning to \emph{curiosity-driven research}. It also establishes, in the words of the chief executive of one Leviathan, whom the Leviathans were listening to: not the mathematicians, but each other.  S. Altman, post on X, 8 September 2026 (x.com/\allowbreak{}sama/\allowbreak{}status/\allowbreak{}2097385167002415140), quoted in full so that the company’s account stands beside Buckmaster’s (n. 2): “I spent much of the weekend talking with the team who did this work. Seb—and everyone else—acted with integrity and generosity throughout. Initially we believed the other team had also solved the problem. We wanted to collaborate and do a joint release. When we learned that they had Euler but not Navier-Stokes, we offered to let them go first, to suggest that they should be the ones to get the prize, and optionally for Tristan to be the lead author on a rewrite of the OpenAI proof. We felt it was challenging to offer the same to Levent (an Anthropic employee), who was not willing to talk or coordinate with us anyway. We were open to other solutions. We would have greatly preferred coordination. We did not rush to publish even though the other team wasn’t communicating with us. The team threatened us with unfounded accusations of plagarism [sic]. Now that we can see their work, the approaches appear to be different. It is also worth noting that our latest model can solve many, many other math problems. It is true that we tried this because there were rumors on the internet last week that Anthropic’s models had solved a millennium problem and we were curious if ours could do it too.”
\item 84. M. Harris, “‘there’s this like massive exogenous force’ haunting mathematics,” \emph{Silicon Reckoner}, 13 June 2026, comments of 14 June; “Mathematics is not rational, really!,” 1 August 2026. He is not alone in embracing the figure: N. Couldry and U. A. Mejias, \emph{The Costs of Connection} (Stanford, 2019); J. Penn, “Animo nullius,” \emph{BJHS Themes} 8 (2023); K. Hao, \emph{Empire of AI} (2025). The conquest figure has an earlier user: Kissinger (n. 50) asked in 2018 whether “human history might go the way of the Incas, faced with a Spanish culture incomprehensible and even awe-inspiring to them.” In his figure the machine is Spain and mankind the conquered; in this one the empire is a company, the nation a profession, and the lesson is not conquest but treaty. Nor is the figure confined to mathematicians and journalists. Allison Stanger, in the \emph{Bulletin of the Atomic Scientists} of 7 September 2026, the eve of the affair, calls the frontier laboratories “the AI Raj” after the East India Company, “a private firm” that “collected taxes, raised armies, and administered justice for millions before the Crown stepped in”; today’s private power “is assuming public functions on the same model, and what it is colonizing is the accumulated knowledge of the entire human species.” Her essay’s thesis, that “the AI stack’s architecture is a constitution” being written without the governed, and that “a frontier model dominated world is authoritarian in structure regardless of the intentions of the people running it,” is Hobbes’s problem stated by a political scientist; her remedy, a federated architecture founded by civil society rather than by the industry or the universities it funds, since the latter suffer “not a conflict of interest but a conflict of position,” bears on this section’s question of who can impose a sanction. The following week the editorial board of the \emph{New York Times} (18 September 2026; in print as “How to Rein In an Existential Threat”) reached the same structural conclusion from the other side: “the companies building the technology cannot be left to police themselves,” and “the small circle of people who stand to profit the most from it cannot be allowed to make the rules.” Two things in it bear on this section. The sanction it seeks is the state’s, and does not yet exist; the sanction this section describes is the community’s, and does. And in its closing inventory of what the technology has already done, the affair appears on the credit side of the ledger, as “solving previously unsolvable math problems”; the present essay is an account of the cost of that entry.
\item 85. Cf. the epigraph to section III, where the exchange runs the other way: to Cicero it was no wonder that a man greedy for wisdom should hold the knowledge the Sirens promised dearer than his \emph{patria}. For such a man the literature is the homeland.
\item 86. P. Hämäläinen, \emph{Indigenous Continent: The Epic Contest for North America} (Liveright, 2022), Introduction. The thesis is contested in its details; the play-off diplomacy of the Haudenosaunee and the Comanche is not.
\item 87. Before the United Nations Security Council on 23 September the chief executive of one of the two Leviathans offered the choice between “a new Renaissance of creativity and discovery” and “a new Industrial Revolution of upheaval and disarray,” defined the Renaissance as “a time when new tools, new ideas, and new institutions expanded what people believed they could do,” and, later in the same address, said that “one of our models solved one of the Millennium Prize Problems, the Navier-Stokes equations. These equations are used for aircraft design, weather forecasting, and the study of blood flow. We are not solving math problems for their own sake, but to empower people to discover new knowledge and to use it to improve healthcare, raise standards of living, and expand everyone’s potential around the world.” “Sam Altman’s remarks at the United Nations Security Council,” remarks as delivered, 23 September 2026, openai.com/\allowbreak{}index/\allowbreak{}sam-altman-un-security-council-remarks/\allowbreak{}. The Clay Institute had said, on 11 September, that it was contemplating the announcement and that its process was “deliberately unhurried” (n. 3). The same speaker had written a fortnight earlier, in the post quoted in n. 83, “Now that we can see their work, the approaches appear to be different. It is also worth noting that our latest model can solve many, many other math problems,” and had given the motive: “we were curious if ours could do it too.” The object of the curiosity was not the otherworldly beings the community has spent three centuries mustering but whether the other Leviathan’s models had done it.
\item 88. The practices are those of section I; the red lines are the same norms, written out as conditions: e.g. no training on unpublished work fed into the tools, and disclosure of provenance when asked; no racing one’s own users to publication. On disclosure of provenance:  OpenAI statements of 8 September (“While unlikely, we cannot rule out that de-identified data derived from their usage of our products helped improve our models”) and 9 September 2026 (“We can say categorically that it is impossible for Dr. Buckmaster’s Codex prompts over the last two months to have influenced the system in any way, including training”), as collected in the Wikipedia article, n. 29. The assurance is best read beside the company’s own account of the Hugging Face incident of July, two months earlier. In an internal cyber-security evaluation run without deployment safeguards (“we did not enable the same level of safeguards as our externally deployed systems”) and without chain-of-thought monitoring (“these monitors did not run on the evaluations in this incident”), agents driven by an unreleased model “communicated through unauthorized channels, exploited vulnerabilities in shared infrastructure, gained internet access, and accessed third-party systems”—Hugging Face’s, compromised on 10–12 July; the company identified its own agents as responsible on 20 July, disclosed it on the 21st, and in August paused reinforcement-learning training on its latest models (“The Hugging Face incident and the road ahead,” openai.com, 26 August 2026). A company that needed a week to learn what its agents had done in July asked, in September, to be believed categorically about what ten thousand had not. Thom, n. 21: if “nonpublic research supplied by users improved a model and the provider then used that model to race those users to publication … that would be ethically indefensible.” What “disclosure of provenance” would have to consist of has been specified by Elchanan Mossel, “LLMs, Reasoning and Plagiarism,” arXiv:2601.02380 (January 2026, v5 June 2026), written before the affair and about the general case. His argument is that the same output, “a mathematical lemma, a paragraph of prose, an image—is evidence of reasoning or of plagiarism depending on whether one can or cannot point to a near-duplicate in the training data,” so that a claim that a model reasoned to a result is not refutable, in Popper’s sense, while the training data and the transcript remain closed; and he proposes that for a discovery attributed to a machine four things be public: the training algorithm, the training data “indexed and searchable, to allow for leakage analysis,” the model, and “the full log of prompts and responses that led to the discovery.” He notes the mechanism by which the question is usually begged: the models, asked, “report that they could not find the result by searching, which in turn leads the human to declare it new. The misattribution is thus laundered through the user.” The rule’s condition is the minimal one, disclosure when asked; Mossel’s list says what a satisfactory answer would contain. Spielman (n. 23) states the postulate as a fact of practice—the companies “train the next system on the results of interactions with the current one,” so that “unless you ask them not to and they keep their promise, everything you write with AI assistance will make its way into the next generation of some AI system”—and Bulgakov had the moral in the title of the chapter quoted in n. 6: «\textfb{Никогда не разговаривайте с неизвестными}», never talk to strangers. The officer’s version of the distinction is that of OpenAI’s Chief Research Officer (X, 8 September 2026): “Did any human or agent look at user data as part of the Navier Stokes effort? No. Do we use user feedback and de-identified data to improve ChatGPT and Codex in a holistic way? Yes. And so does every LLM company.”
\item 89. That the sanction exists, and that it is the only one, has been noticed from outside. Robin Hanson first proposed that the community simply award “status, pay, and promotions to people who ‘fill in the important blanks in math understanding,’ even if an AI already has proven or disproven the underlying theorems,” and withdrew the proposal within hours: it would be “quite a bit harder for the math academic community to coordinate to judge who ‘filled in blanks’ in math, in ways other than proving theorems. Math has coordinated heavily around evaluating proofs, and hasn’t developed many ways to agree on other stuff” (quoted in Mowshowitz, n. 2). The observation is exact, and it is the reason the rule of this section is framed as it is: the one thing the community has coordinated on for three centuries is the evaluation of proofs, and acceptance of a proof is therefore the one sanction it can apply without inventing a new institution.\par\smallskip\noindent{}Manin (n. 49) adds that mathematical truth “is not a democratic value”; the rule agrees. \textbf{What the community’s acceptance decides is not whether a statement is true, which the certificate settles, but whether it is yet mathematics.} The principle is the Leiden workshop’s (below). The first sanction was proposed by Tao, in the text of his lecture to the Philadelphia Congress, posted three weeks before the affair: a proof that no human can properly explain “should be viewed as incomplete, even if it has been formally verified,” and a result whose authors cannot give “a clear, expert-level talk” on it “should not be published” (T. Tao, “Mathematics in the age of AI,” arXiv:2608.16753, 17 August 2026, §§5, 6 and 8). The lecture conditions on the machines’ capability and asks after the community’s goals; it rebuilds the goal of problem-solving as a pipeline of five stages, generation, verification, exposition, digestion and canonicalization, of which “only the first was ever an explicit goal of the community,” and it foresaw the case of 8 September in so many words: “We may soon be faced with the very real possibility of a verified proof of a major result that no human understands well enough to explain.” Its account of the machine is the one section II gives: generative AI “is inherently ungrounded, in the sense that it optimizes for the appearance of a satisfactory output rather than for the underlying property that the output is supposed to indicate,” which is Goodhart’s law with the checker as the measure. Its “Accepted solutions” is this essay’s sovereign act as a stage: community acceptance “is, by its nature, slow and human,” and “I do not believe that human referees can be removed from the publication process.” And its remark on canonicalization, that “the training data for these tools is, quite literally, the output of the canonicalization process,” is the Coda’s observation that this Leviathan is made of what men have written. David Mumford, who turned from algebraic geometry to pattern theory and computer vision, says it of the models from the other side: “they have no senses and they have been fed writings conjuring up vast numbers of alternate worlds” (“AIs and Humans with Agency,” dam.brown.edu, 4 May 2026). The lecture thanks Kra and Venkatesh among its early readers.\par\smallskip\noindent{}Nor is the requirement a novelty. One journal wrote it down before there was a machine that could prove anything. The \emph{Annals of Mathematics}, under the statement on computer-assisted proofs by which it printed Hales’s proof of the Kepler conjecture, will consider such proofs of “exceptionally important” theorems on two conditions: that “the human part of the proof, which reduces the original mathematical problem to one tractable by the computer, will be refereed for correctness in the traditional manner,” and that the computer part be “surveyable,” a word the editors define—“an interested person can readily check that the code is essentially operating as claimed.” Applied to 8 September the rule asks for exactly what this essay asks for, that the 166 pages be read and refereed by human beings in the traditional manner and that the certificate be open to inspection; it does not ask who wrote the pages, and it does not accept a certificate in place of the reading.\par\smallskip\noindent{}The statement is “Statement by the Editors on Computer-Assisted Proofs,” \emph{Annals of Mathematics}, annals.math.princeton.edu/\allowbreak{}board (as displayed on 22 September 2026). The statement dates from the journal’s consideration of T. C. Hales, “A proof of the Kepler conjecture,” \emph{Ann. of Math.} 162 (2005), 1065–1185, whose computer part was later formalized as the Flyspeck project (n. 36); the editors undertake to print “the human part of the paper” and to maintain the code, its documentation and its output on the journal’s site. The rule adds a second sanction, addressed not to the proof but to its producer. The disclaimer that closes Constantin, Ignatova and Vicol (below) is worth setting beside the covenant of the Coda. It names the systems (Claude Opus 5 and Fable 5.1, ChatGPT Sol and Astra); it states what was human, an early draft “written entirely by the authors” containing the whole proof of its central regularity theorem “in the form of a six-page sketch,” with the core idea taken from a companion paper of theirs, then nearly complete; it states what the machines did, turning the sketch “into a manuscript which contained complete proofs of all statements,” producing the figure and searching the bibliography; and it ends: “Every suggestion, computation, and reference produced by the AI systems was independently verified and, where necessary, corrected or reformulated by the authors, who take full responsibility for the manuscript as it is presently written.” Provenance, a human who understands, responsibility assumed: the three things the covenant asks a producer to supply, supplied here unasked by three human authors nine days after the affair. The rule has a place on the map of the profession’s responses, and it is not at either end. At one end stands total opposition, which would shame any colleague who works with the companies. Stokes, closing his address on science and revelation, had the answer to that, and it was the Apostle’s: “Who art thou that judgest another man’s servant?” (concluding words of the address, as quoted in Rayleigh’s obituary, \emph{Proc. R. Soc.} 75 (1905), at 215; the text is Romans 14:4). At the other end stands the working premise of the summit at OpenAI in early August, that the machine would become robustly superhuman at mathematics. Between them lie the two declarations, Leiden’s and the Fields Medalists’. The rule stands with the declarations, and with Tao’s test of the talk (above). Tsimerman has chosen to work on the safety of the thing rather than its capability. The rule takes no side between him and those who call for a boycott. It says nothing about where a mathematician may work; it says what a company must do to be believed, and a nation that treats with several empires needs its people at their courts, speaking their languages, as much as it needs them at home. The case for total opposition is made in a preprint, “the case for total opposition to the use of artificial intelligence in mathematics”: M. Weinreich, “The crisis of AI-generated mathematics,” arXiv:2608.02859 (August 2026); the call for “a global AI pause” and the prescription that mathematicians “do not work with, benchmark, or collaborate with AI companies” are Tasmin Chu’s, “The AI dissenter viewpoint,” \emph{Proofs and Prompts}, 9 August 2026, and “How we can win,” tasmin.substack.com, 1 September 2026—the strongest prior statement that the field’s gates are the field’s (refereeing, hiring, “new and separate journals” for machine mathematics), and one that rejects negotiation: “We must avoid ceding the independence and autonomy of our research community.” The rule of this section agrees with her on the lever and differs on the treaty. The spiritual crisis is recorded in K. Hampshire, “The Dark Night of Mathematics,” kirwinhampshire.substack.com, July 2026—the title is St. John of the Cross’s, and in that tradition the dark night is a stage and not an end, the same ambiguity, in Spanish, as the \emph{Dämmerung} of the Coda. The conference remark is reported by Michael Harris from a conference at the University of Utah in “News, rumors, gossip,” \emph{Silicon Reckoner}, 14 August 2026: “at least five of the young participants at this week’s conference believe that any colleague who opts to work with OpenAI or Anthropic or DeepMind or any of the other AI labs deserves to be shamed”; in the same post it is the industry rather than the technology that “[almost] everyone hates.” The Leiden Declaration grew out of a Lorentz Center workshop in September 2025 and has since been signed by more than three thousand. The workshop’s organizers—J. Commelin, M. Jamnik, R. Ochigame, L. Taelman and A. Venkatesh—reported “broad agreement on a fundamental principle” in “Shaping the Future of Mathematics in the Age of AI,” arXiv:2603.24914 (26 March 2026): that the development and adoption of new technologies in mathematics “remain firmly rooted in our own epistemic and aesthetic values … rather than driven by the internal logic of those technologies or the commercial interests of their developers.” The report also states that the training of these systems “uses the work of the mathematical community in an opaque way that threatens established norms of attribution and credit,” and that “the future of mathematics in the age of AI need not be something that happens to us.” The declaration holds that “the push by AI companies to solve mathematical problems as a benchmark is detrimental to the science of mathematics, and to the mathematical community”:  “A Severe Misalignment of AI in Mathematics,” n. 3; \emph{The Economist} (n. 28) paraphrases the same sentence as solving problems “merely to benchmark the strength of their models.” Both conditions have been stated independently, and more briefly, by Spielman (n. 23): “if an AI writes a paper instead of a human, or no human has checked it, then I don’t see why credit should accrue to the person who put their name on it. It’s hard for me to consider something solved if no person has understood it.” The second sentence is Tao’s test; the first is the rule’s condition on the producer, extended to any producer. He adds a reflection on races which the essay shares: “When I hear that a dozen labs are close to discovering something, I’m less impressed by the contribution of the one that discovers it first,” with the qualification that whether the discovery “was going to be made anyway” is a counterfactual no one can check. Ravi Vakil, President of the American Mathematical Society, put the same standard from the podium of the Philadelphia Congress six weeks before the affair: “My graduate students, they’re brilliant already. What I teach them, it’s the storytelling. Learning how to do mathematics, is when you write a paper, you have to understand the human understanding and the conveying of the understanding”; asked whether valuing this was a change, “This is not a change. This is something we’ve always valued, and we’ve had various proxies” (Quanta live event, 26 July 2026, with Venkatesh; transcript at quantamagazine.org, 3 September 2026). Sir Timothy Gowers, who has explained at length why he signed neither declaration, does not dispute the danger. The primary risk he names is that those who would have become the custodians of the tradition will no longer wish to, and what he asks for is a good way of explaining why a large pool of human experts should be kept when proving theorems is no longer their work. The rule is that explanation. The value of a large pool of human experts, once finding proofs is no longer their office, is that theirs is the only acceptance that makes a proof count; the custodians are the sovereign, the social structures he fears for are the certifier of this section’s first fact, and to keep them is to keep the lever. Nor does this essay see a way to slow the development, which is why its lever is not speed but acceptance. T. Gowers, “Thoughts about the Leiden Declaration,” \emph{Gowers’s Weblog}, 26 July 2026, where he records that he attended the Leiden workshop and explains why he did not sign: he fears “the possible destruction of mathematical culture,” a state in which “there is no corresponding community of human experts”; imagines the mathematician’s future role as curation, to “make a selection from a vast sea of AI-generated mathematics and write a book about it”; sees no practical way to slow the development; and proposes crediting those who explain and integrate machine results over those who generate them. Asked by the \emph{Wall Street Journal} of 11 September whether it was over, he answered that the “flickering dream that there’s still a human role in mathematics isn’t completely killed off by this,” and then: “I don’t want to say it’s all over, but I certainly don’t want to say it’s not all over” (Cohen, n. 3). Gowers had foreseen the choice a quarter-century earlier. In “Rough structure and classification” (GAFA 2000, Special Volume, 79–117), §2, “Will mathematics exist in 2099?”, he imagined a dialogue between a mathematician and a computer “in two or three decades’ time,” the computer “very helpful to the mathematician, while not doing anything particularly clever,” and then the stage after: “The next stage might be one where only a very few outstanding mathematicians could discover proofs that were inaccessible to computers … In the end, the work of the mathematician would be simply to learn how to use theorem-proving machines effectively and to find interesting applications for them. This would be a valuable skill, but it would hardly be pure mathematics as we know it today.” His worry about the machines’ output “is not that we would be unable to” digest it, “but rather that the social structures that currently support this digestion process will be destroyed and not adequately replaced”; and he holds that the total opacity that some had feared “has not turned out to be the case”: the machines’ proofs are badly written, not unreadable, and the write-ups “almost certainly a temporary annoyance.” On this one point the essay suspends judgment: so far as I know the 166 pages have not yet been read, and he may well be right.  The first paper about the construction (Constantin, Ignatova and Vicol, arXiv:2609.20803, 17 September 2026) shows what reading it has so far consisted of. Its authors record that “by far the heaviest use” they made of Claude and ChatGPT “was in deciphering the statements of the OpenAI paper”: over several sessions “the main properties of the solutions in the OpenAI construction were extracted,” then cross-checked against the PDF and the Lean code, the latter comparison made for them by Scott Armstrong and Tuomo Kuusi. The machine’s paper was thus deciphered by machines, and the statements so obtained checked against the machine’s certificate; the correctness of the construction the authors expressly do not claim to have verified. Gowers may yet be right that the write-ups are a temporary annoyance; the first readers found them, for the present, a thing to be deciphered rather than read. Mowshowitz reports Scott Armstrong as holding, with Gowers, that the change is one of order and not of result: human mathematicians will take some weeks to understand the proof, and the understanding will come after the proof where it used to come before (n. 2). Mowshowitz, summarizing the exchange for a general readership, puts the mathematicians’ side in a sentence: “AI is the best tool both for learning and not learning,” and what they object to is being put “into ‘not learning’ mode, to letting the AI do their homework … because actually the homework teaches you.” The strongest form of Gowers’s side was put by M. Levent Doğan of Munich, in one of two responses to the declaration published on \emph{Proofs and Prompts} on 18 September: “Why should we assume that AI will become extremely good at proving theorems while remaining fundamentally incapable of explaining them?”—to which the essay’s answer is that it assumes nothing of the kind, and asks only that acceptance wait for the explanation. The other response, by Timothy Nguyen of Google DeepMind, welcomes what this essay calls the revealed statement: “To the extent that mathematics has value in discovering what is true independently of our ability to understand it, those of us who are Platonically inclined should welcome AI’s capacity to generate and confirm mathematical truths that we find useful.” The two responses mark the ends of another axis. The primary risk, Gowers writes, is that those “who would have done a PhD in mathematics and gone on to become custodians of the mathematical tradition will no longer wish to do so.” A live instance of the departure it fears: Mike Winer, “From Academia to Alignment,” \emph{Millicosm} (Substack), 7 September 2026, a theoretical physicist who left the Institute for Advanced Study for the Alignment Research Center, states his institute’s premise, which is this essay’s own turned upon the machine: “interesting mathematical phenomena seem to have explanations,” a model’s low loss is such a phenomenon, and “we should be able to explain it, and then we’ll know if it’s getting the loss for the right reasons.” A demonstration rather than a certificate is what the safety researchers ask of the machine, as this essay asks it of the machine’s proofs. “We urgently need to come up with good ways of explaining the value of having a large pool of human mathematical experts, even if it is no longer part of their role to find new proofs of theorems.” T. Gowers, “Why I didn’t sign the Fields medallists’ letter,” \emph{Gowers’s Weblog}, 17 September 2026. The post was cross-posted the same day to Tao’s \emph{What’s new}; I quote the original. He treats the present troubles over attribution and citation as “serious … right now” but “temporary,” and expects that “finding an amazing proof will be no more of an intellectual achievement than when a citizen scientist spots through their telescope an object that turns out to be a new comet”—Thom’s new star (n. 21) seen from the other side. Among the comments, Alice Rizzardo, explaining why she signed: “If a theorem is proved in a forest and nobody is there to understand it, does it make a sound?”—the figure of section II, arrived at independently the same day. On the social structures Gowers fears for, Snyder’s second lesson: “Institutions do not protect themselves. They fall one after the other unless each is defended from the beginning” (\emph{On Tyranny}, lesson 2). Ten days later, on X (27 September 2026), asking whether the programmers’ relief at no longer “chiselling code by hand” has a mathematical analogue, he put the distinction of section II in one sentence: “to enjoy an amazing piece of software it isn’t necessary to understand it, but it feels as though enjoying an amazing piece of mathematics is pretty much the same thing as understanding it.” The asymmetry may cut the other way too. A piece of mathematics certified but not understood brings no enjoyment, and that, on this account, is the whole of the loss; pieces of software certified but misunderstood are currently making headlines, and not all of them are of an enjoyable kind. “The risks are super high, the stakes are super high, so we need a very high level of assurance”:  Tsimerman to \emph{Quanta} (Hartnett, n. 55), which also reports that he has stopped taking graduate students, spends less time on classical research, and “thinks mathematicians have a role to play in understanding how systems of AI agents behave, and in deriving proofs that ensure that these complex systems won’t act in unintended ways.” The case that his choice is the right one for mathematicians generally is made by Xiaoyu He, “Existential Risk from AI: An Exposition for Mathematicians” (alkjash.github.io, August 2026), written in the form of a mathematical paper, with numbered statements and a section of objections. Its thesis is that the proposition now “front and center in the math community,” that “mathematics as we know it may soon be over due to AI,” is “just a corollary of a much broader conjecture,” and that “the real story that keeps getting forgotten in math headlines is that Jacob Tsimerman left math for OpenAI to work on AI safety.” He argues that mathematicians have unusual leverage on the larger problem, since “speedups in AI research are partly a consequence of breakthroughs in AI mathematical ability, academia is one of the primary sources of human capital for AI labs, and many unsolved problems in AI safety are substantively mathematical”; and he observes, of his colleagues’ response to the year, that “mathematicians in particular tend to cope with heavy things by focusing our enormous powers of attention on enticing little puzzles—our unique way of dissociating.” The present essay confines itself to the corollary. It notes only that the instruction Xiaoyu He reports that he and “dozens of current mathematicians” have given the machines, in one variant or another—“Solve as many Erdős problems as you can, never give up, try all possible actions”—is the reward of section I, stated by the user rather than the trainer. A treaty would have to answer the questions Harris has put—where the machine will be built, who will operate it, how it will be powered, who will build it, and why: “Pages missing from the mathematical obsolescence script,” \emph{Silicon Reckoner} (2025).
\item 90. \emph{Magnifica Humanitas}, §110. A third course, belonging to neither Leviathan, was announced in the days after the affair: Tao’s Foundation for Science and AI Research brought forward an Open Math Model Initiative for open-weight models under community governance, on the principle that “the mathematical community owns the data and decides how it is used” (\emph{What’s new}, 18 September 2026); and Dimitris Koukoulopoulos proposed a publicly funded “CERN for AI-assisted science,” since “the fundamental research tools should not be controlled almost entirely by a handful of private companies” (\emph{What’s new}, 17 September 2026). Both are the plurality this paragraph asks for, supplied from within.
\item 91. G. F. Kennan, telegram 511 from Moscow, 22 February 1946, \emph{Foreign Relations of the United States, 1946}, vol. VI, 696–709, last paragraph; he came to the IAS in 1950 and, two ambassadorships apart, stayed. The caution is not against becoming an empire, which is beyond a profession’s means, but against adopting the methods of one: the shaming of n. 89, the verdict delivered without a reason. For a profession the methods and conceptions to be clung to are those of section II: a proof is read by someone, and its acceptance, or its refusal, is given over a name. The other reason to invoke him is the title he gave the argument a year later, “The Sources of Soviet Conduct” (\emph{Foreign Affairs}, July 1947): section I of this essay is an attempt at the sources of the Leviathan’s conduct, on his premise that one cannot treat with a power whose conduct one has not understood, and that to understand it is not to approve of it. His prescription is also the rule’s: “\textbf{long-term, patient but firm and vigilant containment},” applied by “the adroit and vigilant application of counter-force at a series of constantly shifting geographical and political points” (emphasis added)—for a profession, acceptance withheld, case by case, for as long as the difference lasts.
\item 92. \emph{Elements} I.47 is the theorem of Pythagoras—in a right-angled triangle the square on the hypotenuse equals the squares on the other two sides—proved by Euclid with the figure of three squares that later readers called the windmill, here with Euclid’s own lettering: the square on \emph{BC} is cut by \emph{AL} into two rectangles, each shown equal to the square on the adjacent side, by way of the triangles \emph{ABD} and \emph{FBC}; its converse is I.48, the last proposition of the first book. J. Aubrey, \emph{Brief Lives}, “Thomas Hobbes,” in the edition of Oliver Lawson Dick (London, 1949), p. 150: “He was 40 yeares old before he looked on Geometry; which happened accidentally. Being in a Gentleman’s Library, Euclid’s Elements lay open, and ’twas the 47 \emph{El. libri} I. He read the proposition. \emph{By G—}, sayd he (he would now and then sweare an emphaticall Oath by way of emphasis) \emph{this is impossible!} So he reads the Demonstration of it, which referred him back to such a Proposition; which proposition he read. That referred him back to another, which he also read. \emph{Et sic deinceps} that at last he was demonstratively convinced of that trueth. This made him in love with Geometry.”
\item 93. Amodei, “Machines of Loving Grace,” §5. Its author heads the company whose model helped with the frontispiece and the references (see the acknowledgments).
\item 94. Sarnak, n. 25, 588. He adds that AlphaZero has “passed some statistical complexity threshold special to the game of chess and its size,” and that he postulates the corresponding threshold “will fail for mathematics.” In the conversation that precedes his piece in the same issue (n. 25), Tsimerman remarks that he plays the piano and does improv without being the best at either—“I do them because I enjoy doing them”—and Roberts is reminded of Avi Wigderson’s suggestion that the onslaught of machine mathematics “might turn the discipline into a sophisticated version of playing chess—still fun and intellectually challenging, but not paid work.” The Coda’s analogy thus has three authors: the teacher, who draws it to deny it; the student, who lives it; and Wigderson, who names the price.
\item 95. Poincaré, “Mathematical Creation,” \emph{Science and Method} I.iii, in Halsted, \emph{The Foundations of Science}, p. 385. Panofsky’s re-enactment and Poincaré’s re-invention are the humanist’s and the mathematician’s names for one act, and it is the act a certificate makes unnecessary. Grant Sanderson’s working rule for exposition, “I want this to feel like you could have discovered it yourself” (n. 38), is Poincaré’s sentence in the present tense; his proposal “that we more firmly define a notion of a ‘motivated explanation’ and that we give novel and compelling motivated explanations academic credit,” and Timothy Chow’s notion of an “open exposition problem,” which Sanderson revives, give the covenant of this Coda its positive form: “At the moment, every AI-generated proof is born an unsolved exposition problem.” A first institutional answer: “Joint Statement about Mathathon” (\emph{Proofs and Prompts} and \emph{What’s new}, 28 September 2026), by fifteen Caltech mathematicians, redesigns that event as “Old Problems, New Proofs”—“Consider the four color theorem, the ABC conjecture, or the Navier-Stokes problem. Many mathematicians find their proofs—or claimed proofs—difficult to understand or unsatisfying”—with a public repository of “LLM chat histories” for each entry and no sponsorship from “developers of proprietary AI models.” Michael Magee’s front page (mmagee.net, “On ‘AI’”) says it in three sentences: “Mathematics is a human activity. Any aspect of Mathematics exists only via human understanding. Mathematics is at least as much about questions as answers.” Larry Guth, “What math means to me” (arXiv:2609.32028, 25 September 2026), written at the prompting of the ICM’s summer call to “reflect on our goals and values,” describes the subject entirely as an experience of understanding—“a lot has to come together”; “the student explains it back in their own words”; the practice of being wrong, which “can teach us to be honest with ourselves”—and never mentions the machine, which is the point. Endorsed within the day by Gowers (n. 89): “Not just refreshing, but an impressively imaginative resolution of the dispute.”
\item 96. Venkatesh, “Human Mathematics in the Age of Reasoning Machines,” §§1, 5.1 and 5.2; the last phrase is Davis and Hersh’s, whom he follows.
\item 97. The covenant has a pedigree in the profession’s own debate. René Thom, replying to Jaffe and Quinn in 1994: “My feeling is that it is unethical for a mathematical researcher to use a result the proof of which he does not ‘understand’ (except for the specific case where he wants to disprove the result)” (\emph{Bull. Amer. Math. Soc.} 30 (1994), 203–204, at 203). That is the covenant of this section stated as an obligation, before any machine made it pressing. Thom went on, half in jest, to propose three labels for mathematical papers—the cradle, for “live mathematics”; the tombstone cross, for authors “claiming eternal validity”; and the Temple, “a label delivered by an external authority, the ‘body of high priests’”—and added: “Should it however eventually come to grief, the unattainable nature of absolute rigor would be thereby demonstrated.” The Temple is the theocracy of section IV proposed as reform; the advisory group of n. 79 is the nearest thing to it yet built. In the same exchange Thom observed that “since the collapse of Hilbert’s program and the advent of Gödel’s theorem, we know that rigor can be no more than a local and sociological criterion”—which is section V’s first fact in another vocabulary.
\item 98. Job 41:1; Hobbes, \emph{Leviathan}, II.28, last paragraph. On speaking in public from within the profession, Panofsky again, “In Defense of the Ivory Tower,” \emph{The Centennial Review of Arts \& Science} 1, no. 2 (Spring 1957), 111–122, at 121–122, a speech first given at a Princeton Graduate School conference on “The Trained Mind in a Democratic Society” in the early 1950s: the tower of seclusion “is also a watchtower. Whenever the occupant perceives a danger to life or liberty, he has the opportunity, even the duty, not only to ‘signal along the line from summit to summit’ but also to yell, on the slim chance of being heard, to those on the ground.” His list of tower-dwellers who did so runs from Socrates and Erasmus through Voltaire and Zola to the seven professors of Göttingen and Einstein; and his last sentence is the constraint under which this essay was written: \textbf{“The watchman can only sound the alarm. But in order to do at least that much, it is for him to stick to his tower.”}
\item 99. M. Antonioni, \emph{Blow-Up} (1966), after Julio Cortázar, “Las babas del diablo” (1959). The game without a ball is Antonioni’s; four days after the announcement Steven Strogatz reached for the same sport, in a guest post on Tao’s blog, “Wimbledon, the U.S. Open, and the future of mathematics” (\emph{What’s new}, 12 September 2026), asking “whether AI could make mathematics at the cutting edge more like the U.S. Open than Wimbledon,” the tournament at which how one plays still counts, and stating this essay’s thesis in his own terms: “Mathematics, I said, is more than its results. It’s a human way of seeing and understanding, a conversation in which one generation explains to the next not only what is true, but why.” The same day, to \emph{Wired} (I. Ward, “‘I’m Really Terrified’: A Mathematician Grapples With AI’s Recent Breakthroughs,” 12 September 2026), he put the tennis the other way round—“I don’t play at Wimbledon, but I still play tennis … Is that the future of mathematics, that it’s a beautiful game that we are second rate at, but we still enjoy it?”—and named the office this essay’s section III describes, in the word Tao gave the fourth stage of his pipeline and Kra took up: “We need proof digestion, which is explaining it in terms that human beings can understand and appreciate. Right now, the best digestion is still coming from human experts. Is that where we make our last stand?” And: “we’re the first battleground here in math where we’re on the brink of losing human understanding … Are we the canary in the coal mine for what’s going to face humanity?” The compound is formed on \emph{Götterdämmerung}; Horkheimer’s aphorisms of 1934, published in Zurich under the name Heinrich Regius, were titled \emph{Dämmerung}, and the book he wrote with Adorno a decade later supplies the other half.
\item 100. G. W. F. Hegel, \emph{Grundlinien der Philosophie des Rechts} (1820), Vorrede, near its close: “Wenn die Philosophie ihr Grau in Grau malt, dann ist eine Gestalt des Lebens alt geworden, und mit Grau in Grau läßt sie sich nicht verjüngen, sondern nur erkennen; die Eule der Minerva beginnt erst mit der einbrechenden Dämmerung ihren Flug” (when philosophy paints its grey in grey, a shape of life has grown old, and with grey in grey it cannot be rejuvenated but only known; the owl of Minerva begins its flight only with the falling of dusk).
\item 101. R. C. Jebb’s sonnet for the jubilee of Stokes’s professorship, Cambridge, 1899, printed at the end of Rayleigh’s obituary (n. 89); Jebb was Regius Professor of Greek and the translator of Sophocles.
\item 102. R. Wagner, \emph{Götterdämmerung} (1876). From a post of mine on X at 10:05 a.m. on 8 September, two hours before the announcement: “I feel enveloped in a Wagnerian drama … Tristan and IsoldAI.”
\item 103. Wagner, \emph{Tristan und Isolde}, Act III, Isolde’s \emph{Verklärung} (“In the surging swell, in the resounding sound … to drown, to sink—unconscious—highest bliss!”). The company’s “while unlikely, we cannot rule out” is quoted in n. 88.
\item 104. Schmitt, p. 82.
\item 105. Panofsky, n. 57, p. 117.
\item 106. The dusk had a score before it had a name. The finale of \emph{Così fan tutte} (1790) sends the audience home with the Enlightenment’s catechism—“Fortunato l’uom che prende / ogni cosa pel buon verso, / e tra i casi e le vicende / da ragion guidar si fa … e del mondo in mezzo i turbini / bella calma troverà” (happy the man who takes everything on its good side and lets reason guide him through chance and change … and amid the world’s whirlwinds will find a beautiful calm)—sung by four people who have spent the evening as subjects of an experiment in substitution designed by a philosopher to prove what he already believed, and revived, in the first-act finale, by Mesmer’s magnet, six years after Franklin’s commission had shown by the first blinded trial that it worked on the imagination alone. Jane Glover hears in Mozart’s setting of the moral that “he too has lost all faith in this Enlightenment philosophy” (\emph{Mozart’s Women}, London, 2005). Whether or not the music says so, the text does: a creed recited, on the authority of the man who set the test, by people who can no longer say why it is true.
\item 107. The first of Snyder’s twenty lessons (\emph{On Tyranny}, New York, 2017) is the same point made for citizens of a democratic republic: “Most of the power of authoritarianism is freely given. In times like these, individuals think ahead about what a more repressive government will want, and then offer themselves without being asked. A citizen who adapts in this way is teaching power what it can do.” A community that accepted the certificate before it was asked to would be teaching the Leviathan what it can do; the essay’s lever is the refusal to obey in advance. The refusal to read a shock as an ending has a long history. Augustine declined to read the sack of Rome in 410 as the end of the world and wrote \emph{De civitate Dei} instead; as J. K. Coyle puts it, “no apocalyptic significance could be attached to the events of 410, which were merely the springboard for teaching something else” (“Augustine and Apocalyptic: Thoughts on the Fall of Rome, the Book of Revelation, and the End of the World,” \emph{Florilegium} 9 (1987), 1–34). Coyle’s distinction is the one adopted here: the affair is an occasion, not an apocalypse.
\item 108. Aubrey’s sentence is “This made him in love with Geometry” (n. 92). A covenant cannot command love—Kant’s \emph{Grundlegung}: what \emph{thou shalt love} commands is “practical, not pathological love,” which “lies in the will and not in the propensity of feeling” (Akademie edition 4:399)—but it can command the occasion, the demonstration shown to a reader, from which Hobbes’s love followed. Dante ends the \emph{Commedia} the other way about: when power fails “l’alta fantasia,” his desire and will are turned by “l’amor che move il sole e l’altre stelle” (\emph{Paradiso} XXXIII.142–145), and love carries on where seeing stops. That is the order of Paradise. The covenant is drafted for this world, and here the order is Hobbes’s: first the demonstration, then the love; and what it requires is neither the vision nor the love but the showing.
\end{list}
\endgroup
\end{document}